\documentclass[10pt,reqno]{amsart}
\usepackage{enumerate}
\usepackage{amssymb}
\usepackage{bm}
\usepackage{xfrac}
\usepackage{graphicx}
\usepackage[centertags]{amsmath}
\usepackage{amsfonts}
\usepackage{amsthm}
\usepackage{a4}
\usepackage{latexsym}
\usepackage{amsmath}
\usepackage{amsfonts}
\usepackage{amssymb,amsbsy}
\usepackage{amscd}
\usepackage{amsthm}
\usepackage{layout}
\usepackage{color}
\usepackage{mathrsfs}
\usepackage{cite}
\usepackage{epstopdf}
\usepackage[font=footnotesize,labelfont=bf]{caption}

\usepackage[noabbrev]{cleveref}
\usepackage[normalem]{ulem}
\usepackage{comment}

\theoremstyle{definition}
\newtheorem{pr}{Problem}
\newtheorem{thm}{Theorem}[section]
\newtheorem{cor}[thm]{Corollary}
\newtheorem{prop}[thm]{Proposition}
\newtheorem{lem}[thm]{Lemma}

\newtheorem{quest}[thm]{Question}

\newtheorem*{thm*}{Theorem}
\newtheorem*{notn*}{Notation}

\theoremstyle{definition}
\newtheorem{dfn}[thm]{Definition}

\theoremstyle{definition}
\newtheorem{rem}[thm]{Remark}

\numberwithin{equation}{section}

\newcommand{\N}{\mathbb{N}}
\newcommand{\T}{\mathcal{T}}
\newcommand{\tT}{{\widetilde{\T}}}
\renewcommand{\tt}{{\tilde{t}}}
\newcommand{\tw}{{\tilde{w}}}

\newcommand{\tW}{{\widetilde{\W}}}
\newcommand{\M}{\mathcal{M}}

\newcommand{\W}{\mathcal{W}}
\renewcommand{\S}{\mathcal{S}}
\newcommand{\barf}{\bar{f}}
\newcommand{\barg}{\bar{g}}
\newcommand{\Q}{\mathbb{Q}}

\newcommand{\e}{\varepsilon}

\newcommand{\supp}{\textnormal{supp}}
\newcommand{\range}{\textnormal{range}}

\renewcommand{\-}{\textnormal{-}}

\newcommand{\al}{\alpha}

\newcommand{\compsep}{\makebox[-1.6pt][s]{\raisebox{-0.03ex}{$\not$}}\makebox[0.6pt][s]{\raisebox{0.2ex}{$\Rightarrow$}}\makebox[9pt][s]{\raisebox{-0.4ex}{\Tiny$\hookleftarrow$}}\;}

\DeclareMathOperator{\minsupp}{minsupp}
\DeclareMathOperator{\maxsupp}{maxsupp}

\DeclareMathOperator{\ran}{range}

\allowdisplaybreaks

\begin{document}

\title[The complete separation of the two finer asymptotic $\ell_p$
structures]{The complete separation of the two finer asymptotic  $\boldsymbol{\ell_{p}}$  structures for
$\boldsymbol{1\le  p<\infty}$}

\author[S. A. Argyros]{S. A. Argyros}
\address{National Technical University of Athens, Faculty of Applied Sciences,
Department of Mathematics, Zografou Campus, 157 80, Athens, Greece}
\email{sargyros@math.ntua.gr}

\author[A. Georgiou]{A. Georgiou}
\address{National Technical University of Athens, Faculty of Applied Sciences,
Department of Mathematics, Zografou Campus, 157 80, Athens, Greece}
\email{ale.grgu@gmail.com}

\author[A. Manoussakis]{A. Manoussakis}
\address{ Technical University of Crete, Department of Environmental
  Engineering , 73100, Chania, Greece}
\email{amanousakis@isc.tuc.gr}

\author[P. Motakis]{P. Motakis}
\address{Department of Mathematics and Statistics, York University, 4700 Keele Street, Toronto, Ontario, M3J 1P3, Canada}
\email{pmotakis@yorku.ca}

\thanks{The fourth author was supported by NSERC Grant RGPIN-2021-03639.}

\begin{abstract}
For $1\le p <\infty$, we present a reflexive Banach space $\mathfrak{X}^{(p)}_{\text{awi}}$, with an unconditional basis, that admits $\ell_p$ as a unique asymptotic model and does not contain any Asymptotic $\ell_p$ subspaces. D. Freeman, E. Odell, B. Sari and B. Zheng have shown that whenever a Banach space not containing $\ell_1$, in particular a reflexive Banach space, admits $c_0$ as a unique asymptotic model then it is Asymptotic $c_0$. These results provide a complete answer to a problem posed by L. Halbeisen and E. Odell and also complete a line of inquiry  of the relation between specific asymptotic structures in Banach spaces, initiated in a previous paper by the first and fourth authors. For the definition of $\mathfrak{X}^{(p)}_{\text{awi}}$ we use saturation with asymptotically weakly incomparable constraints, a new method for defining a norm that remains small on a well-founded tree of vectors which penetrates any infinite dimensional closed subspace.
\end{abstract}

\maketitle
\setcounter{tocdepth}{1}
\tableofcontents

\section{Introduction}

The purpose of this article is to provide an answer to the following problem of Halbeisen and Odell from \cite{HOdell} and is, in particular, the last step towards the complete separation of a list of asymptotic structures from \cite{AM}. Given a Banach space $X$, let $\mathscr{F}_0(X)$ denote the family of normalized weakly null sequences in $X$ and $\mathscr{F}_b(X)$ denote the family of normalized block sequences of a fixed basis, if $X$ has one.

\begin{pr}
Let $X$ be a Banach space that admits a unique asymptotic model with respect to $\mathscr{F}_0(X)$, or with respect to $\mathscr{F}_b(X)$ if $X$ has a basis. Does $X$ contain an Asymptotic $\ell_p$, $1\le p<\infty$, or an Asymptotic $c_0$ subspace?
\end{pr}

The following definition from \cite{AM} provides a more general setting in which we will describe this problem, as well as other previous separation results. A property of a Banach space is called hereditary if it is inherited by all of its closed and infinite dimensional subspaces.

\begin{dfn}
Let (P) and (Q) be two hereditary properties of Banach spaces and assume that (P) implies (Q).
\begin{itemize}
\item[(i)] If (Q)$\not\Rightarrow$(P), i.e., there exists a Banach space satisfying (Q) and failing (P), then we say that (P) is separated from (Q).
\item[(ii)] If there exists a Banach space satisfying (Q) and whose every infinite dimensional closed subspace fails (P), then we say that (P) is completely separated from (Q) and write (Q) \compsep (P).
\end{itemize}
\end{dfn}

We consider properties that are classified into the following three categories: the sequential asymptotic properties, the array asymptotic properties and the global asymptotic properties.

Sequential asymptotic properties are related to the notion of a spreading model from \cite{BS}, which describes the asymptotic behaviour of a sequence in a Banach space. We say that a Banach space admits a unique spreading model with respect to some family of normalized sequences $\mathscr{F}$, if whenever two sequences from $\mathscr{F}$ generate spreading models, then those must be equivalent. If this equivalence happens with some uniform constant, then we say that the space admits a uniformly unique spreading model.

The category of array asymptotic structures concerns the asymptotic behaviour of arrays of sequences $(x^i_j)_j$, $i\in\N$, in a Banach space. Notions that describe this behaviour are those of asymptotic models from \cite{HOdell} and joint spreading models from \cite{AGM}. We define the uniqueness of asymptotic models and the uniform uniqueness of joint spreading models in a similar manner to the uniqueness and uniform uniqueness of spreading models, respectively. Although asymptotic models and joint spreading models are not identical notions, they are strongly related. As B. Sari pointed out, a Banach space $X$ admits a uniformly unique joint spreading model with respect to $\mathscr{F}_b(X)$ or $\mathscr{F}_0(X)$ if and only if it admits a unique asymptotic model with respect to $\mathscr{F}_b(X)$ or $\mathscr{F}_0(X)$, respectively (see, e.g.,\cite[Remark 4.21]{AGLM} or \cite[Proposition 3.12]{AM}). Notably, the property that a Banach space $X$ with a basis admits some $\ell_p$ as a uniformly unique joint spreading model with respect to $\mathscr{F}_b(X)$ can be described by the following statement. The case where this happens with respect to $\mathscr{F}_0(X)$ is given by an easy modification.

\begin{prop}[Lemma \ref{ell_p as uniformly unique joint spreading equivalent form}]
Let $1\le p\le \infty$. A Banach space $X$ with a basis admits $\ell_p$ (or $c_0$ for $p=\infty$) as a uniformly unique joint spreading model with respect to $\mathscr{F}_b(X)$ if and only if there exist constants $c,C>0$ such that for every $\ell\in\N$, any choice of successive families $(F_j)_j$ of normalized blocks in $X$ with $\#F_j=\ell$, there is an infinite subset of the naturals $M=\{m_1<m_2<\ldots\}$  such that for any choice of $x_j\in F_j$, $j\in M$, every $G\subset M$ with $m_k\le G$ and $\#G\le k$, for $k\in\N$, and any choice of scalars $a_j$, $j\in G$, we have
\[
c\|(a_j)_{j\in G}\|_p \le \| \sum_{j\in G}a_jx_j \| \le C\| (a_j)_{j\in G} \|_p.
\]
\end{prop} 
Even though this property is very close to the weaker one that $X$ admits $\ell_p$ or $c_0$ as a uniformly unique spreading model, it was shown in \cite{AM} that these two properties are in fact completely separated for all $1\le p\le \infty$. 

Finally, global asymptotic properties describe the behaviour of finite block sequences that are chosen sufficiently far apart in a space with a basis. We recall the following definition from \cite{MT}.

\begin{dfn}
Let $X$ be a Banach space with a basis $(e_i)_i$ and $1\le p\le \infty$. We say that the basis $(e_i)_i$ of $X$ is asymptotic $\ell_p$ {(asymptotic $c_0$ when $p=\infty$)} if there exist positive constants $D_1$ and $D_2$ such that for all $n\in\N$, there exists $N(n)\in\N$ with the property that whenever $N(n)\le x_1< \cdots< x_n$ are vectors in $X$, then
\[
\frac{1}{D_1}(\sum_{i=1}^n\|x_i\|^p)^{\frac{1}{p}}\le \|\sum_{i=1}^nx_i\|\le D_2 (\sum_{i=1}^n\|x_i\|^p)^{\frac{1}{p}},
\]
{where for $p=\infty$ the above inequality concerns the $\|\cdot\|_\infty$.} Specifically, we say that $(e_i)_i$ is $D$-asymptotic $\ell_p$ {(D-asymptotic $c_0$ when $p=\infty$)} for $D=D_1D_2$.
\end{dfn}

This definition is given with respect to a fixed basis of the space. The coordinate-free notion of Asymptotic $\ell_p$ and $c_0$ spaces was introduced in \cite{MMT}, generalizing the aforementioned one to spaces with or without a basis (note the difference between the terms asymptotic $\ell_p$ and Asymptotic $\ell_p$). Moreover, this property is hereditary and any Asymptotic $\ell_p$ (or $c_0$) space is asymptotic $\ell_p$ (resp. $c_0$) saturated.

Given a Banach space $X$ with a basis, we focus on the following properties, where $1\le p\le \infty$ and whenever $p=\infty$, then $\ell_p$ should be replaced with $c_0$.

\begin{itemize}

\item[(a)$_p$] The space $X$ is Asymptotic $\ell_p$.

\item[(b)$_p$] The space $X$ admits $\ell_p$ as a uniformly unique joint spreading model (or equivalently, a unique asymptotic model, as mentioned above) with respect to $\mathscr{F}_b(X)$.

\item[(c)$_p$] The space $X$ admits $\ell_p$ as a uniformly unique spreading model with respect to $\mathscr{F}_b(X)$.

\item[(d)$_p$] The space $X$ admits $\ell_p$ as a unique spreading model with respect to $\mathscr{F}_b(X)$.

\end{itemize}

Note that it is fairly straightforward to see that the following implications hold for all $1\le p\le \infty$: (a)$_p\Rightarrow$(b)$_p\Rightarrow$(c)$_p\Rightarrow$(d)$_p$. It is also easy to see that (d)$_p\not\Rightarrow$(c)$_p$ for all $1\le p<\infty$. In \cite{BLMS} it was shown that (c)$_p\not\Rightarrow$(b)$_p$ for all $1\le p\le \infty$ and that (b)$_p\not\Rightarrow$(a)$_p$ for all $1<p<\infty$. The latter was also shown in \cite{AGM}, as well as that (b)$_1\not\Rightarrow$(a)$_1$ along with an even stronger result, namely, the existence of a Banach space with a basis satisfying (b)$_1$ and such that any infinite subsequence of its basis generates a non-Asymptotic $\ell_1$ subspace. However, it was proved in \cite{AOST} that (d)$_\infty\Leftrightarrow$(c)$_\infty$ and a remarkable result from \cite{FOSZ} states that (b)$_\infty\Leftrightarrow$(a)$_\infty$ for Banach spaces not containing $\ell_1$. Towards the complete separation of these properties, it was shown in \cite{AM} that (c)$_p$ \compsep (b)$_p$ for all $1\le p\le \infty$ and that (d)$_p$ \compsep (c)$_p$ for all $1\le p <\infty$. Hence, the only remaining open question was whether (b)$_p$ \compsep (a)$_p$ for $1\le p<\infty$. We prove this in the affirmative and, in particular, we show the following.

\begin{thm}
\label{main theorem intro}
For $1\le p<\infty$, there exists a reflexive Banach space $\mathfrak{X}_{\text{awi}}^{(p)}$ with an unconditional basis that admits $\ell_p$ as a uniformly unique joint spreading model with respect to $\mathscr{F}_b(\mathfrak{X}_{\text{awi}}^{(p)})$ and contains no Asymptotic $\ell_p$ subspaces.
\end{thm} 

To construct these spaces, we use a saturation method with asymptotically weakly incomparable constraints. This method, initialized in \cite{AGM}, employs a tree structure, penetrating every subspace of $\mathfrak{X}_{\text{awi}}^{(p)}$, that admits segments with norm strictly less than the $\ell_p$-norm. Thus, we are able to prove that no subspace of $\mathfrak{X}_{\text{awi}}^{(p)}$ is an Asymptotic $\ell_p$ space. This saturation method is different from the method of saturation with increasing weights from \cite{AM}, used to define spaces with no subspaces admitting a unique asymptotic model. It does not seem possible to use the method of increasing weights to construct a space with a unique asymptotic model, i.e., it is not appropriate for showing (b)$_p\;\compsep$\!(a)$_p$. On the other hand, the method of asymptotically weakly incomparable constraints yields spaces with a unique asymptotic model, and thus it cannot be used to show (c)$_p\;\compsep$\!(b)$_p$. {This method will be discussed in detail in Part 1.}

{In the case of $1<p<\infty$, it is possible to obtain a stronger result. Namely, for every countable ordinal $\xi$, the space separating the two asymptotic properties additionally satisfies the property that every block subspace contains an $\ell_1$-tree of order $\omega^\xi$. This is achieved using the attractors method, which was first introduced in \cite{AAT} and later also used in \cite{AMP}. The precise statement of this result is the following.}

\begin{thm}[\cite{AGMM}]
For every $1<p<\infty$ and every infinite countable ordinal $\xi$, there exists a hereditarily indecomposable reflexive Banach space $\mathfrak{X}^{(p)}_\xi$ that admits $\ell_p$ as a uniformly unique joint spreading model with respect to the family of normalized block sequences and whose every subspace contains an $\ell_1$-block tree of order $\omega^\xi$.
\end{thm}
However, in the case of $\ell_1$, we are not able to construct a space whose every subspace contains a well-founded tree which is either $\ell_p$ for some $1<p<\infty$ or $c_0$. This case is more delicate, since as we mentioned, the two properties are in fact equivalent in its dual problem for spaces not containing $\ell_1$. 

The paper is organized as follows: In Section 2, we recall the notions of Schreier families and special convex combinations, and prove some of their basic properties, while Section 3 contains the precise definitions of the aforementioned asymptotic strucures. In Section 4, we recall certain combinatorial results concerning measures on countably branching well-founded trees from \cite{AGM}, which are a key ingredient in the proof that $\mathfrak{X}^{(p)}_{\text{awi}}$ admits $\ell_p$ as a unique asymptotic model for $1\le p <\infty$. We then split the remainder of the paper into two main parts, each dedicated to the definition and properties of $\mathfrak{X}_{\text{awi}}^{(1)}$ and $\mathfrak{X}_{\text{awi}}^{(p)}$ for $p=2$, respectively. The construction of $\mathfrak{X}_{\text{awi}}^{(p)}$ for $1<p<\infty$ and $p\neq 2$ follows as an easy modification of our construction and is omitted. Each of these parts contains an introduction in which we describe the main key points of each construction. Finally, we include two appendices containing variants of the basic inequality, which has been used repeatedly in the past in several related constructions (see, e.g., \cite{AAT}, \cite{AM}, \cite{AMP} and \cite{DM}).

We would like to thank the anonymous referee for providing many helpful suggestions on how to improve the content of our paper.

	\section{Preliminaries}
	
	In this section we recall some necessary definitions, namely,
        the Schreier families $(\S_n)_n$ \cite{AA} and the
        corresponding repeated averages $\{a(n,L):n\in\N,\; L\in[\N]^\infty\}$
        \cite{AMT} which we call $n$-averages, as well as the
        notion of special convex combinations. For a more thorough
        discussion of the above, we refer the reader to \cite{AT}. We begin with some useful notation.

	\begin{notn*}
	By $\N=\{1,2,\ldots\}$ we denote the set of all positive integers. We will use capital letters such as $L,M,N,\ldots$ (resp., lower case letters such as $s,t,u,\ldots$) to denote infinite subsets (resp., finite subsets) of $\N$. For every infinite subset $L$ of $\N$, the notation $[L]^\infty$ (resp., $[L]^{<\infty}$) stands for the set of all infinite (resp., finite) subsets of $L$. For every $s\in[\N]^{<\infty}$, by $|s|$ we denote the cardinality of $s$. For $L\in[\N]^\infty$ and $k\in\N$, $[L]^k$ (resp., $[L]^{\le k}$) is the set of all $s\in[L]^{<\infty}$ with $|s|=k$ (resp., $|s|\le k$). For every $s,t\in[\N]^{<\infty}$, we write $s<t$ if either at least one of them is the empty set, or $\max s<\min t$. Also for  $\emptyset\neq s\in[\N]^{< \infty}$ and $n\in\N$ we write $n<s$ if $n<\min s$.
	
	We shall identify strictly increasing sequences in $\N$ with their corresponding range, i.e., we view every strictly increasing sequence in $\N$ as a subset of $\N$ and conversely every subset of $\N$ as the sequence resulting from the increasing order of its elements. Thus, for an infinite subset $L=\{l_1<l_2<\ldots\}$ of $\N$ and $i\in\N$, we set $L(i)=l_i$ and similarly, for a finite subset $s=\{n_1<\ldots<n_k\}$ of $\N$ and for $1\le i\le k$, we set $s(i)=n_i$.
\end{notn*}

	Finally, throughout the paper, we follow \cite{LT} (see also \cite{AK}) for standard notation and terminology concerning Banach space theory. 	{For $x\in c_{00}(\N)$, we denote $\supp(x)=\{n\in\N:x(n)\neq 0\}$ and by $\ran(x)$ the minimum interval of $\N$ containing $\supp(x)$. Moreover, for $x,y\in c_{00}(\N)$, we write $x<y$ to denote that $\maxsupp(x)<\minsupp(y)$.}

      \addtocontents{toc}{\protect\setcounter{tocdepth}{1}}  
        \subsection{Schreier Families}\label{comb}

For a family $\mathcal{M}$ and a sequence $(E_i)_{i=1}^k$ of finite subsets of $\N$, we say that $(E_i)_{i=1}^k$ is $\mathcal{M}$-admissible if there is $\{m_1,\ldots,m_k\}\in\mathcal{M}$ such that $m_1\le E_1<m_2\le E_2<\cdots<m_k\le E_k$. Moreover, a sequence $(x_i)_{i=1}^k$ in $c_{00}(\N)$ is called $\mathcal{M}$-admissible if $(\supp (x_i))_{i=1}^k$ is $\mathcal{M}$-admissible. In the case where $\M$ is a spreading family (i.e.,  whenever $E=\{m_1,\ldots,m_k\}\in\M$ and $F=\{n_1<\ldots<n_k\}$ satisfy $m_i\le n_i$, $i=1,\ldots,k$, then $F\in\M$), a sequence $(E_i)_{i=1}^k$ is $\M$-admissible if $\{\min E_i:i=1,\ldots,k\}\in\M$, and thus a sequence of vectors $(x_i)_{i=1}^k$ in $c_{00}(\N)$ is $\M$-admissible if $\{\min\supp (x_i):i=1,\ldots,k\}\in\M$.

For $\mathcal{M}$, $\mathcal{N}$ families of finite subsets of $\N$, we define the convolution of $\mathcal{M}$ and $\mathcal{N}$ as follows:
\begin{align*}
\mathcal{M}*\mathcal{N}=\Big\{ & E\subset\N:\text{ there exists an }\mathcal{M}\text{-admissible finite sequence }\\&(E_i)_{i=1}^k\text{ in }\mathcal{N}\text{ such that }E=\cup_{i=1}^kE_i\Big\}\cup\big\{\emptyset\big\}.
\end{align*}

The Schreier families $(\S_n)_{n\in\N}$ are defined inductively as follows:
\[
\S_0=\big\{ \{k\}:k\in\N \big\}\cup\big\{\emptyset\big\}\quad\text{and}\quad \S_1=\big\{E\subset\N:\# E\le \min E\big\}\cup\{\emptyset\}
\]
and if $\S_n$, for some $n\in \N$, has been defined, then
\begin{align*}
\S_{n+1}=\S_1*\S_n=\Big\{ &E\subset\N : E=\cup_{i=1}^kE_i\text{ where } E_1<\ldots<E_k\in\S_n\\&\text{and }k\le\min E_1\Big\}\cup\{\emptyset\}.
\end{align*}
It follows easily by induction that for every $n,m\in\N$,
\[
\S_n*\S_{m}=\S_{n+m}.
\]
Furthermore, for each $n\in\N$, the family $\S_n$ is regular. This means that {it includes the singletons,} it is hereditary, i.e., if $E\in\S_n$ and $F\subset E$, then $F\in\S_n$, it is spreading and finally it is compact, identified as a subset of $\{0,1\}^\N$.

For each $n\in\N$, we also define the regular family
\[
\mathcal{A}_n= \big\{ E\subset\N : \#E\le n \big\}.
\]
Then, for $n,m\in\N$, we are interested in the family $\S_n*\mathcal{A}_m$, that is, the family of all subsets of $\N$ of the form $E=\cup^k_{i=1}E_i$ where $E_1<\ldots<E_k$,  $\#E_i\le m $ for $i=1,\ldots,k$ and $\{\min E_i:1\le i \le k\}\in \S_n$. In fact, any such $E$ is the union of at most $m$ sets in $\S_n$, and moreover, if $m\le E$, then $E\in \S_{n+1}$, as we show next.

\begin{lem}\label{lemma sn conv am}
For every $n,m\in \N$, 
\begin{itemize}
\item[(i)] $\S_n*\mathcal{A}_m\subset \mathcal{A}_m*\S_n$ and
\item[(ii)] if $E\in \S_n*\mathcal{A}_m$ with $m\le E$, then $E\in \S_{n+1}$.
\end{itemize}
\end{lem}

\begin{rem}\label{remark sn conv am}
Let $k,m\in\N$ and $F$ be a subset of $\N$ with $\#F\le km$ and $k\le F$. Set $d=\max\{1,\lfloor \#F/m \rfloor\}$ and define $F_j=\{F(n):n=(j-1)d+1,\ldots,jd\}$ for each $j=1,\ldots,m-1$ and $F_m=F\setminus\cup_{j=1}^{m-1}F_j$. Then, it is immediate to check that $F_j\in\S_1$ for every $i=1,\ldots,m$. 
\end{rem}

\begin{proof}[Proof of Lemma \ref{lemma sn conv am}]
Fix $n,m\in\N$. We prove (i) by induction on $n\in\N$. For $n=1$, let $E\in \S_1*\mathcal{A}_m$, that is, $E=\cup_{i=1}^kE_i$ with $k\le E_1<\ldots<E_k$ and $\#E_i\le m$ for every $i=1,\ldots,k$. Since $\#E\le km$, Remark \ref{remark sn conv am} yields a partition $E=\cup_{j=1}^m F_j$ with $F_j\in \S_1$ for every $j=1,\ldots,m$ and hence $E\in \mathcal{A}_m*\S_1$. 

Suppose that (i) holds for some $n\in\N$ and let $E\in \S_{n+1}*\mathcal{A}_m$. Then $E=\cup_{i=1}^kE_i$ for an $\S_{n+1}$-admissible sequence $(E_i)_{i=1}^k$ with $\#E_i\le m$ for every $i=1,\ldots,m$. Hence $\{\min E_i:i=1,\ldots,k\}=\cup_{j=1}^lF_j$, where $F_j\in \S_n$ for every $j=1,\ldots,l$ and $l\le F_1<\cdots<F_l$. Define, for each $j=1,\ldots,l$, \[G_j=\cup\{E_i:i=1,\ldots,k\text{ and } \min E_i\in F_j\},\]  and note that $G_j\in \S_{n}*\mathcal{A}_m$ since $F_j\in \S_n$. Hence, for every $j=1,\ldots,l$, the inductive hypothesis implies that $G_j\in \mathcal{A}_m*\S_n$, that is, $G_j=\cup_{i=1}^{m_j}G^j_i$ with $m_j\le m$ and $G^j_i\in \S_n$ for all $i=1,\ldots,m_j$. Define \[H=\{\min G^j_i:j=1,\ldots,l,\text{ and }i=1,\ldots,m_j\}.\] {Observe that $H\in\mathcal{S}_1*\mathcal{A}_m$} and apply Remark \ref{remark sn conv am} to obtain a partition $H=\cup_{q=1}^mH_q$ where $H_q\in \S_1$ for every $q=1,\ldots,m$. Finally, define \[\Delta_q=\cup\{G^j_i:j=1,\ldots,l,\;i=1,\ldots,m_j\text{ and } \min G^j_i\in H_q\},\]
for each $q=1,\ldots,m$, and observe that $E=\cup_{q=1}^m\Delta_q$ and that $\Delta_q\in \S_1*S_{n}=S_{n+1}$ since $H_q\in \S_1$ and $G^j_i\in \S_{n}$. Thus, we conclude that $E\in \mathcal{A}_m*\S_{n+1}$.

Finally, note that (ii) is an immediate consequence of (i). 
\end{proof}

\subsection{Repeated Averages} The notion of repeated averages was first defined in \cite{AMT}. The notation we use below, however, is somewhat different and we instead follow the one found in \cite{AT}, namely, $\{a(n,L):n\in\N,\; L\in[\N]^\infty\}$. The $n-$averages $a(n,L)$ are defined as elements of $c_{00}(\N)$ in the following manner.

Let $(e_j)_j$ denote the unit vector basis of $c_{00}(\N)$ and $L\in[\N]^\infty$. For $n=0$ we define $a(0,L)=e_{l_1}$ where $l_1=\min L$. Suppose that $a(n,M)$ has been defined for some $n\in\N$ and every $M\in[\N]^\infty$. We define $a(n+1,L)$ in the following way: \\
We set $L_1=L$ and $L_k=L_{k-1}\setminus\supp( a(n,L_{k-1}))$ for $k=2,\ldots,l_1$ and finally define
\[
a(n+1,L)=\frac{1}{l_1}\big( a(n,L_1) + \cdots + a(n,L_{l_1}) \big).
\]

\begin{rem} Let $n\in\N$ and $L\in[\N]^\infty$. The following properties are easily established by induction.
\begin{itemize}
\item[(i)]	$a(n,L)$ is a convex combination of the unit vector basis of $c_{00}(\N)$.
\item[(ii)] $\|a(n,L)\|_{\ell_1}=1$ and $a(n,L)(k)\ge0$ for all $k\in\N$.
\item[(iii)] $\supp (a(n,L))$ is the maximal initial segment of $L$ contained in $\S_n$.
\item[(iv)] $\|a(n,L)\|_\infty= l_1^{-n}$ where $l_1=\min L$.
\item[(v)] If $\supp(a(n,L))=\{i_1<\ldots<i_d\}$ and $a(n,L)=\sum_{k=1}^da_{i_k}e_{i_k}$, then we have that $a_{i_1}\ge \ldots\ge a_{i_d}$.
\end{itemize}
\end{rem}
A proof of the following proposition can be found in \cite{AT}.
\begin{prop}\label{repeated averages proposition}
	Let $n\in\N$ and $L\in[\N]^\infty$. For every $F\in \S_{n-1}$, we have that
	\[
	\sum_{k\in F}a(n,L)(k)<\frac{3}{\min L}.
	\]
\end{prop}

\subsection{Special Convex Combinations}\label{scc subsection} Here we recall the notion of $(n,\varepsilon)$-special convex combinations, where $n\in\N$ and $\varepsilon>0$ (see \cite{AD2} and \cite{AT}).

\begin{dfn}
For $n\in\N$ and $\varepsilon>0$, a convex combination $\sum_{i\in F}c_ie_i$, of the unit vector basis $(e_i)_i$ of $c_{00}(\N)$ is called an $(n,\varepsilon)$-basic special convex combination (or an $(n,\varepsilon)$-basic s.c.c.) if 
\begin{itemize}
\item[(i)] $F\in\S_n$ and
\item[(ii)] for any $G\subset F$ with $G\in S_{n-1}$, we have that $\sum_{i\in G}c_i<\varepsilon$.
\end{itemize}
We will also call $\sum_{i\in F}c_i^{1/2}e_i$ a $(2,n,\varepsilon)$-basic special convex combination .
\end{dfn}
As follows from Proposition \ref{repeated averages proposition}, every $n$-average $a(n,L)$ is an $(n,3/\min L)$-basic s.c.c. and this yields the following. 
\begin{prop}\label{existenceofaverages}
	Let $M\in[\N]^{\infty}$, $n\in\N$ and $\varepsilon>0$. Then there is a $k\in\N$ such that for any $F\subset M$ such that $F$ is maximal in $\S_n$ and $k\le\min F$, there exists an $(n,\varepsilon)$-basic s.c.c. $x\in c_{00}(\N)$ with $\supp (x)=F$.
\end{prop}

Clearly, this also implies the existence of $(2,n,\varepsilon)$-basic special convex combinations by taking the square roots of the coefficients of an $(n,\varepsilon)$-b.s.c.c.

\begin{dfn}
Let $x_1<\ldots<x_d$ be vectors in $c_{00}(\N)$ and define $t_i=\min\supp (x_i)$, $i=1,\ldots,d$. We say that the vector $\sum_{i=1}^dc_ix_i$ is an $(n,\varepsilon)$-special convex combination (or an $(n,\varepsilon$)-s.c.c.) for some $n\in\N$ and $\varepsilon>0$ if  $\sum_{i=1}^d c_ie_{t_i}$ is an $(n,\varepsilon)$-basic s.c.c.,  and a $(2,n,\varepsilon)$-special convex combination if  $\sum_{i=1}^d c_ie_{t_i}$ is a $(2,n,\varepsilon)$-basic s.c.c.

\end{dfn}

\section{Asymptotic structures} 

Let us recall the definitions of the asymptotic notions that appear in the results of this paper and were mentioned in the introduction. Namely, asymptotic models, joint spreading models and the notions of Asymptotic $\ell_p$ and Asymptotic $c_0$ spaces. For a more thorough discussion, including several open problems and known results, we refer the reader to \cite[Section 3]{AM}.

\begin{dfn}[\cite{HOdell}]
An infinite array of sequences $(x^{i}_j)_j$, $i\in\N$, in a Banach space $X$, is said to generate a sequence $(e_i)_i$, in a seminormed space $E$, as an asymptotic model if for every $\varepsilon>0$ and $n\in\N$, there is a $k_0\in\N$ such that for any natural numbers $k_0\leq k_1<\cdots<k_n$ and any scalars $a_1,\ldots,a_n$ in $[-1,1]$ we have
\[\Big|\big\|\sum_{i=1}^na_ix_{k_i}^{i}\big\| - \big\|\sum_{i=1}^na_ie_{i}\big\|\Big| < \varepsilon.\]
\end{dfn}

A Banach space $X$ is said to admit a unique asymptotic model with respect to a family $\mathscr{F}$ of normalized sequences in $X$ if whenever two infinite arrays, consisting of sequences from $\mathscr{F}$, generate asymptotic models, then those must be equivalent. Typical families under consideration are those of normalized weakly null sequences, denoted  $\mathscr{F}_0(X)$, normalized Schauder basic sequences, denoted  $\mathscr{F}(X)$, or the family of all normalized block sequences of a fixed basis of $X$, if it has one, denoted  $\mathscr{F}_b(X)$.

	\begin{dfn}[\cite{AGLM}]
	Let $M\in[\N]^\infty$ and $k\in\N$. A {plegma} (resp., {strict plegma}) family in $[M]^k$ is a finite sequence $(s_i)_{i=1}^l$ in $[M]^k$ satisfying the following.
	\begin{enumerate}
		\item[(i)] $s_{i_1}(j_1)<s_{i_2}(j_2)$ for every $1\le j_1<j_2\le k$  and $1\le i_1,i_2\le l$.
		\item[(ii)] $s_{i_1}(j)\le s_{i_2}(j)$ $\big($resp., $s_{i_1}(j)< s_{i_2}(j)\big)$ for all $1\le i_1<i_2\le l$ and $1\le j\le k$ .
	\end{enumerate}
	For each $l\in \N$, the set of all sequences $(s_i)^l_{i=1}$ which are plegma families in $[M]^k$ will be denoted by $Plm_l([M]^k)$ and that of the strict plegma ones by $S$-$Plm_l([M]^k)$.
\end{dfn}

\begin{dfn}[\cite{AGLM}]
A finite array of sequences $(x^{i}_j)_j$, $1\leq i\leq l$, in a Banach space $X$, is said to generate another array of sequences $(e_j^{i})_j$, $1\leq i\leq l$, in a seminormed space $E$, as a joint spreading model if for every $\varepsilon>0$ and $n\in\N$, there is a $k_0\in\N$ such that for any $(s_i)_{i=1}^l\in S$-$Plm_{l}([\N]^n)$ with $k_0\le s_1(1)$ and for any $l\times n$ matrix $A=(a_{ij})$ with entries in $[-1,1]$, we have that
\[\Big|\big\|\sum_{i=1}^l\sum_{j=1}^na_{ij}x_{s_i(j)}^{i}\big\| - \big\|\sum_{i=1}^l\sum_{j=1}^na_{ij}e_j^{i}\big\|\Big|<\varepsilon.\]
\end{dfn}

 A Banach space $X$ is said to admit a  uniformly unique joint spreading model with respect to a family of normalized sequences $\mathscr{F}$ in $X$ if there exists a constant $C$ such that whenever two arrays $(x_j^{i})_j$ and $(y_j^{i})_j$, $1\leq i\leq l$, of sequences from $\mathscr{F}$ generate joint spreading models, then those must be $C$-equivalent. Moreover, a Banach space admits a uniformly unique joint spreading model with respect to a family $\mathscr{F}$ if and only if it admits a unique asymptotic model with respect to $\mathscr{F}$ (see, e.g.,  \cite[Remark 4.21]{AGLM} or \cite[Proposition 3.12]{AM}). In particular, if a space admits a uniformly unique joint spreading model with respect to some family $\mathscr{F}$ satisfying certain conditions described in \cite[Proposition 4.9]{AGLM}, then this is equivalent to some $\ell_p$. In order to show that a space admits some $\ell_p$ as a uniformly unique joint spreading model, it may be more convenient to prove (ii) of the following lemma, thereby avoiding the use of plegma families.
 
 \begin{lem}\label{ell_p as uniformly unique joint spreading equivalent form}
 Let $X$ be a Banach space and $\mathscr{F}$ be a family of normalized sequences in $X$. Let also $1\le p<\infty$. The following are equivalent.
 \begin{itemize}
 	\item[(i)] $X$ admits $\ell_p$ as a uniformly unique joint spreading model with respect to the family $\mathscr{F}$.
 	\item[(ii)] There exist constants $c,C>0$ such that for every array $(x^i_j)_j$, $1\le i\le l$, of sequences from $\mathscr{F}$, there is $M=\{m_1<m_2<\ldots\}$, an infinite subset of the naturals, such that for any choice of $1\le i_j\le l$, $j\in M$, every $F\subset M$ with $m_k\le F$ and $|F|\le k$ and any choice of scalars $a_j$, $j\in F$, 
 	\[
 	c \|(a_j)_{j\in F}\|_p \le \big\| \sum_{j\in F} a_j x^{i_j}_{j} \big\| \le C \|(a_j)_{j\in F}\|_p.
 	\]
 \end{itemize}
 \end{lem}
 \begin{proof}
 Note that (i) implies that there are constants $c,C>0$ such that for every array $(x^i_j)_j$, $1\le i\le l$, of sequences from $\mathscr{F}$, there is $N=\{n_1<n_2<\ldots \}$, an infinite subset of the naturals, such that for any $k$, any strict plegma family $(s_i)_{i=1}^l\in S\-Plm_l([N]^k)$ with $n_k\le s_1(1)$ and any $l\times k$ matrix $A=(a_{ij})$ of scalars we have that
  	\[
 	c \|(a_{ij})_{i=1,j=1}^{l,k}\|_p \le \big\| \sum_{i=1}^l\sum_{j=1}^k a_{ij} x^{i}_{s_i(j)} \big\| \le C \|(a_{ij})_{i=1,j=1}^{l,k}\|_p.
 	\]
 Let $N'=\{ n_{2kl}:k\in\N \}$ and observe that for $k_1,\ldots,k_d\in \N$, there is a strict plegma family $(s_i)_{i=1}^l\in S\-Plm_l([N]^d	)$ such that $n_{2k_jl}\in\{s_i(j):i=1,\ldots,l\}$ for all $j=1,\ldots, d$. Hence, we may find $M\subset N'$ satisfying (ii) with constants $c,C$. Finally, by repeating the sequences in the array, it follows easily that (ii) yields (i).
 \end{proof}
 
We recall the main result from \cite{AGLM}, stating that whenever a Banach space admits a uniformly unique joint spreading model with respect to some family satisfying certain stability conditions, then it satisfies a property concerning its bounded linear operators called the Uniform Approximation on Large Subspaces property (see \cite[Theorem 5.17]{AGLM} and \cite[Theorem 5.23]{AGLM}).
 
 \begin{dfn}[\cite{MMT}]
A Banach space $X$ is called Asymptotic $\ell_p$,  $1\leq p<\infty$, (resp., Asymptotic $c_0$) if there exists a constant $C$ such that in a two-player $n$-turn game $G(n,p,C)$, where in each turn $k=1,\ldots,n$, player (S) picks a finite codimensional subspace $Y_k$ of $X$, and then player (V) picks a normalized vector $x_k\in Y_k$, player (S) has a winning strategy to force player (V) to pick a sequence $(x_k)_{k=1}^n$ that is $C$-equivalent to the unit vector basis of $\ell^n_p$ (resp., $\ell_\infty^n)$.
\end{dfn}

Although this is not the initial formulation, it is equivalent and follows from \cite[Subsection 1.5]{MMT}. The typical example of a non-classical Asymptotic $\ell_p$ space is the Tsirelson space  from \cite{FJ}. This is a reflexive Asymptotic $\ell_1$ space, and it is the dual of Tsirelson's original space from \cite{T} which is  Asymptotic $c_0$. Finally, whenever a Banach space is Asymptotic $\ell_p$ or Asymptotic $c_0$, it admits a uniformly unique joint spreading model with respect to $\mathscr{F}_0(X)$ (see, e.g., \cite[Corollary 4.12]{AGLM}).

	 The above definition is the coordinate-free version of the notion of an asymptotic $\ell_p$ Banach space with a basis introduced by V. D. Milman and N. Tomczak-Jaegermann in \cite{MT}.

\begin{dfn}[\cite{MT}]
Let $X$ be a Banach space with a Schauder basis $(e_i)_i$ and $1\le p< \infty$. We say that the Schauder basis $(e_i)_i$ of $X$ is asymptotic $\ell_p$ if there exist positive constants $D_1$ and $D_2$ such that for all $n\in\N$, there exists $N(n)\in\N$ with the property that whenever $N(n)\le x_1< \cdots< x_n$ are vectors in $X$, then
\[
\frac{1}{D_1}(\sum_{i=1}^n\|x_i\|^p)^{\frac{1}{p}}\le \|\sum_{i=1}^nx_i\|\le D_2 (\sum_{i=1}^n\|x_i\|^p)^{\frac{1}{p}}.
\]
Specifically, we say that $(e_i)_i$ is $D$-asymptotic $\ell_p$ for $D=D_1D_2$. The definition of an asymptotic $c_0$ space is given similarly.
\end{dfn}

It is easy to show that if $X$ has a Schauder basis that is asymptotic $\ell_p$, then $X$ is Asymptotic $\ell_p$. Moreover, if $X$ is Asymptotic $\ell_p$, then it contains an asymptotic $\ell_p$ sequence. In particular, note that if $X$ has a Schauder basis and $Y$ is an Asymptotic $\ell_p$ subspace of $X$, then $Y$ contains a further subspace that is isomorphic to an asymptotic $\ell_p$ block subspace.

A noteworthy remark is that sequential asymptotic properties, array asymptotic properties, and global asymptotic properties of a Banach space $X$ can alternatively be interpreted as properties of special weakly null trees. A collection $\{x_A:A\in[\N]^{\leq n}\}$ in $X$ is said to be a normalized weakly null tree of height $n$, if for every $A\in [\N]^{\leq n-1}$, $(x_{A\cup\{j\}})_{j>\max(A)}$ is a normalized weakly null sequence. Such a tree is said to originate from a sequence $(x_j)_j$ if for all $A = \{a_1,\ldots,a_i\}$ we have $x_A = x_{a_i}$. Similarly, a tree $\{x_A:A\in[\N]^{\leq n}\}$ is said to originate from an array of sequences $(x^{(i)}_j)_j$, $1\leq i\leq n$, if for all $A = \{a_1,\ldots,a_i\}$ we have $x_A = x^{(i)}_{a_i}$. Then, $X$ has a uniformly unique $\ell_p$ spreading model if and only if there exists $C>0$ so that every tree $\{x_A:A\in[\N]^{\leq n}\}$ originating from a normalized weakly null sequence $(x_j)_j$ in $X$ has a maximal branch that is $C$-equivalent to the unit vector basis of $\ell_p^n$.  Similarly, $X$ has a unique $\ell_p$ asymptotic model if the same can be said about all trees originating from normalized weakly null arrays in $X$. Finally, a Banach space $X$ is an Asymptotic $\ell_p$ space (or an Asymptotic $c_0$ space if $p=\infty$) if there exists $C>0$ so that every normalized weakly null tree of height $n$ has a maximal branch $x_{\{a_1\}}, x_{\{a_1,a_2\}},\ldots,x_{\{a_1,a_2,\ldots,a_n\}}$ that is $C$-equivalent to the unit vector basis of $\ell_p^n$. For more details see \cite[Remark 3.11]{BLMS}.

\section{Measures on countably branching well-founded trees}\label{measures section}

	In this section we recall certain results from \cite{AGM} concerning measures on countably branching well-founded trees. These will be used to prove that for all $1\le p<\infty$, the space $\mathfrak{X}^{(p)}_{\text{awi}}$ admits $\ell_p$ as a unique asymptotic model. In particular, Proposition \ref{proposition measures 1 : incomparable} and Lemma \ref{splitting lemma} will be used to prove Lemma \ref{final measure lemma}, which is {one of the key ingredients} in the proof of the main result, Theorem \ref{main theorem intro}.
	 
	{
	Let $\T=(A,<_\T)$, where $A$ is a countably infinite set equipped with  a partial order $<_\T$. In the sequel, we use $t\in\T$ instead of $t\in A$. We assume that $<_\T$ is such that there is a unique minimal element in $\T$, and for each $t\in \T$, the set $S_t=\{s\in \T: s\le_\T t \}$ is finite and totally ordered, that is, $\T$ is a rooted tree. We also assume that $\T$ is well-founded, i.e., it contains no infinite totally ordered sets, and countably branching, i.e., every non-maximal node has countably infinite immediate successors.}
	 
	{ Observe that $\tT=(\{S_t :{t\in \T}\},<_\tT)$, where $<_\tT$ denotes inclusion, is also a tree, and that $\T$ is in fact isomorphic to $\tT$ via the mapping $t\mapsto S_t$. Given $t\in\T$, we will denote $S_t$ by $\tt$, identifying it as an element of $\tT$. For each $\tt\in\tT$ we denote by $S(\tt)$ the set of immediate successors of $\tt$ in $\tT$. In particular, if $\tt$ is maximal then $S(\tt)$ is empty. Moreover, for $\tt\in\tT$ we denote $V_\tt = \{\tilde{s}\in\tT:\tt\leq_\tT \tilde{s}\}$, and view $\tT$ as a topological space} with the topology generated by the sets $V_\tt$ and $\tT\setminus V_\tt$, $\tt\in\tT$, i.e., the pointwise convergence topology.  This is a compact metric topology such that for each $\tt\in\tT$, the sets of the form $V_\tt\setminus(\cup_{\tilde{s}\in F}V_{\tilde{s}})$, $F\subset S(\tt)$ finite, form a neighbourhood base of clopen sets for $\tilde t$.  We denote by $\mathcal{M}_+(\tT)$ the cone of all bounded positive measures $\mu:\mathcal{P}(\tT)\to[0,+\infty)$. For $\mu\in\mathcal{M}_+(\tT)$, we define the support of $\mu$ to be the set $\mathrm{supp}(\mu) = \{\tt\in\tT:\mu(\{\tt\})>0\}$.  Finally, we say that a subset $\mathcal{A}$ of $\mathcal{M}_+(\tT)$ is bounded if $\sup_{\mu\in\mathcal{A}}\mu(\tT)<\infty$. 
\begin{prop}\label{proposition measures 1 : incomparable}
				Let $(\mu_i)_i$ be a bounded and disjointly supported sequence in $\mathcal{M}_+(\tT)$. Then for every $\varepsilon>0$, there is an $L\in[\N]^\infty$ and a subset $G_i$ of $\supp(\mu_{i})$ for each $i\in L$,  satisfying the following.
		\begin{itemize}
			\item[(i)] $\mu_{i}(\tT\setminus G_i)\le \varepsilon$ for every $i\in L$.
			\item[(ii)] The sets $G_i$, $i\in L$, are pairwise incomparable.
		\end{itemize}
	\end{prop}
	For the proof, we refer the reader to  \cite[Proposition 3.1]{AGM}.
	
\begin{dfn}
Let $(\mu_i)_i$ {be a sequence} in $\mathcal{M}_+(\tT)$ and $\nu\in\mathcal{M}_+(\mathcal{\tT})$. We say that $\nu$  is the successor-determined limit of $(\mu_i)_i$ if for all $\tt\in \tT$, we have $\nu(\{\tt\}) = \lim_i\mu_i(S(\tt))$. In this case we write $\nu = \mathrm{succ}\-\!\lim_i\mu_i$.
\end{dfn}

Recall that {a bounded sequence $(\mu_i)_i$} in $\mathcal{M}_+(\tT)$ converges in the $w^*$-topology to a $\mu\in\mathcal{M}_+(\tT)$ if and only if for all clopen sets $V\subset \tT$ we have $\lim_i\mu_i(V) = \mu(V)$ if and only if for all $\tt\in\tT$ we have $\lim_i\mu_i(V_\tt) = \mu(V_\tt)$. In this case we write $\mu = w^*\-\lim_i\mu_i$.

{\begin{lem}
	Let $(\mu_i)_i$ be a bounded sequence in $\mathcal{M}_+(\tT)$. There exist a subsequence $(\mu_{i_n})_n$ of $(\mu_i)_i$ and $\nu\in\mathcal{M}_+(\tT)$ with $\nu = \mathrm{succ}\-\!\lim_n\mu_{i_n}$.
\end{lem}
}

\begin{rem}
It is possible for a bounded sequence $(\mu_i)_i$ in $\mathcal{M}_+(\tT)$ to satisfy $w^*\-\lim_i\mu_i\neq \mathrm{succ}\-\!\lim_i\mu_i$. Take for example $\tT = [\mathbb{N}]^{\leq 2}$ (all subsets of $\mathbb{N}$ with at most two elements with the partial order of initial segments) and define $\mu_i = \delta_{\{i,i+1\}}$, $i\in\N$. Then $w^*\-\lim_i\mu_i = \delta_\emptyset$ whereas $\mathrm{succ}\-\!\lim_i\mu_i = 0$.
\end{rem}

Although these limits are not necessarily the same, there is an explicit formula relating $\mathrm{succ}\-\!\lim_i\mu_i$ to $w^*\-\lim_i\mu_i$. 

\begin{lem}
\label{formula w-succ}
Let $(\mu_i)_i$ be a bounded and disjointly supported sequence in $\mathcal{M}_+(\tT)$ such that $w^*\-\lim_i\mu_i = \mu$ exists and for all $\tt\in\tT$  the limit $\nu(\{\tt\}) = \lim_i\mu_i(S(\tt)) $  exists as well. Then for every $\tt\in \tT$ and enumeration $(\tt_j)_j$ of $S(\tt)$, we have
\begin{equation}
\label{formula w-succ formula}
\mu(\{\tt\}) = \nu(\{\tt\}) + \lim_j\lim_i\mu_i\Big(\cup_{k\geq j}(V_{\tt_k}\setminus\{\tt_{k}\})\Big).
\end{equation}
In particular, $\mu(\{\tt\}) = \nu(\{\tt\})$ if and only if the double limit in \eqref{formula w-succ formula} is zero.
\end{lem}

\begin{lem}
\label{splitting lemma}
Let $(\mu_i)_i$ be a bounded and disjointly supported sequence in $\mathcal{M}_+(\tT)$ such that $\mathrm{succ}\-\!\lim_i\mu_i = \nu$  exists. Then there exist an infinite $L\subset\mathbb{N}$ and partitions $A_i$, $B_i$ of $\mathrm{supp}(\mu_i)$, $i\in L$, such that the following are satisfied.
\begin{itemize}

\item[(i)] If for all $i\in L$ we define the measure $\mu_i^1$ by  $\mu_i^1(C)= \mu_i(C\cap A_i)$, then $\nu = w^*\-\lim_{i\in L}\mu_i^1 = \mathrm{succ}\-\!\lim_{i\in L}\mu_i^1$.

\item[(ii)] If for all $i\in L$ we define the measure $\mu_i^2$ by  $\mu_i^2(C)= \mu_i(C\cap B_i)$, then for all $\tt\in\tT$, the sequence $(\mu_i^2(S(\tt)))_i$ is eventually zero. In particular, $\mathrm{succ}\-\!\lim_{i\in L}\mu_i^2 = 0$.

\end{itemize}
\end{lem}

For the proofs, we refer the reader to \cite[Lemma 4.10]{AGM} and \cite[Lemma 4.12]{AGM}.

{\begin{rem}
Although the results from \cite{AGM} were formulated for trees $\mathcal{T}$ defined on infinite subsets of $\N$, this is not a necessary restriction and they can be naturally extended to the more general setting of countably branching well-founded trees.
\end{rem}}


{\Large
\part{ The case of $\boldsymbol{\ell_1}$}
}

\section{Definition of the space $\mathfrak{X}_{\textnormal{awi}}^{(1)}$}\label{section 1}
\label{ell1section}
The method of saturation with asymptotically weakly incomparable constraints, that is used in the construction of both spaces presented in this paper, was introduced in \cite{AGM} where it was shown that (b)$_1\not\Rightarrow$(a)$_1$. There, it was also used to prove an even stronger result, namely, the existence of a Banach space with a basis admitting $\ell_1$ as a unique asymptotic model, and in which any infinite subsequence of the basis generates a non-Asymptotic $\ell_1$ subspace. This method requires the existence of a well-founded tree defined either on the basis of the space or on a family of functionals of its norming set. In this section we define the space $\mathfrak{X}_{\text{awi}}^{(1)}$ by introducing its norm via a norming set, which is a subset of the norming set of a Mixed Tsirelson space $\T[(m_j,\S_{n_j})_j]$ for an appropriate choice of $(m_j)_j$ and $(n_j)_j$ described below. The key ingredient in the definition of this norming set is the notion of asymptotically weakly incomparable sequences of functionals, which is also introduced in this section. This notion will allow the space $\mathfrak{X}_{\text{awi}}^{(1)}$ to admit $\ell_1$ as a unique asymptotic model while at the same time it will force the norm to be small on the branches of a tree, in every subspace of $\mathfrak{X}_{\text{awi}}^{(1)}$, showing that the space does not contain Asymptotic $\ell_1$ subspaces.

\subsection{Definition of the space $\mathfrak{X}_{\textnormal{awi}}^{(1)}$}
	Define a pair of strictly increasing sequences of natural numbers $(m_j)_j$, $(n_j)_j$ as follows:
		\begin{align*}
			m_1 &= 2 & n_1 &= 1 \\
			m_{j+1}&=m_j^{m_j} & n_{j+1}&=2^{2m_{j+1}}n_j
		\end{align*}
	\begin{dfn}
		Let $V_{(1)}$ denote the minimal subset of $c_{00}(\N)$ that
		\begin{itemize}
			\item[(i)] contains {0} and all $\pm e_j^*$, $j\in \N$, and
			\item[(ii)] is closed under the operations $(m_j,\S_{n_j})_j$, i.e., if $j\in\N$ and $f_1<\ldots<f_n$ is an $\S_{n_j}$-admissible sequence (see Section \ref{comb}) in $V_{(1)}\setminus\{0\}$ then $m_j^{-1}\sum_{i=1}^nf_i $ is also in $V_{(1)}$.
		\end{itemize}
	\end{dfn}

	\begin{rem}
		\begin{itemize}
			\item[(i)] If $f\in V_{(1)}\setminus\{0\}$, then either $f\in\{\pm e^*_j:j\in\N\}$, or it is of the form $f=m_j^{-1}\sum_{i=1}^nf_i$ with $f_1<\ldots<f_n$ an $\S_{n_j}$-admissible sequence in $V_{(1)}$ for some $j\in\N$.
			\item[(ii)] As usual, we view the elements of  $V_{(1)}$ as functionals acting on $c_{00}(\N)$, inducing a norm $\|\cdot\|_{V_{(1)}}$. The completion of $(c_{00}(\N),\|\cdot\|_{V_{(1)}})$ is the Mixed Tsirelson space $\T[(m_j,\S_{n_j})_j]$ introduced for the first time in \cite{AD2}. The first space with a saturated norm defined by a countable family of operations is the Schlumprecht space \cite{S} which is a fundamental discovery and was used by W. T. Gowers and B. Maurey \cite{GM} to define the first hereditarily indecomposable (HI) space. 
		\end{itemize}
	\end{rem}
	
	We now recall the notion of tree analysis which appeared for the first time in \cite{AD}. This has become a standard tool in proving upper bounds for the estimations of functionals on certain vectors in Mixed Tsirelson spaces. However, it is the first time where the tree analysis has a significant role in the definition of the norming set $W_{(1)}$. Additionally, it is also a key ingredient in the proof that $\mathfrak{X}^{(1)}_{\text{awi}}$ contains no Asymptotic $\ell_1$ subspaces.
	
	Let $\mathcal{A}$ be a rooted tree. For a node $\alpha\in\mathcal{A}$, we denote by $S(\alpha)$ the set of all immediate successors of $\alpha$, by $|\alpha|$ the height of $\alpha$, i.e., $|\alpha|=\#\{\beta\in\mathcal{A}:\beta<_{\mathcal{A}}\alpha\}$, and finally we denote by $h(\mathcal{A})$ the height of $\mathcal{A}$, that is, the maximum height over its nodes.

	\begin{dfn}\label{tree analysis}
		Let $f\in V_{(1)}\setminus\{0\}$. For a finite tree $\mathcal{A}$, a family $(f_\alpha)_{\alpha\in \mathcal{A}}$ is called a tree analysis of $f$ if the following are satisfied.
		\begin{itemize}
			\item[(i)] $\mathcal{A}$ has a unique root denoted by $0$ and $f_0=f$.
			\item[(ii)] Each $f_\alpha$ is in $V_{(1)}$ and if $\beta<\alpha$ in $\mathcal{A}$ then $\range(f_\alpha)\subset \range(f_\beta)$.
			{
			\item[(iii)] For every maximal node $\alpha\in\mathcal{A}$  we have that $|\alpha|=h(\mathcal{A})$.

			\item[(iv)] For every non-maximal node $\alpha\in\mathcal{A}$, either $f_\alpha$ is the result of some $(m_j,\S_{n_j})$ operation of $(f_\beta)_{\beta\in S(\alpha)}$, i.e., $f_\alpha=m_j^{-1}\sum_{\beta\in S(\alpha)} f_\beta$, or $f_\alpha\in\{\pm e_j^*:j\in\N\}$ and $S(\alpha)=\{\beta\}$ with $f_\beta=f_\alpha$.}
			\item[(v)] For every maximal node $\alpha\in\mathcal{A}$, $f_\alpha\in\{\pm e^*_j:j\in\N\}.$
		\end{itemize}
	\end{dfn}

	\begin{rem}
		\begin{itemize}
			\item[(i)] It follows by minimality that every $f$ in $V_{(1)}\setminus\{0\}$ admits a tree analysis, but it may not be unique. For example, $f=(m_{1}m_2)^{-1}e^*_1$ admits two distinct tree analyses.
			\item[(ii)] The standard definition of a tree analysis does not include \ref{tree analysis} (iii). This property is included for technical reasons, and is used below in the equality of Remark \ref{w f f alpha remark} (i).
		\end{itemize}		
	\end{rem}

	\begin{dfn}\label{weight of functional l1}
		Let $f\in V_{(1)}$.
		\begin{itemize}
			\item[(i)] If $f=0$ or $f\in\{\pm e^*_j:j\in\N\}$, then we define the weight $w(f)$ of $f$ as $w(f)=0$ and $w(f)=1$, respectively.
			\item[(ii)] If $f$ is the result of an $(m_j,S_{n_j})$-operation for some $j\in\N$, then $w(f)=m_j$.
		\end{itemize}
	\end{dfn}
	\begin{rem}\label{remarks on w(t)}
	 It is not difficult to see that $w(f)$, for $f\in V_{(1)}$, is not uniquely determined, i.e., $f$ could be the result of more than one distinct $(m_j,\S_{n_j})$-operations. However, if we fix a tree analysis $(f_\alpha)_{\alpha\in\mathcal{A}}$ of $f$, then for $\alpha\in\mathcal{A}$ with $f_\alpha=(m_{j_\alpha})^{-1}\sum_{\beta\in S(\alpha)}f_\beta$, the tree analysis determines the weight $w(f_\alpha)$, being equal to $m_{j_\alpha}$. Thus, for $f\in V_{(1)}$ and a fixed tree analysis $(f_\alpha)_{\alpha\in\mathcal{A}}$ of $f$, with $w(f_\alpha)$ we will denote the weight $m_{j_\alpha}$ determined by $(f_\alpha)_{\alpha\in\mathcal{A}}$, for every $\alpha\in\mathcal{A}$. In addition, we will denote by $\bar{f}_\alpha$ the pair $(f_\alpha,m_{j_\alpha})$. 
	\end{rem}
	
	\begin{dfn}
		Let $f\in V_{(1)}$ and $(f_\alpha)_{\alpha\in\mathcal{A}}$ be a tree analysis of $f$. Then for $\alpha\in\mathcal{A}$, we define the relative weight $w_f(f_\alpha)$ of $f_\alpha$ as
		\[
		w_f(f_\alpha) = \begin{cases}
			\prod_{\beta<\alpha}w(f_\beta) & \text{if }\alpha\neq0\\
			1 &\text{otherwise.}
		\end{cases}
		\]
	\end{dfn}

	\begin{rem}\label{w f f alpha remark}	

		Let $f\in V_{(1)}$ and $(f_\alpha)_{\alpha\in\mathcal{A}}$ be a tree analysis of $f$.
		\begin{itemize}
			\item[(i)]  For every $k=1,\ldots,h(\mathcal{A})$ 
			\[
			f=\sum_{|a|=k}w_f(f_\alpha)^{-1}f_\alpha.
			\]
			This can be proved by induction and essentially relies on the fact that $(f_\alpha)_{\alpha\in\mathcal{A}}$ satisfies  \ref{tree analysis} (iii).
			\item[(ii)] If $\mathcal{B}$ is a maximal pairwise incomparable subset of $\mathcal{A}$, then
			\[
			f=\sum_{\beta\in\mathcal{B}}w_f(f_\beta)^{-1}f_\beta.
			\]
			\item[(iii)] For every $\alpha\in\mathcal{A}$, whose immediate predecessor $\beta$ in $\mathcal{A}$ (if one exists) satisfies $f_\beta\notin \{\pm e_j^*:j\in\N\}$, we have $w_f(f_\alpha)\ge 2^{|\alpha|}.$
		\end{itemize}					 
	\end{rem}

	Fix an injection $\sigma$ that maps any pair $(f,w(f))$, for $f\in V_{(1)}$ and $w(f)$ a weight of $f$, to some $m_j$ with $m_j>\max\supp (f)\: w(f)$ whenever $f\neq 0$.

	\begin{dfn}\label{definition of the tree}
		Define a partial order $<_\T$ on the set of all pairs $(f,w(f))$ for $f\in V_{(1)}$ and $w(f)$ a weight of $f$, as follows: $(f,w(f))<_\T (g ,w(g))$ either if $f=0$ or if there exist $f_1<\ldots<f_n \in V_{(1)}$ and weights $w(f_1),\ldots,w(f_n)$ such that 
		\begin{itemize}
			\item[(i)] $(f_i)_{i=1}^n$ is $\S_1$-admissible,
			\item[(ii)] $w(f_1)=\sigma(0,0)$ and $w(f_i)=\sigma(f_{i-1},w(f_{i-1}))$ for every $i=2,\ldots,n$,
			\item[(iii)] there are $1\le i_1<i_2\le n$ such that $f=f_{i_1}$ and $g=f_{i_2}$.
		\end{itemize}
	\end{dfn}

	It is easy to see that $<_\T$ induces a tree structure rooted at $\bar{0}=(0,0)$. Let us denote this tree by $\T$ and observe that this is a countably branching well-founded tree, due to \ref{definition of the tree}(i). For $t=(f,w(f))\in\T$, we set $f_t=f$ and $w(t)=w(f)$. 
	
	It is clear that unlike the case where the tree is defined on the basis of the space, here incomparable segments need not necessarily have disjoint supports. This forces us to introduce the notion of essentially incomparable nodes, which was first defined in \cite{AGM}. To this end, we first need to define an additional tree structure that is readily implied by $\T$ via the projection $(f,w(f))\mapsto w(f)$.

	\begin{dfn}
		Define a partial order $<_{\W}$ on $\{ m_j:j\in\N \}$ as follows: $m_i<_\W m_j$ if there exist $t_{1},t_{2}\in \T$ such that $t_1<_\T t_2$, $w(t_1)=m_i$ and $w(t_2)=m_j$.
	\end{dfn}

	As an immediate consequence of the fact that $\T$ is a countably branching well-founded tree, we have that $<_\W$ also defines a tree structure. Let us denote this tree by $\W$ and note that it is also countably branching and well-founded. 

	\begin{rem}
	The above definition implies that if $t_1,t_2\in\T$ are such that $w(t_1)<_\W w(t_2)$, then there exist $t_3,t_4\in\T$ such that $t_3<_\T t_4$, $w(t_3)=w(t_1)$ and $w(t_4)=w(t_2)$. The tree structure of $\T$ implies that $t_3$ is uniquely defined, and we will say that $t_3$ {generates} $w(t_2)$. This is not the case however for $t_4$ and moreover, it is not necessary that $t_3<_\T t_2$.
	\end{rem}

	\begin{dfn} \label{incomparabilities 1}
		\begin{itemize} 
			\item[(i)] A subset $A$ of $\T\setminus\{\bar{0}\}$ is called  essentially incomparable if whenever $t_1,t_2\in A$ are such that $w(t_1)<_\W w(t_2)$,  then for the unique $t_3\in \T$ with $w(t_3)=w(t_1)$ that {generates} $w(t_2)$, we have that $f_{t_3}<f_{t_1}$.
			\item[(ii)] A subset $A$ of $\T$ is called weight incomparable if for any $t_1\neq t_2$ in $A$, $w(t_1)\neq w(t_2)$ and the weights $w(t_1)$ and $w(t_2)$ are incomparable in $\W$.
			\item[(iii)] A sequence $(A_j)_j$ of subsets of $\T$ is called pairwise weight incomparable if for every $j_1\neq j_2$ in $\N$, $t_1\in A_{j_1}$ and $t_2\in A_{j_2}$, $w(t_1)\neq w(t_2)$ and the weights $w(t_1)$ and $w(t_2)$ are incomparable in $\W$.
		\end{itemize}
	\end{dfn}

	\begin{rem}\label{remarks on the various incomparabilites}
		\begin{itemize}
			\item[(i)] If $A$ is an essentially (resp., weight) incomparable subset of $\T$, then every $B\subset A$ is also essentially (resp., weight) incomparable.
			\item[(ii)] Any subsequence of a pairwise weight incomparable sequence in $\T$ is also pairwise weight incomparable.
			\item[(iii)] Any weight incomparable subset of $\T$ is essentially incomparable.
			\item[(iv)] Let $A=\{(f,1):f\in \{\pm e^*_j:j\in\N\} \}$. Then $A$ is essentialy incomparable and additionally, if $B\subset\T$ is essentially incomparable then the same holds for $A\cup B$.
		\end{itemize}		
	\end{rem}

We can finally describe the rule used to define the norming set $W_{(1)}$ of $\mathfrak{X}_{\text{awi}}^{(1)}$, namely, asymptotically weakly incomparable constraints.

	\begin{dfn}\label{incomparabilities 2}
		Let $J$ be an initial segment of $\N$ or $ J =\N$. Then a sequence $(f_j)_{j\in J}$ of functionals with successive supports in $V_{(1)}\setminus\{0\}$ is called asymptotically weakly incomparable (AWI) if each $f_j$ admits a tree analysis $(f_{j,\alpha})_{\alpha\in \mathcal{A}_j}$, $j\in J$, such that the following are satisfied.
		\begin{itemize}
			\item[(i)] There is a partition $\{\barf_j:j\in J\}= C_1^0\cup C_2^0$ such that $C^0_1$ is essentially incomparable and $C^0_2$ is weight incomparable.
			\item[(ii)] For every $k,j\in J$ with $j\ge k+1$, there exists a partition
			\[
				\{\barf_{j,\alpha}:\alpha\in \mathcal{A}_j\text{ and }|\alpha|=k\}=C^{k}_{1,j}\cup C^{k}_{2,j}
			\]
			such that $\cup_{j=k+1}^\infty C^k_{1,j}$ is essentially incomparable and $(C^{k}_{2,j})_{j=k+1}^\infty$ is pairwise weight incomparable.
		\end{itemize} 
	\end{dfn}
	
\begin{figure}[h]\centering
		\includegraphics[scale=0.5]{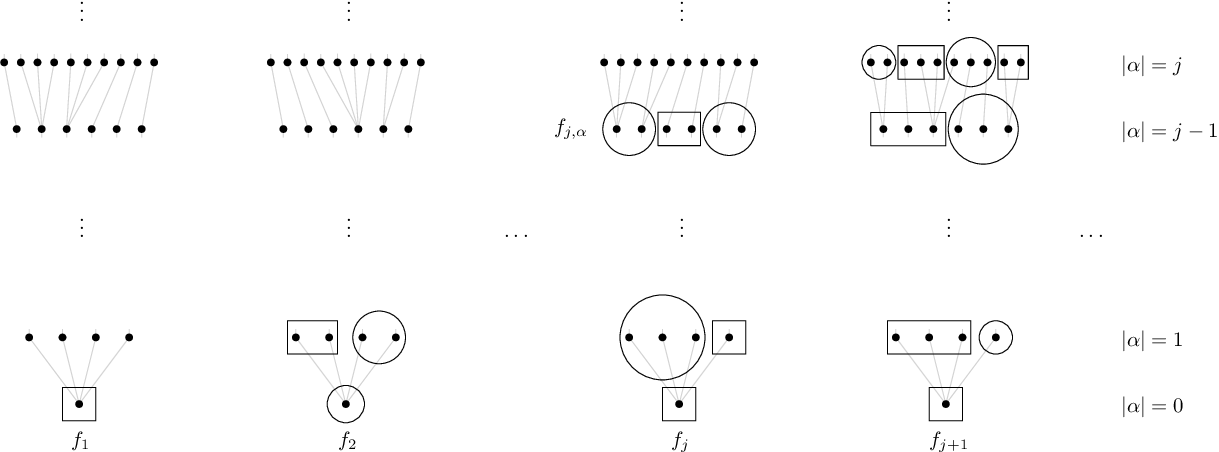}
		\caption{The collection of nodes of a fixed level in rectangles across all tree analyses forms an essentially incomparable subset, while circles across a fixed level form a family of pairwise weight incomparable subsets.}
			\end{figure}

Before defining the space $\mathfrak{X}_{\text{awi}}^{(1)}$ we prove that AWI sequences are stable under taking subsequences and under taking restrictions of functionals to subsets. This fact will imply the unconditionality of the basis of $\mathfrak{X}_{\text{awi}}^{(1)}$.

	\begin{rem}\label{functionals restricted on a subset of the supp}
		Let $f\in V_{(1)}$ and $(f_\alpha)_{\alpha\in\mathcal{A}}$ be a tree analysis of $f$. Let $\Delta$ be a non-empty subset of $\supp(f)$ and set $g=f|_\Delta$. First, note that $g\in V_{(1)}$. Moreover, $(f_\alpha)_{\alpha\in\mathcal{A}}$ naturally induces a tree analysis $(g_\alpha)_{\alpha\in\mathcal{B}}$ for $g$ as follows: $\mathcal{B}=\{\alpha\in\mathcal{A}:\supp(f_\alpha)\cap \Delta\neq\emptyset \}$ and $g_\alpha=f_\alpha|_\Delta$, $\alpha\in\mathcal{B}$. Finally, it is easy to see that $w(g)=w(f)$.
	\end{rem}

	\begin{prop}\label{proposition AWI - subsequence etc}
		Let $J$ be an initial segment of $\N$ or $ J =\N$ and $(f_j)_{j\in J}$ be an AWI sequence in $V_{(1)}$.
		\begin{itemize}
			\item[(i)] Every subsequence of $(f_j)_{j\in J}$ is also an AWI sequence in $V_{(1)}$.
			\item[(ii)] If $\Delta_j$ is a non-empty subset of $\supp (f_j)$ and $g_j=f_j|_{\Delta_j}$, $j\in J$, then $(g_j)_{j\in J}$ is an AWI sequence in $V_{(1)}$. 
			\item[(iii)] If $(g_j)_{j\in J}$ is a sequence in $V_{(1)}$ such that $|g_j|=|f_j|$ for all $j\in J$, then $(g_j)_{j\in J}$ is also AWI.
		\end{itemize}
	\end{prop}
	\begin{proof}
		Let for every $j\in J$, $(f_{j,\alpha})_{\alpha\in \mathcal{A}_j}$ be a tree analysis of $f_j$ with
		\[
			\{ \barf_j:j\in J \}= C^0_1\cup C^0_2
		\]
		and for every $k,j\in J$ with $j> k$
		\[
				\{\barf_{j,\alpha}:\alpha\in \mathcal{A}_j\text{ and }|\alpha|=k\}=C^{k}_{1,j}\cup C^{k}_{2,j},
		\]
		witnessing that $(f_j)_{j\in J}$ is AWI. We will define the desired partitions proving the cases (i)--(iii).
		
		To prove (i), let $N$ be a subset of $J$ and define
		\[
			F^0_i=\{\barf_j:j\in N\}\cap C^0_i,\quad i=1,2.
		\]
		Then $\{\barf_j:j\in N \} = F^0_1\cup F^0_2$, where $F^0_1$ is essentially incomparable and $F^0_2$ is weight incomparable. For the remaining part, let $k\in N$ and note that for $N_k=\{j\in N:j\ge k\}$,
		$\cup_{j\in N_k} C^{k}_{1,j}$ is essentially incomparable and $(C^{k}_{2,j})_{j\in N_k}$ is pairwise weight incomparable.   

		To prove (ii), Remark \ref{functionals restricted on a subset of the supp} implies that $g_j\in V_{(1)}$, $w(g_j)=w(f_j)$,  and we let $(g_{j,\alpha})_{\alpha\in\mathcal{B}_j}$ be the tree analysis of $g$ induced by $(f_{j,\alpha})_{\alpha\in \mathcal{A}_j}$, $j\in J$. Define
		\[ 
			F^0_i=\{\barg_j:j\in J\text{ and }g_j=f_j|_{\Delta_j}\text{ with }\barf_j\in C^0_i \}, \quad i=1,2,
		\]
		and observe that $\{\barg_j:j\in J \}=F^0_1\cup F^0_2$. Moreover, for $j\in J$, $\supp (g_j)\subset \supp (f_j)$ and $w(g_j)=w(f_j)$, and hence whenever $g_i\neq g_j$ are in $F^0_1$ with $w(g_i)<_\W w(g_j)$ we have $w(f_i)<_\W w(f_j)$, implying that the generator $t_3\in\T$ of $w(f_j)=w(g_j)$ with $w(t_3)=w(f_i)=w(g_i)$ is such that $f_{t_3}<f_i$, and thus $f_{t_3}<g_i$. This yields that $F^0_1$ is essentially incomparable. Clearly, $F^0_2$ is weight incomparable. Next, for $k,j\in J$ with $j>k$, define
		\[
			F^k_{i,j}=\{\barg_{j,\alpha}:g_{j,\alpha}=f_{j,\alpha}|_{\Delta_j}\text{ with }\barf_{j,\alpha}\in C^k_{i,j} \}, \quad i=1,2.
		\]
		Note that for each $k\in J$, $(F^{k}_{2,j})_{j=k+1}^\infty$ is pairwise weight incomparable, and the proof that $\cup_{j=k+1}^\infty F^k_{1,j}$ is essentially incomparable is identical to that for $F^1_0$. Finally, the proof of (iii) is similar that of (ii).
	\end{proof}

	\begin{dfn}
		Let $W_{(1)}$ be the smallest subset of $V_{(1)}$ that is symmetric, contains the singletons and for every $j\in\N$ and every $\S_{n_j}$-admissible AWI sequence $(f_i)_{i=1}^n$ in $W_{(1)}$ we have that $m_j^{-1}\sum_{i=1}^nf_i$ is in $W_{(1)}$. Moreover, let $\mathfrak{X}_{\text{awi}}^{(1)}$ denote the completion of $c_{00}(\N)$ with respect to the norm induced by $W_{(1)}$.
	\end{dfn}
	
\begin{rem}\label{remarks on the norming set W_{(1)}}
\begin{itemize}
\item[(i)] The norming set $W_{(1)}$ can be defined as the increasing union of a sequence $(W^n_{(1)})_{n=0}^\infty$, where $W^0_{(1)}=\{\pm e^*_k:k\in\N\}\cup\{0\}$ and
\begin{align*}
W^{n+1}_{(1)}=W^n_{(1)}\cup\Big\{\frac{1}{m_j}\sum_{l=1}^df_l:\;&j,d\in\N\text{ and }(f_l)_{l=1}^d \text{ is an } \\& \S_{n_j}\text{-admissible AWI sequence in }W^n_{(1)}\Big\}.
\end{align*}
\item[(ii)] Note that Remark \ref{remarks on the various incomparabilites} (iv) implies that any sequence of singletons is AWI. Hence, we have that
\[
	W^{1}_{(1)}=W^0_{(1)}\cup\Big\{\frac{1}{m_j}\sum_{k\in E}\epsilon_ke^*_k:\;j\in\N, \; E\in\mathcal{S}_{n_j} \text{ and } \epsilon_k\in\{-1,1\} \text{ for }k\in E \Big\}.
\]
\item[(iii)] Proposition \ref{proposition AWI - subsequence etc} yields that the standard unit vector basis of $c_{00}(\N)$ forms an $1$-unconditional Schauder basis for $\mathfrak{X}_{\text{awi}}^{(1)}$.
		\end{itemize}
	\end{rem}

\section{Outline of proof}
Although unconditionality of the basis of $\mathfrak{X}_{\text{awi}}^{(1)}$ is almost immediate, it is not however straightforward to show that $\mathfrak{X}_{\text{awi}}^{(1)}$ admits $\ell_1$ as an asymptotic model. Indeed, this requires Lemma \ref{final measure lemma}, which is based on the combinatorial results concerning measures on well-founded trees of Section 4, which first appeared in \cite{AGM}. This lemma yields that for any choice of successive families $(F_j)_j$ of normalized blocks in $\mathfrak{X}_{\text{awi}}^{(1)}$ and for any $\varepsilon>0$, we may pass to a subsequence $(F_j)_{j\in M}$ and find a family $(G_j)_{j\in M}$ of subsets of $W_{(1)}$ such that for any choice of $x_j\in F_j$, $j\in M$, there is a $g_j\in G_j$ with $g_j(x_j)>1-\varepsilon$ so that $(g_j)_{j\in M}$ is AWI. Thus, we are able to prove, employing Lemma \ref{ell_p as uniformly unique joint spreading equivalent form}, the aforementioned result.

To prove the non-existence of Asymptotic $\ell_1$ subspaces in $\mathfrak{X}_{\text{awi}}^{(1)}$, we start with the notion of exact pairs. This is a key ingredient in the study of Mixed Tsirelson spaces, used for the first time by Th. Schlumprecht \cite{S}. 
	
	\begin{dfn}
	We call a pair $(x,f)$, where $x\in\mathfrak{X}_{\text{awi}}^{(1)}$ and $f\in W_{(1)}$, an $m_j$-exact pair if the following hold.
	\begin{itemize}
	\item[(i)] $\|x\|\le 3$, $f(x)=1$ and $w(f)=m_j$.
	\item[(ii)] If $g\in W_{(1)}$ with $w(g)<w(f)$, then $|g(x)|\le 18w(g)^{-1}$.
	\item[(iii)] If $g\in W_{(1)}$ with $w(g)>w(f)$, then $|g(x)|\le 6(m_{j}^{-1}+m_j w(g)^{-1}) $.
	\end{itemize}
	If additionally for every $g\in W_{(1)}$ that has a tree analysis $(g_\alpha)_{\alpha\in\mathcal{A}}$ such that $w(g_\alpha)\neq m_j$ for all $\alpha\in\mathcal{A}$, we have $|g(x)|\le 18m^{-1}_j$, then we call $(x,f)$ a strong exact pair.
	\end{dfn}
That is, roughly speaking, for an exact pair $(x,f)$ the evaluation of a functional $g$ in $W_{(1)}$, on $x$, admits an upper bound depending only on the weight of $g$. In the case of an $m_j$-strong exact pair $(x,f)$, any $g$ in $W_{(1)}$ with a tree analysis $(g_\alpha)_{\alpha\in\mathcal{A}}$ such that $w(g_\alpha)\neq m_j$, has negligible evaluation on $x$.   We will consider certain exact pairs which we call standard exact pairs (SEP) (see Definition \ref{SEP definition}) and which we prove to be strong exact pairs. It is the case that such pairs can be found in any block subspace of $\mathfrak{X}_{\text{awi}}^{(1)}$, and this is used to prove the reflexivity of $\mathfrak{X}_{\text{awi}}^{(1)}$ as well as the following proposition which yields the non-existence of Asymptotic $\ell_1$ subspaces.
	
	\begin{prop}\label{not contain as l1 intro} Given $0<c<1$ there is $n\in\N$ so that in any block subspace $Y$ there is a sequence $(x_1,f_1),\ldots,(x_n,f_n)$ of SEPs, where $x_i\in Y$, $i=1,\ldots,n$, with $\barf_1<_\T\ldots<_\T \barf_n$ such that $\|x_1+\cdots+x_n\|<c\:n$.
	\end{prop}
	
	 To this end, we first employ the following lemma that highlights the importance of the asymptotically weakly incomparable constraints.
	
	\begin{lem}\label{ekf intro}
	Let $(x_1,f_1),\ldots,(x_n,f_n)$ be SEPs with $\barf_1<_\T\ldots<_\T \barf_n$. Then, for any $f\in W_{(1)}$ with a tree analysis $(f_\alpha)_{\alpha\in \mathcal{A}}$ and $k\in\N$, the number of $f_i$'s, $i=1,\ldots,n$, such that there exists $\alpha\in\mathcal{A}$ with $|\alpha|=k$, $w(f_i)=w(f_\alpha)$ and $\supp(x_i)\cap \supp(f_\alpha)\neq\emptyset$, is at most $ek!$, where $e$ denotes Euler's number.
	\end{lem}
	
\begin{figure}[h]
		\includegraphics[scale=0.62]{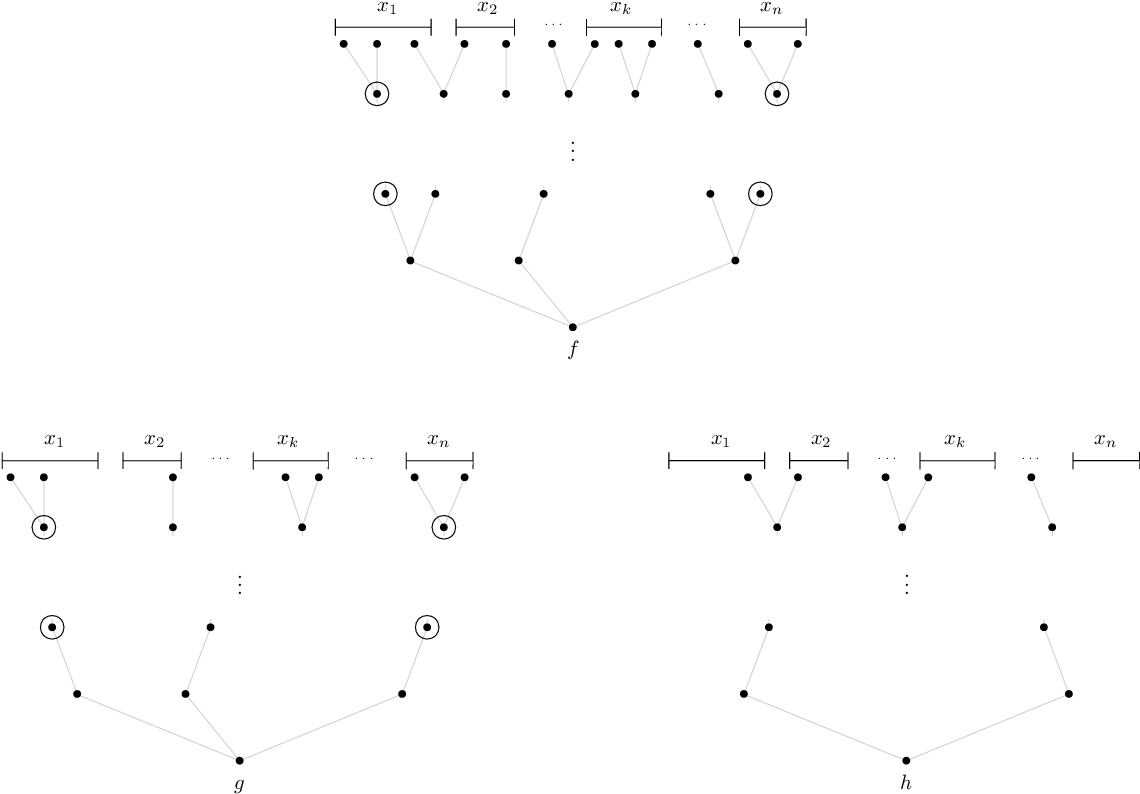}
			\caption{The tree analysis of $f$ and the induced tree analyses of $g$ and $h$. The circled nodes $\alpha$ are such that $w(f_\alpha)=w(f_k)$ and $\supp(x_k)\cap\supp(f_\alpha)\neq\emptyset$ for some $k\in\{1,\ldots,n\}$.}
			\label{ghfigure}
	\end{figure}

	Then, we consider a sequence of standard exact pairs $(x_1,f_1),\ldots,(x_n,f_n)$ with $\barf_1<_\T\ldots<_\T \barf_n$ and fix $0<c<1$. Pick an $m\in\N$ such that $3/2^m<c$. For $f\in W_{(1)}$, with a tree analysis $(f_\alpha)_{\alpha\in \mathcal{A}}$, we consider partitions $f=g+h$ and $g=g_1+g_2$ as follows: First, set
		\[
		G=\cup\{\range (x_k)\cap \range (f_\alpha): k\in\{1,\ldots,n\}\text{ and }\alpha\in\mathcal{A}\text{ with }w(f_\alpha)=w(f_k)\}.
		\]
		and define  $g=f|_G$ and $h=f|_{\N\setminus G}$ (see Figure \ref{ghfigure}).

\begin{figure}[h]
		\includegraphics[scale=0.62]{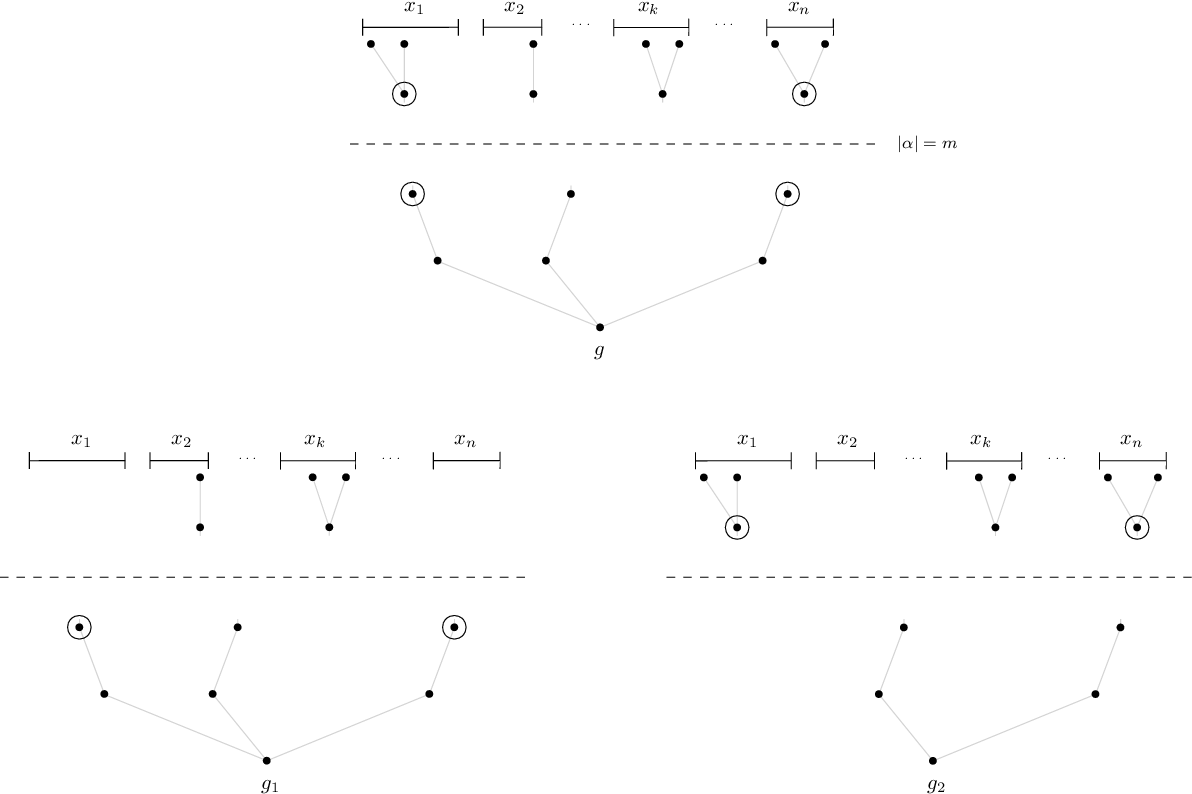}
			\caption{We consider the $m$-th level of the induced tree analysis of $g$. Nodes $\alpha$ with  $w(f_\alpha)=w(f_k)$ and $|\alpha|\le m$ are used to define $g_1$, while such nodes of height greater than $m$ define $g_2$, restricted on each $x_k$ for $k=1,\ldots,n$.		}
			\label{g1g2figure}
	\end{figure}

To define $g_1$, consider the tree analysis $(g_\alpha)_{\alpha\in\mathcal{A}_g}$ of $g$ that is induced by $(f_\alpha)_{\alpha\in\mathcal{A}}$, i.e., $g_\alpha=f_\alpha|_G$ for $\alpha\in\mathcal{A}_g$ and $\mathcal{A}_g=\{\alpha\in\mathcal{A}:\supp(f_\alpha)\cap G\neq\emptyset\}$. 	Then, we define
		\begin{align*}
\mathcal{B}^1_k=\{\alpha\in\mathcal{A}_g:&\;|a|\le m,\; w(f_\alpha)=w(f_k)\text{ and } w(f_\beta)\neq w(f_k)\text{ for all }\beta<\alpha\text{ in }\mathcal{A}_g\}
\end{align*}
for $k=1,\ldots,n$,
\[
G_1=\cup_{k=1}^n\cup\left\{\supp(g_\alpha)\cap \supp(x_k):\alpha\in\mathcal{B}^1_k\right\},
\]
and finally $g_1 = g|_{G_1}$ (see Figure \ref{g1g2figure}). Observe that Lemma \ref{ekf intro} implies that 
 \begin{equation}\label{intro1 eq 1}
 		 	\#\left\{ k\in\{1,\ldots,n\} : g_1(x_k)\neq 0\right\} \le \ell = e \sum_{k=1}^mk!.
 \end{equation}
Moreover, the induced tree analysis $(h_\alpha)_{\alpha\in\mathcal{A}_h}$ of $h$ is such that $w(h_\alpha)\neq w(f_k)$ for all $k=1,\ldots,n$, and therefore, the fact that $(x_k,f_k)$ are strong exact pairs yields
\begin{equation}\label{intro1 eq2}
|h\left(x_k\right)|\le \frac{18}{w(f_k)}, \quad k=1,\ldots,n.
\end{equation}
Considering a further partition of $g_2|_{\supp(x_k)}$, we show that 
\begin{equation}\label{intro1 eq3}
|g_2\left(x_k\right)|\le \frac{18}{w(f_k)}+\frac{3}{2^m}, \quad k=1,\ldots,n.
\end{equation}	
Hence, \eqref{intro1 eq 1}, \eqref{intro1 eq2} and \eqref{intro1 eq3} imply
\[
\left|f\left( \frac{1}{n} \sum_{k=1}^n x_k \right)\right| \le \frac{36+\ell}{n}+\frac{3}{2^m}.
\]	
Then, our choice of $m$ yields Proposition \ref{not contain as l1 intro} for sufficiently large $n$ and $w(f_1)$, where $w(f_1)$ is chosen appropriately to deal with the case where $f\in\{\pm e_j^*:j\in\N\}$.

	Assuming that $Y$ is a $C$-asymptotic $\ell_1$ block subspace of $\mathfrak{X}^{(1)}_{\text{awi}}$, we pick a sequence of standard exact pairs $(x_1,f_1),\ldots,(x_n,f_n)$ in $Y\times W_{(1)}$, satisfying the conclusion of Proposition \ref{not contain as l1 intro} and derive a contradiction.
	
	The remainder of this part of the paper is organized as follows. In Section \ref{asmodel section} we prove that $\mathfrak{X}_{\text{awi}}^{(1)}$ admits $\ell_1$ as a unique asymptotic model. Next, in Section \ref{sep section} we prove existence and properties of standard exact pairs. The final section of this part contains the results leading up to the proof that $\mathfrak{X}^{(1)}_{\text{awi}}$ does not contain Asymptotic $\ell_1$ subspaces.

\section{Asymptotic models generated by block sequences of $\mathfrak{X}_{\textnormal{awi}}^{(1)}$}

\label{asmodel section}

	We show that the space $\mathfrak{X}_{\text{awi}}^{(1)}$ admits a unique asymptotic model, or equivalently, a uniformly unique joint spreading model with respect to $\mathscr{F}_b(\mathfrak{X}_{\text{awi}}^{(1)})$, that is equivalent to the unit vector basis of $\ell_1$. {The key ingredient in the proof is the following lemma concerning bounded positive measures on the tree of initial segments of $\T$.
	}

	{\begin{rem}
		Let us first recall some notation from Section \ref{measures section}. We denote by $\tT$ the tree of initial segments of $\T$ equipped with the partial order induced by inclusion, and consider the isomorphism $t\mapsto\tt =\{s\in\T:s\le_\T t\}$, between $\T$ and $\tT$. Similarly, by $\tW$ we denote the tree of initial segments of $\W$ and  consider the isomorphism $w\mapsto \tw = \{v\in\W:v\le_\W w\}$ between $\W$ and $\tW$. Finally, for $t\in\T$, we set $\tilde{w}(t)=\{\tw\in\tW: w\le_\W w(t) \}$.
	\end{rem}}
	
	\begin{lem}\label{final measure lemma}
		Let $(\mu_i)_i$ be a bounded finitely and disjointly supported sequence in $\mathcal{M}_+(\tT)$. Assume that the sets $\cup\{\supp (f_t):\tt\in\supp(\mu_i)\}$, $i\in\N$, are disjoint. Then, for every $\varepsilon>0$, there exists an infinite subset of the natural numbers $L$ and for each $i\in L$ subsets $G^1_i$, $G^2_i$ of $\tT$  such that
		\begin{itemize}
			\item[(i)] $G^1_i$, $G^2_i$ are disjoint subsets of $\supp(\mu_i)$ for every $i\in L$,
			\item[(ii)] $\mu_i(\tT\setminus G^1_i\cup G^2_i)<\varepsilon$ for every $i\in L$,
			\item[(iii)] $\{t\in\T:\tt\in \cup_{i\in L}G^1_i\}$ is essentially incomparable and
			\item[(iv)] if $F_i^2=\{t\in\T:\tt\in G_i^2\}$, $i\in L$, then the sequence $(F_i^2)_{i\in L}$ is pairwise weight incomparable.
		\end{itemize}
	\end{lem}
	\begin{proof}
 {Passing to a subsequence if necessary, we may assume that the (unique) root of $\tT$ is not in the support of any $\mu_i$, $i\in\N$,
 $\mathrm{succ}\-\!\lim_i\mu_{i}$ exists and that there exist partitions $\mathrm{supp}(\mu_i)=A_i\cup B_i$, $i\in\N$, satisfying the conclusion of Lemma \ref{splitting lemma}.} Define for each $i\in\mathbb{N}$, the measures $\mu_i^1$, $\mu_i^2\in\mathcal{M}_+(\tT)$ given by $\mu_i^1(C) = \mu_i(A_i\cap C)$  and $\mu_i^2(C) = \mu_i(B_i\cap C)$, and let $\nu = w^*\-\lim_i\mu^1_i = \mathrm{succ}\-\!\lim_i\mu^1_i$. Pick a finite subset $F$ of $\tT$ such that $\nu(\tT\setminus F) <\varepsilon/2$. Then, $\nu = w^*\-\lim_i\mu^1_i$ implies that $\lim_i\mu_i^1(\tT) = \nu(\tT)$ and thus, since $\nu =  \mathrm{succ}\-\!\lim_i\mu^1_i$, we have
\begin{equation*}
\lim_i\Big|\mu^1_i(\tT) - \mu^1_i(\cup_{\tilde t\in F}S(\tilde t))\Big| = \Big|\nu(\tT)  - \lim_i\sum_{\tilde t\in F}\mu_i^1(S(\tilde t))\Big| = \nu(\tT\setminus F)<\frac{\varepsilon}{2}.
\end{equation*}
Hence, we can find $i_0\in\mathbb{N}$ such that for all $i\geq i_0$ we have
\begin{equation}
\label{successor is close enough}
\Big|\mu_i(A_i) - \mu_i\Big(A_i\cap\big(\cup_{\tilde t\in F}S(\tilde t)\big)\Big)\Big| = \Big|\mu^1_i(\tT) - \mu^1_i(\cup_{\tilde t\in F}S(\tilde t))\Big| <\frac \varepsilon 2.
\end{equation}
{We set $\Sigma=\sigma(\{t\in\T:\tt\in F\})$ and
\begin{align*}
R = \{ r\in\T :  w(r)\in\Sigma \text{ and there is } s\in\T    
 \text{ with }w(s)\in\Sigma\text{ such that }r <_\T s\} 
\end{align*}
Note that $\Sigma$ and $R$ are finite, since $F$ is finite. Thus, using the fact that  the sets $\cup\{\supp (f_t):\tt\in\supp(\mu_i)\}$ for $i\in\N$ are disjoint, find $i_1\in\N$ with $i_1\geq i_0$ so that
\begin{equation}
\label{avoid all supports}
\cup_{r\in R}\supp(f_r) < \supp (f_t), \text{ for all }\: \tilde t\in\cup_{i\geq i_1}\mathrm{supp}(\mu_i^1)
\end{equation}
For $G_i^1 = A_i\cap(\cup_{\tt\in F}S(\tt))$, \eqref{successor is close enough} implies that $|\mu_i(A_i) - \mu_i(G_i^1)|<\varepsilon/2$, $i\geq i_1$. We will show that $\{t\in\T:\tt\in\cup_{i\geq i_1}G_i^1\}$ is essentially incomparable, i.e., that (iii) is satisfied. To this end, first observe that if $\tt\in\cup_{i\geq i_1}G_i^1$, then $w(t)\in \Sigma$. Let $\tt_1,\tt_2\in\cup_{i\geq i_1}G_i^1$ with $w(t_1)<_\W w(t_2)$. It is immediate that if $t_3\in\T$ is the generator of $w(t_2)$ with $w(t_3)=w(t_1)$, then $t_3\in R$ and hence \eqref{avoid all supports} implies that $f_{t_3}<f_{t_1}$, proving the desired result.

For the remaining part of the proof, recall the root of $\tT$ avoids the supports of all $\mu_i^2$, $i\ge i_1$. This implies that every $\tt\in\cup_{i\ge i_1} \supp(\mu_i^2)$ is the successor of some node in $\tT$. Then, since for all $i\ge i_1$ the set $B_i = \supp(\mu_i^2)$ is finite (as a subset of the finite support of $\mu_i$) and for each $\tt\in\tT$ the sequence $(\mu_i^2(S(\tt)))_{i\ge i_1}$ is eventually zero, we may pass to a subsequence so that for all $i_1\le i<j$ we have $\{w(t):\tt\in \supp(\mu_i^2)\}\cap \{w(t):\tt\in\supp(\mu_j^2)\} = \emptyset$.
We can therefore define the bounded sequence of disjointly supported measures $(\nu_i)_{i\ge i_1}$ on $\widetilde{\mathcal{W}}$ given by 
$\nu_i(\{\tw\}) = \mu_i^2(\{\tilde t\in\tT: \tilde{w}(t) = \tw\})$. Hence, applying Proposition \ref{proposition measures 1 : incomparable} and passing to a subsequence, we obtain a subset $E_i$ of $\supp(\nu_i)$ such that $\nu_i(\tW\setminus E_i)<\varepsilon/2$ and the sets $E_i$, $i\ge i_1$, are pairwise incomparable. It is easy to verify that if $G^2_i=\{\tt\in B_i:\tw(t)\in E_i \}$ and $F^2_i=\{t\in\T:\tt\in G^2_i\}$, $i\ge i_1$, then $(F^2_i)_{i\ge i_1}$ is pairwise weight incomparable and $|\mu_i(B_i) - \mu_i(G_i^2)| = \mu^2_i(\tT\setminus G^2_i)<\varepsilon/2$ for every $i\ge i_1$.}
\end{proof}	
\begin{lem}\label{unique l1 as model lemma 1}
		Let $x\in \mathfrak{X}_{\text{awi}}^{(1)}$, $f\in W_{(1)}$ and a tree analysis $(f_\alpha)_{\alpha\in\mathcal{A}}$ of $f$ such that $f_\alpha(x)\ge 0$ for every $\alpha\in\mathcal{A}$. Let $\varepsilon_1,\ldots,\varepsilon_{h(\mathcal{A})}$ be positive reals and $G_i$ be a subset  of $\{ \alpha\in \mathcal{A}:|\alpha|=i \}$ such that $\sum_{\alpha\in G_i}w_f(f_\alpha)^{-1}f_\alpha(x)>f(x)-\varepsilon_i$ for every $1\le i\le h(\mathcal{A})$, and $f(x)>\sum_{i=1}^{h(\mathcal{A})}\varepsilon_i$. Then, there exists a $g\in W_{(1)}$ satisfying the following conditions.
		\begin{itemize}
			\item[(i)] $\supp (g)\subset \supp (f)$ and $w(g)=w(f)$.
			\item[(ii)]  $g(x)>f(x)-\sum_{i=1}^{h(\mathcal{A})}\varepsilon_i$.
			\item[(iii)] $g$ has a tree analysis $(g_\alpha)_{\alpha\in \mathcal{A}_g}$ such that for every $\alpha\in \mathcal{A}_g$, there is a unique $\beta\in G_{|\alpha|}$ with $\supp (g_\alpha)\subset\supp (f_\beta)$ and $w(g_\alpha)=w(f_\beta)$.
		\end{itemize}
	\end{lem}
	\begin{proof}
		Let $\mathcal{A}_k$ denote the set of all nodes in $\mathcal{A}$ such that $|\alpha|=k$, $1\le k\le h(\mathcal{A})$. {We define $g$ by constructing the tree analysis $(g_\alpha)_{\alpha\in \mathcal{A}_g}$. First, define by induction $B_1=\mathcal{A}_1\setminus G_1$ and for $2\le k \le h(\mathcal{A})$:
		\[
			B_k = \{ \alpha \in \mathcal{A}_k : \alpha\notin G_{k}\text{ or there is a } \beta\in B_{k-1}\text{ such that }  \alpha\in S(\beta) \}.
		\] 
		It follows easily that $\alpha\in B_k$ if and only if there exists $\beta\le\alpha$ such that $\beta\notin G_{|\beta|}$.
		Let $\mathcal{C}_g= \mathcal{A}_{h(\mathcal{A})}\setminus {B}_{h(\mathcal{A})}$. Note that $f_\alpha\in\{\pm e_j^*:j\in\N\}$ for every $\alpha\in\mathcal{C}_g$, and let $\Delta_g=\cup\{\supp (f_\alpha):\alpha\in\mathcal{C}_g \}$. Then $g=f|_{\Delta_g}$ and $(g_\alpha)_{\alpha\in\mathcal{A}_g}$ is the tree analysis induced by $(f_\alpha)_{\alpha\in \mathcal{A}}$. }

		Observe that, by construction, $g$ satisfies (i) and (iii). To see that it also satisfies (ii), we show by induction that for every $1\le k \le h(\mathcal{A})$
		\begin{equation}\label{equation}
			\sum_{  \alpha\in\mathcal{A}_k\setminus {B}_k  }\frac{f_\alpha(x)}{w_f(f_\alpha)} > f(x) - \sum_{i=1}^k \varepsilon_i.
		\end{equation}
		This indeed proves (ii), since the left hand side of (\ref{equation}) for $k=h(\mathcal{A})$ is equal to $g(x)$. We now prove \eqref{equation} by induction. Assume that the inequality holds for some $1\leq k <h(\mathcal{A})$. Then, for every $\alpha\in\mathcal{A}_{k}\setminus {B}_{k}$, we have
		\[
			f_\alpha(x)=\sum_{ \beta\in S(\alpha)\cap G_{k+1} }\frac{f_\beta(x)}{w(f_\alpha)}+\sum_{ \beta\in S(\alpha)\setminus G_{k+1} }\frac{f_\beta(x)}{w(f_\alpha)}
		\] 
		and 
		\[
			\sum_{\alpha\in\mathcal{A}_k\setminus {B}_k}\sum_{\beta\in S(\alpha)\setminus G_{k+1} }\frac{f_\beta(x)}{w_f(f_\alpha)w(f_\alpha)} = \sum_{\alpha\in\mathcal{A}_k\setminus {B}_k}\sum_{\beta\in S(\alpha)\setminus G_{k+1} }\frac{f_\beta(x)}{w_f(f_\beta)}<\varepsilon_{k+1}.
		\]
		Hence
\begin{align*}
\sum_{\alpha\in\mathcal{A}_k\setminus {B}_k}\sum_{\beta\in S(\alpha)\cap G_{k+1} }\frac{f_\beta(x)}{w_f(f_\alpha)w(f_\alpha)}&=
\sum_{\alpha\in\mathcal{A}_k\setminus {B}_k}\frac{f_{\alpha}(x)}{w_{f}(f_{\alpha})}-
  \sum_{\alpha\in\mathcal{A}_k\setminus {B}_k}\sum_{\beta\in S(\alpha)\setminus G_{k+1} }\frac{f_\beta(x)}{w_f(f_{\beta})}\\
  &>( f(x) - \sum_{i=1}^k \varepsilon_i)-\varepsilon_{k+1}
  \end{align*}
which, along with the previous inequality, proves the desired result since 
		\[
			\{ \beta\in\mathcal{A}:\beta\in S(a)\setminus G_{k+1}\text{ for some }\alpha\in \mathcal{A}_k\setminus {B}_k \}	= \mathcal{A}_{k+1}\setminus {B}_{k+1}.
		\]
	\end{proof}	

	\begin{lem}\label{unique l1 as model lemma 2}
		Let $(x^1_j)_j,\ldots,(x^l_j)_j$ be normalized block sequences in $\mathfrak{X}_{\text{awi}}^{(1)}$. For every $\varepsilon>0$, there exists an $L\in[\N]^\infty$ and a $g^i_j\in W_{(1)}$ with $g^i_j(x^i_j)>1-\varepsilon$, $1\le i\le l$ and $j\in L$, such that for any choice of $1\le i_j\le l$ the sequence $(g^{i_j}_j)_{j\in L}$ is AWI.
	\end{lem}
	\begin{proof}
		Let $(\varepsilon_k)_{k=0}^\infty$ be a sequence of positive reals such that $\sum_{k=0}^\infty\varepsilon_k<\varepsilon/2$. For every $1\le i\le l$ and $j\in \N$, pick an $f^i_j\in W_{(1)}$ and a tree analysis $(f^i_{j,\alpha})_{\alpha\in \mathcal{A}^i_j}$ of $f^i_j$ such that $f^i_j(x^i_j)>1-\varepsilon/2$ and $f^i_{j,\alpha}(x^i_j)> 0$ for every $\alpha\in\mathcal{A}^i_j$. For $1\le i\le l$ and $j\in \N$, we set $t_{j}^i=\barf_j^i$ and $t_{j,\alpha}^i=\barf_{j,\alpha}^i$, $\alpha\in\mathcal{A}_j^i$. We will choose, by induction, an $L\in [\N]^\infty$ and, for every $1\le i\le l$, $j\in L$ and $k\in\N$, a subset $G^{k,i}_j$ of $\{ \alpha\in\mathcal{A}^i_j:|\alpha|=k \}$ satisfying the following conditions. For $k\in\N$, we set {$L_{>k}=\{j\in L:j>k\}$}.
		\begin{itemize}
			\item[(i)]  For every $j\in L$, there is a partition $\{ t^{i}_j:i=1,\ldots,l\}=C^{0}_{1,j}\cup C^{0}_{2,j}$ such that $\cup_{j\in L}C^{0}_{1,j}$ is essentially incomparable and $(C^{0}_{2,j})_{j\in L}$ is pairwise weight incomparable.
			\item[(ii)] For every $1\le i\le l$, $k\in\N$ and {$j\in L_{>k}$}, there is a partition $G^{k,i}_j=G^{k,i}_{1,j}\cup G^{k,i}_{2,j}$ such that for any choice of $1\le i_j \le l$, $\cup_{j\in L_{>k} } \{t^{i_j}_{j,\alpha} : \alpha\in G^{k,i_j}_{1,j} \}$ is essentially incomparable and $(\{ t^{i_j}_{j,\alpha} : \alpha\in G^{k,i_j}_{2,j} \})_{j\in L_{>k}}$ is pairwise weight incomparable.
			\item[(iii)] For every $i=1,\ldots,l$, $j\in L$ and $k\in\N$ with $k\le h(\mathcal{A}^i_j)$  \[\sum_{\alpha\in G^{k,i}_{j}} w_{f^i_j}(f^i_{j,\alpha})^{-1}f^i_{j,\alpha}(x^i_j)>f^i_j(x^i_j)-\varepsilon_k.\]
		\end{itemize}
		Observe then that (iii) and an application of the previous lemma yield, for every $1\le i \le l$ and $j\in L$, a functional $g^i_j\in W_{(1)}$ such that 
		\[
			g^i_j(x^i_j)>f^i_j(x^i_j)-\sum_{k=1}^\infty\varepsilon_k>1-\varepsilon.
		\]
		Fix a choice of $1\le i_j\le l$, $j\in L$. Then, (i) implies that $\{ t^{i_j}_j:j\in L\}\cap(\cup_{j\in L}C^0_{1,j})$ is essentially incomparable, and that $\{ t^{i_j}_j:j\in L\}\cap(\cup_{j\in L}C^0_{2,j})$ is weight incomparable. Finally, (ii) and Lemma \ref{unique l1 as model lemma 1} (iii)  yield that, for every $1\le i\le l$, $k\in\N$ and $j\in L_{>k}$, there is a partition 
		\[
			\{\barg^{i_j}_{j,\alpha}:a\in \mathcal{A}^{i_j}_j\text{ and }|\alpha|=k\}=C^{k,i_j}_{1,j}\cup C^{k,i_j}_{2,j}
		\]
		such that $\cup_{j\in L_{>k}}C^{k,i_j}_{1,j}$ is essentially incomparable and $(C^{k,i_j}_{2,j})_{j\in L_{>k}}$ is pairwise weight incomparable. Hence $g^i_j$, $1\le i\le l$ and $j\in L$, satisfy the desired conditions.

		To obtain $L$, let us first assume that $\sup_{i,j}h(\mathcal{A}^i_j)=+\infty$ (if $\sup_{i,j}h(\mathcal{A}^i_j)<+\infty$, then a finite version of the same proof works). {Moreover, passing to a subsequence, we may further assume that $\max_i h(\mathcal{A}_j^i)>k$ whenever $j>k$, for $j,k\in\N$.} Define, for each $j\in\N$, the measure $\mu^0_j$ on $\tT$ given by
		\[
			\mu_j^0=\sum_{i=1}^lf^i_j(x^i_j)\delta_{\tilde{t}_{j}^i }.
		\]
		Applying Lemma \ref{final measure lemma}, we obtain an $L_0\in[\N]^\infty$ such that, for every $j\in L_0$, there exist disjoint subsets $G^0_{1,j}$ and $G^0_{2,j}$ of $\supp (\mu^0_j)$ so that the following hold.
		\begin{itemize}
			\item[($\alpha_0$)] $\mu^0_j( \tT \setminus G^0_{1,j}\cup G^0_{2,j} ) < \varepsilon_0$ for every $j\in L_0$.
			\item[($\beta_0$)] Define $C^0_{1,j}=\{t\in\T:\tt\in G^0_{1,j}\}$ for $j\in L_0$. Then $\cup_{j\in L_0}C^0_{1,j}$ is essentially incomparable.
			\item[($\gamma_0$)] Define $C^0_{2,j}=\{t\in\T:\tt\in G^0_{2,j}\}$ for $j\in L_0$. Then the sequence $(C^0_{2,j})_{j\in L_0}$ is pairwise weight incomparable.
		\end{itemize}
		Note that ($\alpha_0$) implies that $\supp (\mu^0_j)=G^0_{1,j}\cup G^0_{2,j}$ since $f^i_j(x^i_j)>1-\varepsilon/2$, that is, $\{\tt^i_j:i=1,\ldots,l \}=C^0_{1,j}\cup C^0_{2,j}$.
We proceed by induction on $\N$. Suppose we have chosen $L_0,\ldots,L_{k-1}$ and $G^0_{1,j_0},G^0_{2,j_0},\ldots,G^{k-1}_{1,j_{k-1}},G^{k-1}_{2,j_{k-1}}$ for some $k\in \N$ and every $j_i\in L_i$, for $i=0,\ldots,k-1$. Set $L_k^0=\{j\in L_{k-1}:h(\mathcal{A}^i_j)< k\text{ for all }1\le i\le l\}$. Then, for each {$j\in L_{k-1}\setminus L^0_{k}$}, define the following measure on $\tT$
		\[
			\mu^{k}_j=\sum_{i=1}^l\sum_{\substack { \alpha\in\mathcal{A}^i_j \\ |\alpha|=k } } \frac{f^i_{j,\alpha}(x^i_j)}{w_{f^i_j}(f^i_{j,\alpha})}\delta_{\tilde{t}^i_{j,\alpha}}.
		\]
		Again, applying Lemma \ref{final measure lemma} yields an {$L^1_{k}\in [L_{k-1}\setminus L^0_k]^\infty$} and disjoint subsets $G^{k}_{1,j}$ and $G^{k}_{2,j}$ of $\supp (\mu^{k}_j)$, $j\in L^1_k$, such that
		\begin{itemize}
			\item[($\alpha_{k}$)] $\mu^k_j( \tT \setminus G^k_{1,j}\cup G^k_{2,j} ) < \varepsilon_k$ for every $j\in L^1_k$,
			\item[($\beta_{k}$)] $\{t\in\T:\tt\in \cup_{j\in {L^1_k}}G^k_{1,j}\}$ is essentially incomparable and
			\item[($\gamma_{k}$)] the sequence $(\{t:\tt\in G^k_{2,j}\})_{j\in {L^1_k}}$ is pairwise weight incomparable.
		\end{itemize}
		{Then, set $L_k=L_k^0\cup L_k^1$ and $G^k_{i,j}=\{\alpha\in\mathcal{A}_j^i:|\alpha|=k\}$, for $1\le i\le l$ and $j\in L_{k}^0$.}
		Finally, choose $L$ to be a diagonalization of $(L_k)_{k}$, i.e., $L(k)\in L_k$ for  $k\in\N$. Observe that $(\beta_k)$ and $(\gamma_k)$ imply (ii), while $(\alpha_k)$ implies (iii).
	\end{proof}

	\begin{prop}\label{unique l1 asymptotic model}
		The space $\mathfrak{X}_{\text{awi}}^{(1)}$ admits a unique asymptotic model, with respect to $\mathscr{F}_b(\mathfrak{X}_{\text{awi}}^{(1)})$, equivalent to the unit vector basis of $\ell_1$.
	\end{prop}
	\begin{proof}
		Equivalently, we will show that $\mathfrak{X}_{\text{awi}}^{(1)}$ admits $\ell_1$ as a uniformly unique joint spreading model with respect to $\mathscr{F}_b(\mathfrak{X}_{\text{awi}}^{(1)})$. To this end, let $(x^1_j)_j,\ldots,(x^l_j)_j$ be normalized block sequences in $\mathfrak{X}_{\text{awi}}^{(1)}$. Passing to a subsequence, we may assume that $\supp (x^{i_1}_j) < \supp (x^{i_2}_{j+1})$ for every $i_1,i_2=1,\ldots, l$ and $j\in \N$. Fix $\varepsilon>0$ and apply Lemma \ref{unique l1 as model lemma 2} to obtain an $L\in[\N]^\infty$ and a functional $g^i_j\in W_{(1)}$, for each $1\le i\le l$ and $j\in L$, such that
		\begin{itemize}
			\item[(i)] $\supp (g^i_j)\subset \supp (x^i_j)$ and $g^i_j(x^i_j)>1-\varepsilon$, for all $1\le i\le l$ and $j\in L$, and
			\item[(ii)] the sequence $(g^{i_j}_j)_{j\in L}$ is AWI for any choice of $1\le i_j\le l$, $j\in L$.
		\end{itemize}
Fix a choice of $1\le i_j\le l$, $j\in L$, and let $k\in\N$ and $F\subset L$ with $L(k)\le F$ and $|F|\le k$. Note that $(g^{i_j}_{j})_{j\in F}$ is an $\S_1$-admissible sequence in $W_{(1)}$ and is in fact AWI, as implied by (ii). Hence $g=1/2\sum_{j\in F} g^{i_j}_{j}$ is in $W_{(1)}$ and thus, for any choice of scalars $(a_{j})_{j\in F}$, we calculate
		\[
			\left\| \sum_{j\in F} a_{j} x^{i_j}_{j} \right\|  \ge \left\| \sum_{j\in F} |a_{j}| x^{i_j}_{j} \right\| \ge g \left( \sum_{j\in F} |a_{j}| x^{i_j}_{j} \right) 
			 \ge \frac{1-\varepsilon}{2} \sum_{j\in F} |a_{j}|.
		\]	
		Then, Lemma \ref{ell_p as uniformly unique joint spreading equivalent form} yields the desired result.
	\end{proof}

\section{Standard exact pairs}
\label{sep section}

We pass to the study of certain basic properties of Mixed Tsirelson spaces which have appeared in several previous papers (see \cite{AD2} and \cite{AM}). The goal of this section is to define the standard exact pairs in $\mathfrak{X}_{\text{awi}}^{(1)}$ and present their basic properties. In the next section, we will use the existence of sequences of such pairs in any block subspace of $\mathfrak{X}_{\text{awi}}^{(1)}$ to show that it is not Asymptotic $\ell_1$. The proof of the properties of the standard exact pairs are based on the definition of an auxiliary space and the basic inequality which are given in Appendix A.

\subsection{Special Convex Combinations}
We return our attention to special convex combinations, defined in Section \ref{scc subsection}. These types of vectors are used to prove the presence of standard exact pairs in every block subspace of $\mathfrak{X}_{\text{awi}}^{(1)}$.

\begin{rem}\label{remark (n,e)-scc from blocks}
Let $(x_k)_k$ be a block sequence in $\mathfrak{X}_{\text{awi}}^{(1)}$. Then Proposition \ref{repeated averages proposition} implies that, for every $\e>0$, $n,m\in \N$ and $M\in[\N]^\infty$ there exist $F\subset M$ with $m\le F$ and scalars $(a_k)_{k\in F}$ such that $\sum_{k\in F}a_kx_k$ is a $(n,\e)$-s.c.c.
\end{rem}

\begin{lem}
\label{good lower ell1}
Let $(x_k)_k$ be a normalized block sequence in $\mathfrak{X}_{\text{awi}}^{(1)}$. For every $\varepsilon>0$, there exists $M\in[\N]^\infty$ such that for every $j\in\N$, every $\S_{n_j}$-admissible sequence $(x_k)_{k\in F}$ with $F\subset M$ and any choice of scalars $(a_k)_{k\in F}$ we have
\[
\Big\|\sum_{k\in F}a_kx_k\Big\|\ge \frac{1-\varepsilon}{m_j}\sum_{k\in F}|a_k|.
\]
\end{lem}
\begin{proof}
Apply Lemma \ref{unique l1 as model lemma 2} to obtain $M\in[\N]^\infty$ and an $f_k\in W_{(1)}$ with $f_k(x_k)>1-\varepsilon$, for each $k\in M$, such that $(f_k)_{k\in M}$ is AWI. We may also assume that $\supp (f_k)\subset \supp (x_k)$, $k\in M$. Pick an $F\subset M$ such that $(x_k)_{k\in F}$ is $S_{n_j}$-admissible. Then, $(f_k)_{k\in F}$ is $\S_{n_j}$-admissible and clearly $(f_k)_{k\in F}$ is AWI. Hence, $f={m_j}^{-1}\sum_{k\in F}f_k$ is in $W_{(1)}$ and we calculate
\[
\|\sum_{k\in F}a_kx_k\| = \|\sum_{k\in F}|a_k|x_k\|  \ge f(\sum_{k\in F}|a_k|x_k)\ge \frac{1-\varepsilon}{m_j}\sum_{k \in F}|a_k|.
\]
\end{proof}

\begin{prop}
\label{very good ell1 vectors}
Let $Y$ be a block subspace of $\mathfrak{X}_{\text{awi}}^{(1)}$. Then, for every $n\in\N$ and $\e>0$,  there exists a $(n,\e)$-s.c.c. $x = \sum_{k=1}^mc_kx_k$ with $\|x\| > 1/2$, where $x_1,\ldots,x_m$ are in the unit ball of $Y$.
\end{prop}

\begin{proof}
Towards a contradiction assume that the conclusion is false. That is, for any $\mathcal{S}_n$-admissible sequence $(x_k)_{k=1}^m$ in the unit ball of $Y$ such that the vector $x = \sum_{k=1}^mc_kx_k$ is a $(n,\e)$-s.c.c. we have that $\|x\| \leq 1/2$. 

Start with a normalized block sequence $(x^{0}_{k})_{k}$ in $Y$ and
pass to a subsequence satisfying the conclusion of Lemma \ref{good
  lower ell1} {for $\varepsilon=1/2$}. Using  the choice of the sequence $(m_{k})_{k}$, we may find $j\in\N$ such that
\begin{equation}\label{very good ell1 vectors inequality}
2^{n_j/n}\ge 4m_j.
\end{equation}
Set $d=\lfloor n_j/n\rfloor$ and, using Remark \ref{remark (n,e)-scc from blocks}, define inductively block sequences $(x^i_k)_k$, $i=1,\ldots,d$, such that for each $i=1,\ldots,d$ and $k\in \N$ there is an $\mathcal{S}_n$-admissible sequence $(x_m^{i-1})_{m\in F_k^i}$ and coefficients $(c_m^{i})_{m\in F_k^i}$ such that $\tilde x_k^i = \sum_{m\in F_k^i}c_m^{i}x_m^{i-1}$ is a $(n,\e)$-s.c.c. and $x_k^i = 2\tilde x_k^i$.

Using the negation of the desired conclusion, it is straightforward to check by induction that $\|x_k^i\| \leq 1$ for every $i=1,\ldots,d$ and $k\in\N$. Moreover, note that each vector $x_k^i$ can be written in the form
$$x_k^i = 2^i\sum_{m\in G_k^i}d_m^ix_m^0$$
for some subset $G_k^i $ of $\N$ such that $(x_m^0)_{m\in G_k^i}$ is $\mathcal{S}_{ni}$-admissible and  $\sum_{m\in G_k^i}d_m^i = 1$. As the sequence $(x^0_k)_k$ satisfies the conclusion of Lemma \ref{good lower ell1}, we deduce that 
\begin{equation*}
1\geq \|x_1^d\| \geq \frac{2^d}{2m_j}  > \frac{2^{{n_j}/{n}}}{4m_j},
\end{equation*}
since $n_j-n<dn$, and this contradicts \eqref{very good ell1 vectors inequality}.
\end{proof}

\begin{prop}
\label{scc are ris}
Let $x = \sum_{i=1}^mc_ix_i$ be a $(n,\e)$-s.c.c. in $\mathfrak{X}_{\text{awi}}^{(1)}$ with $\|x_i\| \leq  1$, $i=1,\ldots,m$, and $f\in\W_{(1)}$ with $ w(f) = m_j$ such that $n_j<n$. Then we have
\begin{equation*}
 |f(x)| \leq \frac{1+2\e w(f)}{w(f)}.
\end{equation*}
\end{prop}

\begin{proof}
Let $f = m^{-1}_{j}\sum_{l=1}^df_l$, where $(f_l)_{l=1}^d$ is an $\mathcal{S}_{n_j}$-admissible AWI sequence in $W_{(1)}$, and define
\begin{align*}
A = \big\{i\in\{1,\ldots,m\}: & \text{ there is at most one } 1\leq l\leq 	d \\& \text{ such that } \range(x_i)\cap\range(f_l)\neq\emptyset\big\}.
\end{align*}
Note that $|f(x_i)| \leq 1/m_{j}$, for each $i\in A$, and hence
\begin{equation}
\label{scc are ris eq1}
\left|f\left(\sum_{i=1}^mc_ix_i\right)\right| \leq \frac{1}{m_{j}}\sum_{i\in A}c_i + \sum_{i\notin A}c_i.
\end{equation}
Set $B = \{1,\ldots,m\}\setminus A$. The spreading property of the Schreier families implies that the vectors $(x_i)_{i\in B\setminus\{\min(B)\}}$ are $\mathcal{S}_{n_{j}}$-admissible. Moreover, the singleton $\{x_{\min B}\}$ is $\mathcal{S}_0$-admissible. Thus, $\sum_{i\in B\setminus\min B}c_i < \e$ and $c_{\min B}<\e$. Applying this to \eqref{scc are ris eq1} immediately yields the desired conclusion.
\end{proof}

\subsection{Rapidly Increasing Sequences}
These sequences are a standard tool in the study of HI and related constructions. They are the building blocks of standard exact pairs.

\begin{dfn}\label{RIS definition}
	Let $C\ge1$, $I$ be an interval of $\N$ and $(j_k)_{k\in I}$ be a strictly increasing sequence of naturals. A block sequence $(x_k)_{k\in I}$ in $\mathfrak{X}_{\text{awi}}^{(1)}$ is called a $(C,(j_k)_{k\in I})$-rapidly increasing sequence (RIS) if
	\begin{itemize}
		\item[(i)] $\|x_k\|\le C$ for every $k\in I$,
		\item[(ii)] $\max\supp (x_{k-1})\le \sqrt{m_{j_k}}$ for every $k\in I\setminus\{\min I\}$ and
		\item[(iii)] $|f(x_k)|\le C/w(f)$ for every $k\in I$ and $f\in W_{(1)}$ with $w(f)<m_{j_k}$.
	\end{itemize}
\end{dfn}

\begin{prop}
\label{seminormalized RIS}
Let $Y$ be a block subspace of $\mathfrak{X}_{\text{awi}}^{(1)}$ and $C>2$. Then there exists a strictly increasing sequence $(j_k)_{k\in\N}$ of naturals and a $(C,(j_k)_{k\in\N})$-RIS $(x_k)_{k\in\N}$ in $Y$ such that $1/2< \|x_k\| \le 1$, for all $k\in\N$.
\end{prop}

\begin{proof}
We define the sequences $(j_k)_k$ and $(x_k)_k$ inductively as follows. First, choose $x_1$, using Proposition \ref{very good ell1 vectors}, to be a $(n_1,m_{1}^{-2})$-s.c.c. $x_1$ in $Y$ with $1/2< \|x_1\|\le 1$ and set $j_1=1$. Suppose that we have chosen $j_1,\ldots,j_{k-1}$ and $x_1,\ldots,x_{k-1}$ for some $k\in\N$. Then, choose $j_k\in\N$ with $j_k>j_{k-1}$ and $\sqrt{m_{j_k}} > \max\supp(x_{k-1})$ and use Proposition \ref{very good ell1 vectors} to find an $(n_{j_k},m^{-2}_{j_k})$-s.c.c. $x_k$ in $Y$ with $\min\supp(x_k)>\max\supp(x_{k-1})$ and $1/2<\|x_k\| \leq 1$. Proposition \ref{scc are ris} then yields that $x_k$ satisfies (iii) of Definition \ref{RIS definition} and hence we conclude that the sequences $(j_k)_{k\in\N}$ and $(x_k)_{k\in\N}$ satisfy the desired conclusion.
\end{proof}

\subsection{Standard Exact Pairs}
We are ready to define standard exact pairs and prove their existence in every block subspace of  $\mathfrak{X}_{\text{awi}}^{(1)}$.

\begin{dfn}\label{SEP definition}
	Let $C\ge1$ and $j_0\in\N$. We call a pair $(x,f)$, for $x\in \mathfrak{X}_{\text{awi}}^{(1)}$ and $f\in W_{(1)}$, a $(C,m_{j_0})$-standard exact pair (SEP) if there exists a $(C,(j_k)_{k=1}^n)$-RIS $(x_k)_{k=1}^n$ with $j_0<j_1$ such that
	\begin{itemize}
			\item[(i)] $x=m_{j_0}\sum_{k=1}^na_kx_k$ and $\sum_{k=1}^na_kx_k$ is a $(n_{j_0},m_{j_0}^{-2})$-s.c.c.,

		\item[(ii)] $x_k$ is a $(n_{j_k},m_{j_k}^{-2})$-s.c.c. and $1/2< \|x_k\|\le 1$ for every $k=1,\ldots,n$ and
		\item[(iii)] $f=m_{j_0}^{-1}\sum_{k=1}^nf_k$ where $(f_k)_{k=1}^n$ is an  $S_{n_{j_0}}$-admissible AWI sequence in $W_{(1)}$ with {$f_k(x_k)>1/4$}, for every $k=1,\ldots,n$.
	\end{itemize}
\end{dfn}

The following proposition is an immediate consequence of the definition of standard exact pairs, the existence of seminormalized rapidly increasing sequences in every block subspace of $\mathfrak{X}_{\text{awi}}^{(1)}$, as follows from , and  applied to a sequence.

\begin{prop}\label{SEP existence}
	Let $Y$ be a block subspace of $\mathfrak{X}_{\text{awi}}^{(1)}$. Then, for every $C>2$ and $j_0,m\in\N$, there exists a $(C,m_{j_0})$-SEP $(x,f)$ with $x\in Y$ and $m\le \min\supp (x)$.
\end{prop}
\begin{proof}{
	Applying Proposition \ref{seminormalized RIS}, we obtain a $(C,(j_k)_{k\in\N})$-RIS $(x_k)_{k\in\N}$ in $Y$ such that $m\le\minsupp(x_1)$ and $1/2< \|x_k\| \le 1$, $k\in\N$, with $j_0<j_1$. Then, applying Lemma \ref{unique l1 as model lemma 2} for $\varepsilon=1/2$ and passing to a subsequence, we obtain an AWI sequence $(f_k)_{k\in\N}$ in $W_{(1)}$ so that $f_k(x_k)>(1-\varepsilon)/2=1/4$, $k\in\N$. We may assume that $\supp(f_k)\subset\supp(x_k)$, $k\in\N$. Remark \ref{remark (n,e)-scc from blocks} then yields the desired SEP.}
\end{proof}

\begin{dfn}\label{Ifdef}
Let $I$ be an interval of $\N$ and $(x_k)_{k\in I}$ be a block sequence in $\mathfrak{X}_{\text{awi}}^{(1)}$. For every $f\in W_{(1)}$, we define the sets $I_f=\{k\in I:\supp (x_k)\subset\range (f)\}$, {$J_f=I_f\cap \{k\in I:\supp (x_k)\cap\supp (f)\neq\emptyset\}$} and $I'_f=\{k\in I:\supp (x_k)\cap\range (f)\neq\emptyset\}$.

If $(x,f)$ is a $(C,m_{j_0})$-SEP and $g\in W_{(1)}$, then when we write $I_g^x$ {or $J_g^x$} we mean $I_g$ {or $J_g$}, respectively, with respect to the sequence $(x_k)_{k=1}^n$ as in Definition \ref{SEP definition}.
\end{dfn}

\begin{rem}
Let $I$ be an interval of $\N$ and $(x_k)_{k\in I}$ be a block sequence in $\mathfrak{X}_{\text{awi}}^{(1)}$. Then, for every $f\in W_{(1)}$, the following hold.
\begin{itemize}
	\item[(i)] $I_f$ is a finite subset of $I$ and $\# \{k\in I:\supp (x_k)\cap \range (f)\neq \emptyset\} \le \#I_f+2$.
	\item[(ii)] If $f=m_j^{-1}\sum_{l=1}^df_l$, then $\cup_{l=1}^dI_{f_l}\subset I_f$.
	\item[(iii)] If there exists $k\in I$ such that $\range (f)\subsetneq \range (x_k)$, then $I_f=\emptyset$.
\end{itemize}
\end{rem}

\begin{prop}\label{evaluations SEP}
For every $(C,m_{j_0})$-SEP $(x,f)$ the following hold.
\begin{itemize}
	\item[(i)] For every $g\in W_{(1)}$ 
	\[
		\big|g(x)\big|\le 
	\begin{cases}
		\frac{2C}{m_{j_0}},\quad\quad & g=\pm e_i^*\text{ for some }i\in\N\\
			2C[\frac{1}{m_{j_0}}+\frac{m_{j_0}}{w(g)}],\quad\quad & w(g)\ge m_{j_0} \\
		 \frac{6C}{w(g)},\quad\quad & w(g)<m_{j_0}
	\end{cases}
	\]
	\item[(ii)] If $g\in W_{(1)}$ with a tree analysis $(g_\alpha)_{\alpha\in\mathcal{A}}$ such that $I^x_{g_\alpha}=\emptyset$ for all $\alpha\in\mathcal{A}$ with $w(g_\alpha)= m_{j_0}$, then
	\[|g(x)|\le\frac{6C}{m_{j_0}}.\]
	\end{itemize}
\end{prop}

For the proof we refer the reader to Appendix A.

{
\begin{rem}\label{evaluations SEP remark}
	Proposition \ref{evaluations SEP} (ii) remains valid if we replace $I^x_{g_\alpha}$ with $J^x_{g_\alpha}$.
\end{rem}
}

\begin{cor}
The space $\mathfrak{X}_{\text{awi}}^{(1)}$ is reflexive.
\end{cor}
\begin{proof}
The unit vector basis of $c_{00}(\N)$ forms an unconditional Schauder basis for $\mathfrak{X}_{\text{awi}}^{(1)}$ and it is also boundedly complete, since the space admits a unique $\ell_1$ asymptotic model. Hence, it suffices to show that $\mathfrak{X}_{\text{awi}}^{(1)}$ does not contain $\ell_1$.
{
To this end, suppose that $\mathfrak{X}_{\text{awi}}^{(1)}$ contains $\ell_1$ and in particular, from James' $\ell_1$ distortion theorem \cite{James distortion}, there is a normalized block sequence $(x_k)_k$ in $\mathfrak{X}_{\text{awi}}^{(1)}$ such that for $0<\varepsilon<1/2$
\[
\left\| \sum_{k=1}^na_kx_k\right\|\ge (1-\varepsilon)\sum_{k=1}^n|a_k|	
\]
for all $n\in\N$ and any choice of scalars $a_1,\ldots,a_n$. Choose $j_0\in\N$ such that $12/m_{j_0}<1-\varepsilon$. Let also $y_1<\ldots<y_n$, where each $y_i$ is a special convex combination of $(x_k)_k$ for all $i=1,\ldots,n$, such that $x=m_{j_0}\sum_{i=1}^na_ky_k$ is a $(3,m_{j_0})$-SEP (note that $\|y_i\|\ge 1-\varepsilon>1/2$ for all $i=1,\ldots,n$). Then, Proposition \ref{evaluations SEP} yields that $\|x\|\le 12$ and, since $\|x\|=m_{j_0}\|\sum_{i=1}^na_ky_k\|\ge m_{j_0}(1-\varepsilon)$, we derive a contradiction. }
\end{proof}

\section{The space $\mathfrak{X}_{\textnormal{awi}}^{(1)}$ does not contain asymptotic $\ell_1$ subspaces}

In this last section of the first part of the paper we show that $\mathfrak{X}_{\text{awi}}^{(1)}$ does not contain Asymptotic $\ell_1$ subspaces. It is worth pointing out that unlike the constructions in \cite{AGM}, we are not able to prove the existence of a block tree which is either $c_0$ or $\ell_p$, for some $1<p<\infty$, of height greater or equal to $\omega$, in any subspace of $\mathfrak{X}_{\text{awi}}^{(1)}$.

\begin{dfn}
		We say that a sequence $(x_1,f_1),\ldots,(x_n,f_n)$, with $x_i\in\mathfrak{X}_{\text{awi}}^{(1)}$ and $f_i\in W_{(1)}$ for $i=1,\ldots,n$, is a dependent sequence if each pair $(x_i,f_i)$ is a $(3,m_{j_i})$-SEP and $\bar{f}_1<_\T\ldots<_\T \bar{f}_n$.
	\end{dfn}

	\begin{dfn}
		Given a dependent sequence $(x_1,f_1),\ldots,(x_n,f_n)$, for $f\in W_{(1)}$ with a tree analysis $(f_\alpha)_{\alpha\in\mathcal{A}}$ and each $1\le k\le h(\mathcal{A})$ define
		\begin{align*}
			D^k_f = \big\{\alpha\in\mathcal{A}:\:&|\alpha|=k \text{ and there exists } 1\le i\le n\text{ such that}\\ 
			& w(f_\alpha)=w(f_i)\text{ and }\supp (f_\alpha) \cap \range (f_i)\neq \emptyset\big\}	
		\end{align*}
		and 
		\[
			E^k_f=\big\{i\in \{1,\ldots,n \}:w(f_\alpha)=w(f_i)\text{ for some }\alpha\in D^k_f \big\}.	
		\]
	\end{dfn}

	\begin{rem}\label{key observation - upper bound of Ekf}
		Let $f_1,\ldots,f_n,f$ be as in the above definition and fix $k\in\N$. If $f_\alpha$ and $f_\beta$ are such that $\alpha,\beta\in D^k_f$ and $w(f_\alpha)<w(f_\beta)$, then $w(f_\alpha)<_{\W} w(f_\beta)$ since $w(f_\alpha)=w(f_{i_1})$ and $w(f_\beta)=w(f_{i_2})$ for some $1\le i_1<i_2\le n$. This implies that $\{\bar{f}_\alpha,\bar{f}_\beta \}$ is not essentially incomparable. Indeed, if it were essentially incomparable, then $f_{i_1}<f_\alpha$ and this contradicts the fact that $\supp (f_\alpha) \cap \range(f_{i_1})\neq\emptyset$ in the definition of $D^k_f$.
	\end{rem}

	\begin{prop}\label{upper bound of Ekf}
		Let $(x_1,f_1),\ldots,(x_n,f_n)$ be a dependent sequence and $f\in W_{(1)}$. Then $\# E^k_f\le e k!$ for every $k\in\N$ (where $e$ denotes Euler's number).
	\end{prop}
	\begin{proof}
		Denote by $(a_k)_k$ the sequence satisfying the recurrence relation $a_1=2$ and $a_k=ka_{k-1}+1$, $k\ge 2$. We will show that $\# E^k_f\le a_k$ for every $k\in\N$. Note that this yields the desired result since $a_k=\sum_{j=0}^kk!/j!\le e k!$.

		Let $(f_\alpha)_{\alpha\in\mathcal{A}}$ be a tree analysis of $f$. We proceed by induction. For $k=1$, the definition of $W_{(1)}$ and in particular that of AWI sequences, yields a partition 
		\[
			\{ \bar{f}_\alpha:\alpha\in \mathcal{A}\text{ and }|\alpha|=1 \} = C^0_1\cup C^0_2
		\]
		such that $C^0_1$ is essentially incomparable and $C^0_2$ is weight incomparable. Then, note that Remark \ref{key observation - upper bound of Ekf} implies that
		\begin{equation}\label{bound ekf eq1}
		\# \{ w(f_\alpha):\alpha \in D^1_f\text{ and }\bar{f}_\alpha\in C^0_1 \}\le 1.
		\end{equation}
		 Moreover, since $C^0_2$ is weight incomparable, we  have that
		\begin{equation}\label{bound ekf eq2}
			\# \{ w(f_\alpha):\alpha \in D^1_f\text{ and }\bar{f}_\alpha\in C^0_2 \}\le 1	
		\end{equation}
		and hence \eqref{bound ekf eq1} and \eqref{bound ekf eq2} imply that $\#E^1_f\le 2$. 
		
		Assume that for some $k\in\N$ we have $\#E^k_g\le a_k$ for all functionals $g$ in $W_{(1)}$, with respect to the dependent sequence $(x_1,f_1),\ldots,(x_n,f_n)$. We will show that $\#E^{k+1}_f\le a_{k+1}$. Let $\{\alpha\in\mathcal{A}:|\alpha|=1 \}=\{\alpha_1,\ldots,\alpha_d \}$ where $f_{\alpha_1}<\ldots< f_{\alpha_d}$ and consider the tree analyses $(f_\alpha)_{\alpha\in\mathcal{A}_i}$, where $\mathcal{A}_i=\{\alpha\in\mathcal{A}:\alpha_i\le \alpha \}$ for $1\le i \le d$. The fact that $f$ is in $W_{(1)}$, i.e., $(f_{\alpha_i})_{i=1}^d$ is AWI, implies that there exist partitions
		\[
			\{ \bar{f}_\alpha:\alpha\in \mathcal{A}_i\text{ with }|\alpha|=k \} = C^{k}_{1,i}\cup C^{k}_{2,i}, \quad i\ge k+1,
		\]
		such that $\cup_{i=k+1}^d C_{1,i}^{k}$ is essentially incomparable and $(C_{2,i}^{k})_{i=k+1}^d$ is pairwise weight incomparable. {Here $|\alpha|$ is the height of $\alpha$ in the tree $\mathcal{A}_i$.} Then, using Remark \ref{key observation - upper bound of Ekf} and arguing as in the previous paragraph, we have 
		\begin{equation}\label{bound ekf eq3}
			\# \{ w(f_\alpha):\alpha\in D^{k+1}_f\text{ and }\alpha \in \cup_{i=k+1}^d C_{1,i}^{k} \}\le 1.
		\end{equation}
		Moreover, it follows easily that $D^{k+1}_f\cap C^{k}_{2,i_0}\neq\emptyset$ for at most one $k<i_0\le d$ and thus if such an $i_0$ exists we have
		\[
		\#\cup_{i=k+1,i\neq i_0}^d E^{k}_{f_{\alpha_{i}}}\le 1.
		\]
		If no such $i_0$ exists we have
		\[
		\#\cup_{i=k+1}^d E^{k}_{f_{\alpha_{i}}}\le 1.
		\]
		In any case, since the inductive hypothesis yields that $\# E^{k}_{f_{\alpha_{i_0}}}\le a_k$, we have
		\[
		\#\cup_{i=k+1}^d E^{k}_{f_{\alpha_{i}}}\le a_k+1.
		\]
		Note that the inductive hypothesis also implies that 
		\[
		\# E^{k}_{f_{\alpha_{i}}}\le a_k,\quad1\le i\le k
		\]
		and hence, since $E^{k+1}_f=\cup_{i=1}^dE^{k}_{f_{\alpha_i}}$, we  conclude that 
		\[
			\# E^{k+1}_f \le \sum_{i=1}^{k}\# E^{k}_{f_{\alpha_{i}}} + \# \cup_{i=k+1}^d E^{k}_{f_{\alpha_{i}}}\le ka_k+a_k+1	.
		\] 
		This completes the inductive step and the proof. 
	\end{proof} 
	
	{
	\begin{lem}\label{last before last as l1}
		Let $(x,f_0)$ be a $(3,m_{j_0})$-SEP. If $f$ is a functional in $W_{(1)}$ with a tree analysis $(f_\alpha)_{\alpha\in\mathcal{A}}$ and		\[
			\mathcal{B}=\{\alpha\in\mathcal{A}:w(f_\alpha)=m_{j_0}\text{ and }w(f_\beta)\neq m_{j_0}\text{ for every }\beta<\alpha\},
		\]
		then there exists a partition $\ran(f)=G\cup D$ such that
		\begin{itemize}
			\item[(i)] $|f|_D(x)|\le 18/m_{j_0}$ and
			\item[(ii)] $|\sum_{\alpha\in \mathcal{B}}f_\alpha|_G(x)|\le 3$.
		\end{itemize}
	\end{lem}
	\begin{proof}
		Let $x=m_{j_0}\sum_{k=1}^na_kx_k$ for some $(3,(j_k)_{k=1}^n)$-RIS $(x_k)_{k=1}^n$, and set \[I^x_{f_\alpha}=\{k\in\{1,\ldots,n\}:\supp (x_k)\subset\range ({f_\alpha})\}, \quad \alpha\in\mathcal{A}.\] For every $\alpha\in\mathcal{B}$ and every $k\in I^x_{f_\alpha}$, Definition \ref{RIS definition} (iii) implies that
		\begin{equation}\label{lemma SEP decomp eq}
			|f_\alpha(a_kx_k)|\le \frac{3a_k}{m_{j_0}}
		\end{equation}
	since $j_0<j_k$.	Set $G=\cup\{\range(x_k):k\in \cup_{\alpha\in\mathcal{B}}I^x_{f_\alpha}\}$ and $D=\ran(f)\setminus G$. Then, \eqref{lemma SEP decomp eq} immediately yields that $G$ satisfies (ii). To see that $D$ satisfies (i), note that if $\alpha\in\mathcal{A}$ with $w(f_\alpha)=m_{j_0}$ and $\supp(f_\alpha)\cap D=\emptyset$, there exists $\beta\in\mathcal{B}$ such that $\beta \le \alpha$ and $J^x_{f_\alpha|D}\subset J^x_{f_\beta|D}$. However, it is easy to see that $J^x_{f_\beta|D}=\emptyset$ and thus $J^x_{f_\alpha|D}=\emptyset$. Hence (i) follows from Proposition \ref{evaluations SEP} (ii) and Remark \ref{evaluations SEP remark}.
	\end{proof}		
	}

	\begin{prop}\label{last step towards non-Asymtptotic ell 1}
				For every $0<c<1$, there exists $d\in\N$ such that for any dependent sequence $(x_1,f_1),\ldots,(x_n,f_n)$ where $d\le n$, and any $f\in W_{(1)}$ we have
		\[
		\left| f(\frac{1}{n}\sum_{i=1}^nx_i) \right|<c.
		\]
	\end{prop}	
	\begin{proof} First, pick an  $m\in\N$ such that $3/ 2^m<c$ and fix a dependent sequence $(x_1,f_1),\ldots,(x_n,f_n)$. Let $f\in W_{(1)}$ with $(f_\alpha)_{\alpha\in\mathcal{A}}$ be a tree analysis of $f$ and set
		\[
		G=\cup\{\range (x_k)\cap \range (f_\alpha): k\in\{1,\ldots,n\}\text{ and }\alpha\in\mathcal{A}\text{ with }w(f_\alpha)=w(f_k)\}
		\]
and $H=\N\setminus G$. Let $g=f|_G$ and $h=f|_{H}$. Then, consider the tree analysis $(g_\alpha)_{\alpha\in\mathcal{A}_g}$ for $g$, {induced by $(f_\alpha)_{\alpha\in\mathcal{A}}$}, and define
\begin{align*}
\mathcal{B}^1_k=\{\alpha\in\mathcal{A}_g:&\;|a|\le m,\; w(f_\alpha)=w(f_k)\text{ and} \\
&w(f_\beta)\neq w(f_k)\text{ for all }\beta<\alpha\text{ in }\mathcal{A}_g\}
\end{align*}
for $k=1,\ldots,n$ and 
\[
G_1=\cup_{k=1}^n\cup\{\supp(g_\alpha)\cap \supp(x_k):\alpha\in\mathcal{B}^1_k\}.
\]
Let $g_1=g|_{G_1}$ and $g_2=g|_{\N\setminus G_1}$. {Recall Remark \ref{w f f alpha remark} (ii) and observe} that for $k=1,\ldots,n$,
\[
g_1(x_k)=\sum_{\alpha\in\mathcal{B}^1_k}\frac{1}{w_g(g_\alpha)}g_\alpha(x_k)\quad\text{ and }\quad g_2(x_k)=\sum_{\alpha\in\mathcal{B}^2_k}\frac{1}{w_g(g_\alpha)}g_\alpha(x_k)
\]
where
\begin{align}\label{b2k}
\mathcal{B}^2_k=\{\alpha\in\mathcal{A}_g:&\;|a|> m,\; w(f_\alpha)=w(f_k)\text{ and} \nonumber\\
&w(f_\beta)\neq w(f_k)\text{ for all }\beta<\alpha\text{ in }\mathcal{A}_g\}.
\end{align}

  Consider the tree analysis $(h_\alpha)_{\alpha\in \mathcal{A}_h}$ of $h$, {induced by $(f_\alpha)_{\alpha\in\mathcal{A}}$}. Note that, for every $\alpha$ in $\mathcal{A}_h$ and $k=1,\ldots,n$ such that
$w(h_\alpha)=w(f_k)$, we have $\range(h_\alpha)\cap
\range(x_k)=\emptyset$, and hence $k\not\in I_{h_{\alpha}}$.
Proposition \ref{evaluations SEP} (ii) then
implies that for every $k=1,\ldots,n$
\[
 |h(x_k)|\le \frac{18}{w(f_k)}.
\]
Thus, we obtain
\begin{equation}
\label{not As ell 1 eq 1} |h(\frac{1}{n}\sum_{k=1}^nx_k)|\le \frac{18}{n}.
\end{equation}

		Next, we apply Lemma \ref{last before last as l1} for $g_2$ and each $(x_k,f_k)$, $k=1,\ldots,n$, to obtain partitions $\supp(g_2)\cap\supp(x_k)=G^2_k\cup D^2_k$ such that
		 \begin{itemize}
		 	\item[(a)] $|g|_{D^2_k}(x_k)|\le 18/w(f_k)$ and
		 	\item[(b)] $|\sum_{\beta\in\mathcal{B}^2_k}g_\beta|_{G^2_k}(x_k)|\le 3$.
		 \end{itemize}
		 Then, (b) and Remark \ref{w f f alpha remark} (iii) yield that
		 \[
		 |g_2|_{G_k^2}(x_k)|=|\sum_{\beta\in\mathcal{B}^2_k}w_g(g_\beta)^{-1}g_\beta|_{G^2_k}(x_k)|\le \sum_{\beta\in\mathcal{B}^2_k}2^{-m}|g_\beta|_{G^2_k}(x_k)|\le \frac{3}{2^m}
		 \]
		 and hence using (a) we obtain
		 \begin{equation}\label{not As ell 1 eq 2}
			|g_2(\frac{1}{n}\sum_{k=1}^nx_k)|\le \frac{1}{n}\sum_{k=1}^n\frac{18}{w(f_k)}+\frac{3}{2^m}\le \frac{18}{n}+\frac{3}{2^m}.
		 \end{equation}
		 
		 Finally, observe that it follows immediately from Proposition \ref{upper bound of Ekf} that
		 \[
		 	\{ k\in\{1,\ldots,n\} : g_1(x_k)\neq 0\} \le \ell = e\sum_{k=1}^mk!
		 \]
		 and thus, by Proposition \ref{evaluations SEP} (i),
		 \begin{equation}
		 \label{not As ell 1 eq 3} |g_1(\frac{1}{n}\sum_{k=1}^nx_k)| \le \frac{\ell}{n}6.
		 \end{equation}
		 
		 Then for $d$ such that
		 \[
		 \frac{36+6\ell}{d}+\frac{3}{2^m} < c,
		 \]
		 \eqref{not As ell 1 eq 1}, \eqref{not As ell 1 eq 2} and \eqref{not As ell 1 eq 3} yield the desired result.	\end{proof}

	\begin{prop}\label{does not contain as l1 subspaces}
		The space $\mathfrak{X}_{\text{awi}}^{(1)}$ does not contain Asymptotic $\ell_1$ subspaces.
	\end{prop}
	\begin{proof}
		Suppose that $\mathfrak{X}_{\text{awi}}^{(1)}$ contains a $C'$-Asymptotic $\ell_1$ subspace $Y$. By standard arguments, for every $\varepsilon>0$, there exists a block subspace of $Y$ which is $C'+\varepsilon$ Asymptotic $\ell_1$. Passing to a further block subspace, we may assume that $Y$ is block and $C$-asymptotic $\ell_1$ in the sense of \cite{MT}, i.e., $Y$ admits a Schauder basis $(y_i)_i$, which is a block subsequence of $(e_i)_i$, such that for every $n\in\N$ there exists $N(n)\in\N$ with the property that whenever $N(n)\le x_1\le \ldots\le x_n$ are blocks of $(y_i)_i$ then
\begin{equation}\label{Y is a block asymptotic l1 eq}
\frac{1}{C}\sum_{k=1}^n\|x_k\| \le \left\|\sum_{k=1}^nx_k\right\|.
\end{equation}		
Applying Proposition \ref{last step towards non-Asymtptotic ell 1} for $c=1/2C$, we obtain $n\in\N$ such that for any dependent sequence $(x_1,f_1),\ldots,(x_n,f_n)$ we have
\[
\left\|\frac{x_1+\cdots+x_n}{n}\right\|< \frac{1}{2C}.
\]
{We apply Proposition \ref{SEP existence} iteratively to construct a dependent sequence in $Y$ as follows: We find $x_1\in Y$ with $N(n)\le \supp(x_1)\cup\supp(f_1)$, $w(f_1)=\sigma(\bar{0})$ and set $\bar{f}_1=(f_1,\sigma(\bar{0}))$, and for $1<k\le n$, we find $x_k\in Y$ with $w(f_k)=\sigma(\bar{f}_{k-1})$, and set $\bar{f_k}=(f_k,\sigma(\bar{f}_{k-1}))$.} Note that the sequence $(f_k)_{k=1}^n$ is $\S_1$-admissible since $n\le N(n)$. Then, \eqref{Y is a block asymptotic l1 eq} implies that
\[
\left\|\frac{x_1+\cdots+x_n}{n}\right\|\ge \frac{1}{2C}, 
\]
since $\|x_k\|> 1/2$ for each $k=1,\ldots,n$ as follows from Definition \ref{SEP definition}, which is a contradiction.
	\end{proof}

	\begin{quest}
		Let $\xi<\omega_1$ and $1<p\le \infty$. Does there exist a Banach space $X$ with a Schauder basis admitting a unique $\ell_1$ asymptotic model such that any block subspace of $X$ contains an $\ell_p$ (or $c_0$ if $p=\infty$) block tree of height greater or equal to $\omega^\xi$.
	\end{quest}

{\Large
\part{ The case of $\boldsymbol{\ell_p}$ for $\boldsymbol{1<p<\infty}$}
}

\section{Introduction}

In this second part, we treat the case of $1<p<\infty$ and in particular that of $p=2$. The cases where $p\neq 2$ follow as an easy modification. The definition of $\mathfrak{X}^{(2)}_{\text{awi}}$ and the the proofs of its properties are for the most part almost identical to those of $\mathfrak{X}^{(1)}_{\text{awi}}$. We start with the 2-convexification of a Mixed Tsirelson space and define a countably branching well-founded tree on its norming set. Then, employing the notion of asymptotically weakly incomparable constraints we define the norming set $W_{(2)}$ of $\mathfrak{X}^{(2)}_{\text{awi}}$. To prove that the space admits $\ell_2$ as a unique asymptotic model, we use Lemma \ref{ell_p as uniformly unique joint spreading equivalent form} by first applying the combinatorial results of Section 4, in a manner similar to that of Section \ref{asmodel section}, and prove lower $\ell_2$ estimates for arrays of block sequences of $\mathfrak{X}^{(2)}_{\text{awi}}$, by passing to a subsequence. Then, a result similar to \cite[Proposition 2.9]{DM} shows that any block sequence of $\mathfrak{X}^{(2)}_{\text{awi}}$ also has an upper $\ell_2$ estimate. Finally, to prove that $\mathfrak{X}^{(2)}_{\text{awi}}$ does not contain Asymptotic $\ell_2$ subspaces, just like in Part 1, we show that any block subspace contains a vector, that is an $\ell_2$-average of standard exact pairs, with arbitrarily small norm. The existence of standard exact pairs follows again from similar arguments, while the proof that these are strong exact pairs requires a variant of the basic inequality, which we include in Appendix B. In particular, for a block subspace $Y$ and $0<c<1$, we show that there is a sequence of standard exact pairs $(x_1,f_1),\ldots,(x_n,f_n)$ in $Y$ such that $\bar{f}_1\le_\T\ldots\le_\T \bar{f}_n$ and $\|x_1+\cdots+x_n\|<c\:\sqrt{n}$. To prove this, we consider the evaluation of an $f$ in $W_{(2)}$ on such a sequence and partition $f$ into $g+h$ and then $g$ into $g_1+g_2$ as in the proof of Proposition \ref{last step towards non-Asymtptotic ell 1}. An upper bound for $h$ follows from the fact that standard exact pairs are strong exact pairs, while that of $g_1$ is again an immediate consequence of Lemma \ref{upper bound of Ekf}. Finally, for $g_2$, unlike the case of Part 1, we cannot estimate its action on each $x_k$, $k=1,\ldots,n$ using similar arguments. Instead, we need to carefully apply the Cauchy-Schwarz inequality to provide an upper estimate for its action on $x_1+\cdots+x_k$. We demonstrate this in Lemma \ref{lemma for g2}.

\section{The space $\mathfrak{X}_{\textnormal{awi}}^{(2)}$}
\label{ell2section}

Define a pair of strictly increasing sequences of natural numbers $(m_j)_j$, $(n_j)_j$ as follows:
\begin{align*}
	m_1 &= 4 & n_1 &= 1 \\
	m_{j+1}&=m_j^{m_j} & n_{j+1}&=2^{2m_{j+1}}n_j
\end{align*}

	\begin{dfn}
		Let $V_{(2)}$ denote the minimal subset of $c_{00}(\N)$ that
		\begin{itemize}
			\item[(i)] contains $0$ and all $\pm e_j^*$, $j\in \N$ and
			\item[(ii)] whenever $f_1<\ldots<f_n$ is an $\S_{n_j}$-admissible sequence in $V_{(2)}\setminus\{0\}$ for some $j\in\N$ and $\lambda_1,\ldots,\lambda_n\in\Q$ with $\sum_{i=1}^n\lambda_i^2\le 1$, then $m_j^{-1}\sum_{i=1}^n\lambda_if_i $ is in $V_{(2)}$.
		\end{itemize}
	\end{dfn}
	
	The notion of the weight $w(f)$ of a functional $f$ in $V_{(2)}$ is identical to that in Section \ref{section 1}. We also define, in a similar manner, the notion of tree analysis of a functional in $V_{(2)}$, taking into account the $\ell_2$ version of the $(m_j,\S_{n_j})$-operations, in the definition of $V_{(2)}$. Again, it follows from minimality that every $f$ in $V_{(2)}\setminus\{0\}$ admits a tree analysis and finally, for a functional $f$ in $V_{(2)}\setminus\{0\}$ admitting a tree analysis $(f_\alpha)_{\alpha\in\mathcal{A}}$, we define $w_f(f_\alpha)$ as in Definition \ref{weight of functional l1}.
	
	\begin{dfn}
		Let $f\in V_{(2)}$ with a tree analysis $(f_\alpha)_{ \alpha \in \mathcal{A} }$.
		\begin{itemize}
		\item[(i)] Let $\beta\in\mathcal{A}$ with $\beta\neq\emptyset$. Then, if $\alpha\in\mathcal{A}$ is the immediate predecessor of $\beta$, we will denote by $\lambda_\beta$ the coefficient of $f_\beta$ in the  normal form of $f_\alpha$, i.e,
		\[
		f_\alpha =m_j^{-1} \sum_{\beta\in S(\alpha)}\lambda_\beta f_\beta
		\]
		where $S(\alpha)$ denotes the set of immediate successors of $\alpha$ and $w(f_\alpha)=m_j$.
\item[(ii)] For each $\beta\in\mathcal{A}$ we define
		\[
			\lambda_{f,\beta} =
			\begin{cases}
				\prod_{\alpha<\beta}\lambda_\alpha, &\quad \beta\neq\emptyset \\ 
				1 
			\end{cases}
		\]
		\end{itemize}
	\end{dfn}

	\begin{rem}\label{w f f alpha remark 2}	Let $f\in V_{(2)}$ with a tree analysis $(f_\alpha)_{\alpha\in\mathcal{A}}$.
		\begin{itemize}
		\item[(i)] For every $k=1,\ldots,h(\mathcal{A})$ 
		\[
		f=\sum_{|a|=k}\frac{\lambda_{f,\alpha}}{w_f(f_\alpha)}f_\alpha.
		\]
		\item[(ii)] If $\mathcal{B}$ is a maximal pairwise incomparable subset of $\mathcal{A}$, then
		\[
		f=\sum_{\alpha\in\mathcal{B}}\frac{\lambda_{f,\alpha}}{w_f(f_\alpha)}f_\alpha.
		\]
		\item[(iii)] For every $\alpha\in\mathcal{A}$, whose immediate predecessor $\beta$ in $\mathcal{A}$ (if one exists) satisfies $f_\beta\notin \{\pm e_j^*:j\in\N\}$, we have $w_f(f_\alpha)\ge 4^{|\alpha|}.$
		\item[(iv)] If $\mathcal{B}$ is a pairwise incomparable subset of $\mathcal{A}$, then 
		\[
\sum_{\alpha\in\mathcal{B}}\lambda_{f,\alpha}^2\le1.
		\]
\end{itemize}				
		 
	\end{rem}

	Next, as in Section \ref{section 1}, we define a tree $\T$ on the set of all pairs $(f,w(f))$, for $f\in V_{(2)}$ and $w(f)$ a weight of $f$ and consider the trees $\tT$, $\W$ and $\tW$, which are induced by $\T$, and are defined identically to those in Section \ref{section 1}. These are countably branching well-founded trees. Finally, let us recall all three incomparability notions of Definition \ref{incomparabilities 1}, as well as the notion of asymptotically weakly incomparable (AWI) sequences in Definition \ref{incomparabilities 2}.

	\begin{dfn}
		Let $W_{(2)}$ be the smallest subset of $V_{(2)}$ that is symmetric, contains the singletons and whenever $j\in\N$, $f_1<\ldots<f_n$ is an $\S_{n_j}$-admissible AWI sequence in $V_{(2)}$ and $\lambda_1,\ldots,\lambda_n\in\Q$ with $\sum_{i=1}^n \lambda_i^2\le1$, then $m_j^{-1}\sum_{i=1}^n\lambda_if_i\in W_{(2)}$. Denote by $\mathfrak{X}_{\text{awi}}^{(2)}$ completion of $c_{00}(\N)$ with respect to the norm induced by $W_{(2)}$.
	\end{dfn}
	
\begin{rem}\label{remarks on the norming set W l2 version}
\begin{itemize}
\item[(i)] The norming set $W_{(2)}$ can be defined as an increasing union of a sequence $(W^n_{(2)})_{n=0}^\infty$, where $W^0_{(2)}=\{\pm e^*_i:i\in\N\}$ and
\begin{align*}
W^{n+1}_{(2)}=W^n_{(2)}\cup\Big\{\frac{1}{m_j}\sum_{i=1}^m&\lambda_if_i:j,m\in\N,\; (\lambda_i)_{i=1}^m\subset\Q\text{ with }\sum_{i=1}^m\lambda_i^2\le 1\text{ and}\\&(f_i)_{i=1}^m \text{ is an }  \S_{n_j}\text{-admissible AWI sequence in }W^n_{(2)}\Big\}.
\end{align*}
\item[(ii)] Proposition \ref{proposition AWI - subsequence etc} yields that the standard unit vector basis of $c_{00}(\N)$ forms an $1$-unconditional Schauder basis for $\mathfrak{X}_{\text{awi}}^{(2)}$.
		\end{itemize}
	\end{rem}

The following lemma is a result similar to \cite[Proposition 2.9]{DM}, in which we prove upper $\ell_2$ estimates for block sequences of $\mathfrak{X}^{(2)}_{\text{awi}}$.

\begin{prop}\label{upper ell2}
For any block sequence $(x_k)_k$ in $\mathfrak{X}_{\text{awi}}^{(2)}$, any finite subset $F$ of the naturals and $f\in W_{(2)}$ we have
\[
|f(\sum_{k\in F}x_k)| \le 2\sqrt{2}(\sum_{k\in F}\|x_k\|^2)^{\frac{1}{2}}.
\]
\end{prop}
\begin{proof} Recall from Remark \ref{remarks on the norming set W l2 version} that $W_{(2)}=\cup_{n=0}^\infty W^n_{(2)}$. We will show by induction that for every $n\in\N$, every $f\in W^n_{(2)}$ and any finite subset $F$ of $\N$ we have
\[
|f(\sum_{k\in F}x_k)| \le 2\sqrt{2}(\sum_{k\in F}\|x_k\|^2)^{\frac{1}{2}}.
\]
Clearly, this holds for all $f\in W^0_{(2)}$. Hence, let us assume that it also holds for all functionals in $W^n_{(2)}$ for some $n\ge 0$ and fix $f\in W^{n+1}_{(2)}$. Then $f=m_j^{-1}\sum_{=1}^m\lambda_if_i$ where $(f_i)_{i=1}^m$ is an $S_{n_j}$-admissible AWI sequence in $W^n_{(2)}$ and $\lambda_1,\ldots,\lambda_m\in \mathbb{Q}$ with $\sum_{i=1}^m\lambda_i^2\le1$. Define \[I_k=\{i\in\{1,\ldots,m\}:\supp(x_k)\cap\range(f_i)\neq\emptyset\},\quad k\in F,\]
$F_1=\{ k\in F :\#I_k\le1\}$ and $F_2=F\setminus F_1$. We also define \[K_i=\{k\in F_1:\supp(x_k)\cap\range(f_{i})\neq\emptyset\},\quad i=1,\ldots,m.\] Note that if $k\in F_1$, then $k\in K_i$ for at most one $i\in\{1,\ldots,m\}$. Thus, using the inductive hypothesis and the Cauchy-Schwarz inequality, we have
\begin{align}\label{part 2 eq 1}
|f(\sum_{k\in F_1}x_k)|&=m_j^{-1}|\sum_{i=1}^m\lambda_if_i(\sum_{k\in K_i}x_k)| \nonumber\\
&\le\frac{2\sqrt{2}}{m_j}\sum_{i=1}^m|\lambda_i|(\sum_{k\in K_i}\|x_k\|^2)^\frac{1}{2} \nonumber \\
&\le \sqrt{2} (\sum_{i=1}^m\lambda_i^2)^{\frac{1}{2}}(\sum_{i=1}^m\sum_{k\in K_i}\|x_k\|^2)^{\frac{1}{2}} \nonumber \\
&\le \sqrt{2} (\sum_{k\in F_1}\|x_k\|^2)^\frac{1}{2}
\end{align}
Moreover, for each $k\in F_2$, it is easy to see that
\begin{equation}\label{part 2 eq 2}
|m_j^{-1}\sum_{i\in I_k}\lambda_if_i(x_k)|\le (\sum_{i\in I_k}\lambda_i^2)^\frac{1}{2}\|x_k\|.
\end{equation}
Observe that for each $i\in\{1,\ldots,m\}$ there are at most two $k$'s in $F_2$ such that $\supp(x_k)\cap\range(f_{i})\neq\emptyset$ and thus applying the Cauchy-Schwarz inequality and \eqref{part 2 eq 2} we have
\begin{align}\label{part 2 eq 3}
|f(\sum_{k\in F_2}x_k)|&=m_j^{-1}|\sum_{i=1}^m\lambda_if_i(\sum_{k\in F_2}x_k)|=m_j^{-1}|\sum_{k\in F_2}\sum_{i\in I_k}\lambda_if_i(x_k)| \nonumber\\
&\le\sum_{k\in F_2}(\sum_{i\in I_k}\lambda_i^2)^{\frac{1}{2}}\|x_{k}\|\nonumber\\
&\le(\sum_{k\in F_2}\sum_{i\in I_k}\lambda_i^2)^{\frac{1}{2}}(\sum_{k\in F_2}\|x_k\|^2)^\frac{1}{2}\nonumber\\
&\le \sqrt{2}(\sum_{k\in F_2}\|x_k\|^2)^\frac{1}{2}
\end{align}
Finally, \eqref{part 2 eq 1} and \eqref{part 2 eq 3} yield the desired result.
\end{proof}

\section{Asymptotic models generated by block sequences of $\mathfrak{X}_{\textnormal{awi}}^{(2)}$}

In this section we prove that $\mathfrak{X}_{\text{awi}}^{(2)}$ admits $\ell_2$ as a unique asymptotic model. 	This follows as an easy modification of the results of Section \ref{asmodel section}, which yield lower $\ell_2$ estimates, combined with the upper $\ell_2$ estimates of Proposition \ref{upper ell2}. Let us first recall Proposition \ref{final measure lemma}, and note that this in fact holds for the trees defined in the previous section. Applying this we obtain the following variant of Lemma \ref{unique l1 as model lemma 1}, using similar arguments.	

\begin{lem}\label{unique l2 as model lemma 1}
		Let $x\in \mathfrak{X}_{\text{awi}}^{(2)}$, $f\in W_{(2)}$ and a tree analysis $(f_\alpha)_{\alpha\in\mathcal{A}}$ of $f$ such that $f_\alpha(x)>0$ for every $\alpha\in\mathcal{A}$. Let $\varepsilon_1,\ldots,\varepsilon_{h(\mathcal{A})}$ be positive reals and $G_i$ be a subset  of $\{ \alpha\in \mathcal{A}:|\alpha|=i \}$ such that \[\sum_{\alpha\in G_i}\frac{\lambda_{f,\alpha}}{w_f(f_\alpha)}f_\alpha(x)>f(x)-\varepsilon_i\] for $i=1,\ldots, h(\mathcal{A})$, and $f(x)>\sum_{i=1}^{h(\mathcal{A})}\varepsilon_i$. Then, there exists a $g\in W_{(2)}$ such that
		\begin{itemize}
			\item[(i)] $\supp (g)\subset \supp (f)$ and $w(g)=w(f)$.
			\item[(ii)]  $g(x)>f(x)-\sum_{i=1}^{h(\mathcal{A})}\varepsilon_i$.
			\item[(iii)] $g$ has a tree analysis $(g_\alpha)_{\alpha\in \mathcal{A}_g}$ such that, for every $\alpha\in \mathcal{A}_g$ with $|\alpha|=i$, there is a unique $\beta\in G_i$ such that $\supp( g_\alpha)\subset\supp (f_\beta)$ and $w(g_\alpha)=w(f_\beta)$.
		\end{itemize}
	\end{lem}

	\begin{lem}\label{unique l2 as model lemma 2}
		Let $(x^1_j)_j,\ldots,(x^l_j)_j$ be normalized block sequences in $\mathfrak{X}_{\text{awi}}^{(2)}$. For every $\varepsilon>0$, there exists an $L\in[\N]^\infty$ and a $g^i_j\in W_{(2)}$ with $g^i_j(x^i_j)>1-\varepsilon$ for $ i=1,\ldots, l$ and $j\in L$, such that for any choice of $i_j\in\{1,\ldots,l\}$ the sequence $(g^{i_j}_j)_{j\in L}$ is AWI.
	\end{lem}
	\begin{proof}
The proof is similar to that of Proposition \ref{unique l1 as model lemma 2} with $\mu_j^k$ defined as
\[
\mu_j^k= \sum_{i=1}^l\sum_{\substack{\alpha\in\mathcal{A}^i_j\\|\alpha|=k}}\frac{\lambda_{f^i_j,\alpha}f^i_{j,\alpha}(x^i_j)}{w_{f^i_j}(f^i_{j,\alpha})}\delta_{\
{\bar{f}^i_{j,\alpha}}}
\]
and applying Lemma \ref{unique l2 as model lemma 1} instead of \ref{unique l1 as model lemma 1}.
	\end{proof}

	\begin{prop}\label{unique asymptotic model}
		The space $\mathfrak{X}_{\text{awi}}^{(2)}$ admits a unique asymptotic model, with respect to $\mathscr{F}_b(\mathfrak{X}_{\text{awi}}^{(2)})$, equivalent to the unit vector basis of $\ell_2$.
	\end{prop}
	\begin{proof}
	Let $(x^1_j)_j,\ldots,(x^l_j)_j$ be normalized block sequences in $\mathfrak{X}_{\text{awi}}^{(2)}$. Working as in the proof of Proposition \ref{unique l1 asymptotic model} applying Lemma \ref{unique l2 as model lemma 2}, we have that, passing to a subsequence, for any choice of $1\le i_j\le l$ for $j\in\N$, any $F\in \S_1$ and any choice of scalars $(a_j)_{j\in F}$,  there is a functional $g\in W_{(2)}$ with \[g=\frac{1}{4}\sum_{j\in F}\frac{a_j}{(\sum_{j\in F}a_j^2)^{\frac{1}{2}}} g^{i_j}_{j}\]
such that $g^{i_j}_j(x^{i_j}_j)\ge 1-\varepsilon$ and $\supp(g^{i_j}_j)\subset\supp(x^{i_j}_j)$ for $j\in F$. Hence, we calculate
		\begin{equation}\label{lower ell 2 unique asymptotic}
			\left\| \sum_{j\in F} a_{j} x^{i_j}_{m_j} \right\|   \ge g \left( \sum_{j\in F} a_{j} x^{i_j}_{m_j} \right) 
			 \ge \frac{1-\varepsilon}{4} (\sum_{j\in F} a_{j}^2)^{\frac{1}{2}}.
		\end{equation}
		Moreover, Lemma \ref{upper ell2} implies that
		\begin{equation}\label{upper ell 2 unique asymptotic}
		\left\| \sum_{j\in F} a_{j} x^{i_j}_{j} \right\|\le 2\sqrt{2}(\sum_{j\in F} a_{j}^2)^{\frac{1}{2}}
		\end{equation}
Thus \eqref{lower ell 2 unique asymptotic}, \eqref{upper ell 2 unique asymptotic} and Lemma \ref{ell_p as uniformly unique joint spreading equivalent form} yield the desired result.
	\end{proof}

By the above proposition, $\mathfrak{X}_{\text{awi}}^{(2)}$ cannot contain an isomorphic copy of $c_0$ or $\ell_1$. Therefore, by James' theorem \cite{Jrefl} for spaces with an unconditional basis we obtain the following.
	\begin{prop}
The space $\mathfrak{X}_{\text{awi}}^{(2)}$ is reflexive.
\end{prop}

\section{Standard exact pairs}

The definitions of rapidly increasing sequences and standard exact pairs in $\mathfrak{X}^{(2)}_{\text{awi}}$ are almost identical to these in Part 1. We show that standard exact pairs are in fact strong exact pairs. This requires a variant of the basic inequality that we prove in Appendix B. 

\begin{dfn}\label{RIS definition l2	}
	Let $C\ge1$, $I\subset\N$ be an interval and $(j_k)_{k\in I}$ be a strictly increasing sequence of naturals. A block sequence $(x_k)_{k\in I}$ in $\mathfrak{X}_{\text{awi}}^{(2)}$ is called a $(C,(j_k)_{k\in I})$-rapidly increasing sequence (RIS) if
	\begin{itemize}
		\item[(i)] $\|x_k\|\le C$ for every $k\in I$,
		\item[(ii)] $\max\supp (x_{k-1})\le \sqrt{m_{j_k}}$ for every $k\in I\setminus\{\min I\}$ and
		\item[(iii)] $|f(x_k)|\le C/w(f)$ for every $k\in I$ and $f\in W_{(2)}$ with $w(f)<m_{j_k}$.
	\end{itemize}
\end{dfn}

\begin{dfn}\label{SEP definition l2}
	Let $C\ge1$ and $j_0\in\N$. We call a pair $(x,f)$, where $x\in \mathfrak{X}_{\text{awi}}^{(2)}$ and $f\in W_{(2)}$, a $(2,C,m_{j_0})$-standard exact pair (SEP) if there exists a $(C,(j_k)_{k=1}^n)$-RIS $(x_k)_{k=1}^n$ with $j_0<j_1$ such that
	\begin{itemize}
			\item[(i)] $x=m_{j_0}\sum_{k=1}^na_kx_k$ and $\sum_{k=1}^na_kx_k$ is a $(2,n_{j_0},m_{j_0}^{-4})$-s.c.c.,

		\item[(ii)] $x_k$ is a $(2,n_{j_k},m_{j_k}^{-4})$-s.c.c. and $1/2<\|x_k\|\le 1$ for every $k=1,\ldots,n$,
		\item[(iii)] $f=m_{j_0}^{-1}\sum_{k=1}^nf_k$ where $f_k\in W_{(2)}$ with $f_k(x_k)>1/4$ for every $k=1,\ldots,n$
		\item[(iv)] and $48m^2_{j_0}\leq \min\supp(x)$.
	\end{itemize}
\end{dfn}

The proof of the following proposition, which demonstrates the existence of SEPs in any subspace of $\mathfrak{X}_{\text{awi}}^{(2)}$, is similar to that of Proposition \ref{SEP existence} and is omitted.

\begin{prop}\label{SEP existence l2}
	Let $Y$ be a block subspace of $\mathfrak{X}_{\text{awi}}^{(2)}$. Then, for every $C>2$ and $j_0,m\in\N$, there exists a $(2,C,m_{j_0})$-SEP $(x,f)$ with $x\in Y$ and $m\le \min\supp (x)$.
\end{prop}

\begin{prop}\label{sep l2}
For every $(2,C,m_{j_0})$-SEP $(x,f)$ the following hold.
\begin{itemize}
	\item[(i)] For every $g\in W_{(2)}$ 
	\[
		\big|g(x)\big|\le 
	\begin{cases}
			4C[\frac{1}{m_{j_0}}+\frac{m_{j_0}}{w(g)}],\quad\quad & w(g)\ge m_{j_0} \\
		 \frac{12C}{w(g)},\quad\quad & w(g)<m_{j_0}
	\end{cases}
	\]
	\item[(ii)] If $g\in W_{(2)}$ with a tree analysis $(g_\alpha)_{\alpha\in\mathcal{A}}$ such that $I^x_{g_\alpha}=\emptyset$ for all $\alpha\in\mathcal{A}$ with $w(g_\alpha)= m_{j_0}$, then
	\[|g(x)|\le\frac{6C}{m_{j_0}}.\]
	\end{itemize}
\end{prop}
\begin{proof}
We refer the reader to Appendix B.
\end{proof}

\section{The space $\mathfrak{X}_{\textnormal{awi}}^{(2)}$ does not contain asymptotic $\ell_2$ subspaces}

To prove that $\mathfrak{X}_{\textnormal{awi}}^{(2)}$ contains no asymptotic $\ell_2$ subspaces, we use almost identical arguments as in the case of $\mathfrak{X}^{(1)}_{\text{awi}}$. In particular, we show that any block subspace contains a vector, that is an $\ell_2$-average of standard exacts pairs, with arbitrarily small norm. Again, this requires Lemma \ref{upper bound of Ekf}. However, in this case, we employ Lemma \ref{lemma for g2} to carefully calculate certain upper bounds, using the Cauchy-Schwarz inequality.

\begin{dfn}
		We say that a sequence $(x_1,f_1),\ldots,(x_n,f_n)$, where $x_i\in\mathfrak{X}_{\text{awi}}^{(2)}$ and $f_i\in W_{(2)}$ for $i=1,\ldots,n$, is a dependent sequence if each $(x_i,f_i)$ is a $(2,3,m_{j_i})$-SEP and $\bar{f}_1<_\T\ldots<_\T \bar{f}_n$.
	\end{dfn}

\begin{lem}\label{lemma for g2}
Let $(x,f)$ be a $(2,3,m_j)$-SEP where $x=m_j\sum_{k=1}^na_kx_k$ and let $g_1<\cdots<g_m\in W_{(2)}$ with $w(g_i)=m_j$ and $I_{g_i}^x \neq \emptyset$ for all $i=1,\ldots,m$. Then, for any choice of scalars $\lambda_1,\ldots,\lambda_m$, we have
\[
|\sum_{i=1}^m\lambda_ig_i(x)|\le (4\sqrt{2}+1)(\sum_{i=1}^m\lambda_i^2)^{\frac{1}{2}}.
\]
\end{lem}
\begin{proof}
For each $i=1,\ldots,m$, let
\[
g_i=\frac{1}{m_j}\sum_{\ell\in L_i} \lambda_{i\ell} g^i_\ell,\quad \sum_{\ell\in L_i}\lambda_{i\ell}^2\le 1
\]
and define
\[
K_1 = \{ k\in\{1,\ldots,n\} : k\in \cup_{i=1}^m\cup_{\ell\in L_i} I^x_{g^i_\ell} \},\quad K_2=\{1,\ldots,n\}\setminus K_1.
\]
Then, Lemma \ref{upper ell2} and the Cauchy-Schwarz inequality imply that
\begin{align}\label{g2 eq 1}
|\sum_{i=1}^m\lambda_ig_i (m_j\sum_{k\in F_1}a_kx_k)|&=|\sum_{i=1}^m\lambda_im_j^{-1}\sum_{\ell\in L_i}\lambda_{i\ell}g^i_\ell (m_j\sum_{k\in I^x_{g_\ell^i}}a_kx_k)|  \nonumber\\
&\le 2\sqrt{2}|\sum_{i=1}^m\lambda_i\sum_{\ell\in L_i}\lambda_{i\ell}( \sum_{k\in I^x_{g^i_\ell}}a_k^2)^{\frac{1}{2}}| \nonumber \\
& \le 2\sqrt{2} |\sum_{i=1}^m\lambda_i(\sum_{\ell\in L_i}\lambda_{i\ell}^2)^{\frac{1}{2}}( \sum_{k\in \cup_{\ell \in L_i}I^x_{g^i_\ell}}a_k^2)^{\frac{1}{2}}| \nonumber \\
& \le 2\sqrt{2} (\sum_{i=1}^m\lambda_i^2)^{\frac{1}{2}}( \sum_{k\in K_1}a_k^2)^{\frac{1}{2}}.
\end{align}
For each $k=1,\ldots,n$, let
\[
x_k = \sum_{q\in Q_k} b_{kq} y^k_q,\quad \sum_{q\in Q_k}b_{kq}^2\le 1
\]
Define for each $i=1,\ldots,m$ and $\ell \in L_i$
\[
\begin{split}
M^\ell_i  &= \{k\in K_2:\text{ there is }q\in Q_k\text{ with }\supp(y^k_q)\subset\range(g^i_\ell)\}\text{ and for }k\in M_i\\
N_{i\ell}^k &= \{q\in Q_k:\supp(y^k_q)\subset\range(g^i_\ell)\}\\
\end{split}
\]
Also for  $k\in K_2$ define
\[
O_k = \{q\in Q_k:\text{ there are }i\in\{1,\dots,m\}\text{ and }\ell\in L_i\text{ with }q\in N_{i\ell}^k\}.
\]
Finally, also define
\[
F_1=\cup_{i=1}^m\cup_{\ell\in L_i}\cup_{k\in M^\ell_i}\cup_{q\in N_{i\ell}^k} \supp(y^k_q), \quad F_2 = \N\setminus F_1.
\]
Note that the sets $N^k_{i\ell}$, $i\in\{1,\ldots,m\}$, $\ell\in L_i$, $k\in M^\ell_i$, are pairwise disjoint with union $\cup_{k\in K_2}O_k$. Applying Lemma \ref{upper ell2} and the Cauchy-Schwarz inequality once again, we have
\begin{align}\label{g2 eq 2}
|\sum_{i=1}^m\lambda_i g_i |_{F_1}(m_j\sum_{k\in K_2}a_kx_k)| & = |\sum_{i=1}^m\lambda_i \sum_{\ell\in L_i}\lambda_{i\ell} g^i_\ell(\sum_{k\in M^\ell_i}\sum_{q\in N_{i\ell}^k}a_kb_{kq}y^k_q)\nonumber|\\
&\leq 2\sqrt 2\sum_{i=1}^m\lambda_i \sum_{\ell\in L_i}\lambda_{i\ell} (\sum_{k\in M^\ell_i}\sum_{q\in N_{i\ell}^k}a_k^2b^2_{kq})^{1/2}\nonumber\\
&\leq 2\sqrt 2\sum_{i=1}^m\lambda_i (\sum_{\ell\in L_i}\lambda_{i\ell}^2)^{1/2}(\sum_{\ell\in L_i} \sum_{k\in M^\ell_i}\sum_{q\in N_{i\ell}^k}a_k^2b^2_{kq})^{1/2}\nonumber\\
&\leq 2\sqrt 2\sum_{i=1}^m\lambda_i (\sum_{\ell\in L_i} \sum_{k\in M^\ell_i}\sum_{q\in N_{i\ell}^k}a_k^2b^2_{kq})^{1/2} \nonumber\\
&\leq 2\sqrt 2(\sum_{i=1}^m\lambda^2_i)^{1/2}(\sum_{i=1}^m\sum_{\ell\in L_i} \sum_{k\in M^\ell_i}\sum_{q\in N_{i\ell}^k}a_k^2b^2_{kq})^{1/2} \nonumber\\
& = 2\sqrt 2 (\sum_{i=1}^m\lambda^2_i)^{1/2} (\sum_{k\in K_2}a_k^2\sum_{q \in O_k}b^2_{kq})^{1/2}\nonumber\\
&\leq 2\sqrt 2 (\sum_{i=1}^m\lambda^2_i)^{1/2}  (\sum_{k\in K_2}a_k^2)^{1/2}
\end{align}
For each $i=1,\ldots,m$ and $k\in K_2$, define
\begin{align*}
Q^i_k = \{ q\in Q_k:\text{there is an }\ell\in L_i\text{ such that }&\supp(y^k_q)\cap\supp(g^i_\ell)\neq\emptyset \\ &\text{and } \supp(y^k_q)\not\subset\supp(g^i_\ell)\}.
\end{align*}
Observe that, since $(g^i_\ell)_{\ell\in L_i}$ is $\S_{n_j}$-admissible, $(y_q^k)_{Q^i_k}$ is $\S_{n_j+1}$-admissible for all $i=1,\ldots,m$, and Proposition \ref{repeated averages proposition} thus implies that
\[
\sum_{q\in Q^i_k} b_{kq}^2< \frac{3}{\min\supp (x_k)}.
\]
For $i\in\{1,\ldots,m\}$ put $K_2^i=\{k\in K_2: \ran(g_i)\cap\ran x_k\neq\emptyset\}$. The condition $i\in\{1,\ldots,m\}$, $I_{g_i}^x \neq \emptyset$, for $i\in\{1,\ldots,m\}$ implies that each $k\in K_2$ is in at most two sets $K_2^i$.
We then calculate
\begin{align}\label{g2 eq 3}
|\sum_{i=1}^m\lambda_i g_i |_{F_2}(m_j\sum_{k\in K_2}a_kx_k)| & = |\sum_{i=1}^m\lambda_ig_i(m_j\sum_{k\in K_2\cap I_{g_i}}a_k\sum_{q\in Q_k^i}b_{kq}y^k_q)\nonumber| \\
& \le 2\sqrt{2}m_j\sum_{i=1}^m\lambda_i(\sum_{k\in K_2\cap I_{g_i}}a^2_k\sum_{q\in Q_k^i}b^2_{kq})^{1/2}\nonumber\\
& \leq 2\sqrt{2}m_j(\sum_{i=1}^m\lambda_i^2)^{1/2}(\sum_{i=1}^m\sum_{k\in K_2\cap I_{g_i}}a^2_k\frac{3}{\min\supp(x_k)})^{1/2}\nonumber\\
&\leq 4\sqrt{3}m_j(\sum_{i=1}^m\lambda_i^2)^{1/2}(\sum_{k\in K_2}a_k^2\frac{1}{\min\supp(x_k)})^{1/2}\nonumber\\
&\leq (\sum_{i=1}^m\lambda_i^2)^{1/2}\frac{4\sqrt{3}m_j}{\min\supp(x)^{1/2}}\nonumber\\
&\leq (\sum_{i=1}^m\lambda_i^2)^{1/2}\text{ (by Definition \ref{SEP definition l2} (iv))}
\end{align}
Hence, \eqref{g2 eq 1}, \eqref{g2 eq 2} and \eqref{g2 eq 3} yield the desired result.
\end{proof}

	\begin{prop}\label{last step towards non-Asymtptotic ell 2}
				For every $0<c<1$, there exists $d\in\N$ such that whenever $d\le n$ and $(x_1,f_1),\ldots,(x_n,f_n)$ is a dependent sequence  then
		\[
		\|\frac{1}{\sqrt{n}}\sum_{i=1}^nx_i\|<c.
		\]
	\end{prop}	
	\begin{proof} Pick an $m\in\N$ such that 
	\begin{equation}\label{choose m}
	2^{-m+3}<c
	\end{equation} and fix a dependent sequence $(x_1,f_1),\ldots,(x_n,f_n)$. Let $f\in W\setminus W_0$ and consider the partitions $f=h+g$ and $g=g_1+g_2$ as in the proof of Proposition \ref{last step towards non-Asymtptotic ell 1}. Then, the same arguments and Proposition \ref{sep l2} yield that
	\begin{equation}
\label{h eq} |h(\frac{1}{\sqrt{n}}\sum_{k=1}^nx_k)|\le \frac{18}{\sqrt{n}}.
\end{equation}
Moreover, Proposition \ref{upper bound of Ekf} again implies that
		 \[
		 	\#\{ k\in\{1,\ldots,n\} : g_1(x_k)\neq 0\} \le \ell = e\sum_{k=1}^mk!
		 \]
		 and thus, by Proposition \ref{upper ell2} and Proposition \ref{sep l2} (i),
		 \begin{equation}
		 \label{g1 eq} |g_1(\frac{1}{\sqrt{n}}\sum_{k=1}^nx_k)| \le 2\sqrt{2} \frac{\sqrt \ell}{\sqrt{n}}24 = 48\sqrt{\frac{2\ell}{n}}.
		 \end{equation}

 		 Finally, we treat $g_2$ differently from Proposition \ref{last step towards non-Asymtptotic ell 1}. Recall that for $k=1,\ldots,n$,
		 \[\mathcal{B}_k^2 = \{\alpha\in\mathcal{A}_{f}:|\al| > m, w(f_{\alpha}) = w(f_k),\text{ and } w(f_{\beta})\neq w(f_k)\text{ for }\beta<\alpha\text{ in }\mathcal{A}_{f}\}.\]
		 Define
		 \[G_2 = \cup_{k=1}^n\cup\{\ran(x_k)\cap\ran(f_{\alpha}):\alpha\in\mathcal{B}_k^2\},\]
		 so that $g_2 = g|_{G_2}$. We further split $G_2$ as follows
		 \[G_2^1 = \cup_{k=1}^n\cup\{\supp(x_k)\cap\supp(f_{\alpha}):\alpha\in\mathcal{B}_k^2\text{ and }I^{x_k}_{f_\alpha} = \emptyset\}\text{ and }G_2^2 = G_2\setminus G_2^1.\]
Proposition \ref{sep l2} (ii) implies that for $k\in\{1,\ldots,n\}$,
\begin{equation*}
|g_2|_{G_2^1}(x_k)| \leq \frac{18}{w(f_k)}
\end{equation*}
and thus
		\begin{equation}\label{g21 eq}
		|g_2|_{G_2^1}(\frac{1}{\sqrt{n}}\sum_{k\in K_1}x_k)|\le \frac{18}{\sqrt{n}}.
		\end{equation}	
To complete the computation, we need to evaluate the action of $g_2|_{G_2^2}$. To that end, for $s=m+1,m+2,\ldots$ and $k\in\{1,\ldots,n\}$ put
\[\mathcal{B}_{k,s}^2 = \{\alpha\in\mathcal{B}_k^2:|\al| = s\},\]
	 so that for each $s>m$, the sets $\mathcal{B}_{k,s}^2$, $k\in\{1,\ldots,n\}$ are pairwise disjoint and the set $\cup_{k=1}^n\mathcal{B}_{k,s}^2$ is pairwise incomparable. We use Lemma \ref{lemma for g2} and the definition of $G_2^2$ to calculate
	 \begin{align}
	 \label{g22 eq}
	 |g_2|_{G_2^2}(\frac{1}{\sqrt n}\sum_{k=1}^nx_k)| &= |\frac{1}{\sqrt{n}}\sum_{k=1}^n\sum_{\alpha\in\mathcal{B}^2_k}\frac{\lambda_{f_\alpha}}{w_{f}(f_{\alpha})}f_{\alpha}|_{G^2_2}(x_{k})\nonumber|\\
	 &\leq \frac{4\sqrt 2+1}{\sqrt{n}}\sum_{k=1}^n(\sum_{\alpha\in\mathcal{B}^2_k}\frac{\lambda^2_{f_\alpha}}{w_{f}(f_{\alpha})^2})^{1/2}\nonumber\\
	 &\leq (4\sqrt 2+1)(\sum_{k=1}^n\sum_{\al\in\mathcal{B}_{k}^2}\frac{\lambda^2_{f_\alpha}}{w_{f}(f_{\alpha})^2})^{1/2}\nonumber\\
	 &\leq(4\sqrt 2+1)(\sum_{s=m+1}^\infty\frac{1}{4^{s}}\sum_{\al\in\cup_{k=1}^n\mathcal{B}_{k,s}^2}{\lambda^2_{f_\alpha}})^{1/2}\nonumber\\
	 &\leq(4\sqrt{2} + 1)(\sum_{s=m+1}^\infty\frac{1}{4^s})^{1/2} = \frac{4\sqrt{2}+1}{2^m\sqrt3} \leq \frac{4}{2^m}.
	 \end{align}
	 Then, \eqref{h eq}, \eqref{g1 eq}, \eqref{g21 eq}, \eqref{g22 eq} yield that
		\[
		|f(x)|\le \frac{36+48\sqrt{2\ell}}{\sqrt n} + \frac{4}{2^m} \leq \frac{36+48\sqrt{2\ell}}{\sqrt n} + \frac{c}{2}
		\]
		 and thus for $d$ such that
		 \[
		 \frac{36+48\sqrt{2\ell}}{\sqrt d} < \frac{c}{2},
		 \]
		 we have the desired result.
	\end{proof}

	\begin{prop}
		The space $\mathfrak{X}_{\text{awi}}^{(2)}$ does not contain Asymptotic $\ell_2$ subspaces.
	\end{prop}
	\begin{proof}
It is an immediate consequence of Proposition \ref{last step towards non-Asymtptotic ell 2}, using similar arguments as in Proposition \ref{does not contain as l1 subspaces}.
	\end{proof}

	\begin{rem}
	Unlike the case of $\ell_1$, for every $1<p<\infty$, it is in fact possible to define a reflexive Banach space with a Schauder basis, admitting a unique $\ell_p$ asymptotic model with respect to the family of normalized block sequences, whose any block subspace contains an $\ell_1$ block tree of height $\omega^\xi$. Such a space can be defined using the attractors method, which was first introduced in \cite{AAT} and later used in \cite{AMP}.
	\end{rem}

\section{Appendix A}

In this section we prove the properties of standard exact pairs in $\mathfrak{X}_{\text{awi}}^{(1)}$, given in Proposition \ref{evaluations SEP}. This requires three steps. First, we need to define an auxiliary space which is also a Mixed Tsirelson space. Then, on the special convex combinations of its basis, we give upper bounds on the evaluations of the functionals in its norming set $W_{\text{aux}}^{(1)}$. Finally for a standard exact pair $(x,f)$, via the basic inequality, we reduce the upper bounds of the evaluations of functionals in $W_{(1)}$ acting on $x$, to the corresponding one of a functional $g$ in $W_{\text{aux}}^{(1)}$ on a normalized special convex combination of the basis of the auxiliary space.

\subsection{The Auxiliary Space}

\begin{dfn}
	Let $W_{\text{aux}}^{(1)}$ be the minimal subset of $c_{00}(\N)$ such that
	\begin{itemize}
		\item[(i)] $\pm e_i$ is in $W_{\text{aux}}^{(1)}$ for all $i\in\N$ and
		\item[(ii)] for every $j\in\N$ and every $\S_{n_j+1}$-admissible sequence of functionals $(f_i)_{i=1}^d$ in $W_{\text{aux}}^{(1)}$ we have that $f=m_{j}^{-1}\sum_{i=1}^df_i$ is in $W_{\text{aux}}^{(1)}$.
	\end{itemize}
\end{dfn}

The purpose of the following two lemmas is to provide upper bounds for the norms of linear combinations of certain vectors in the auxiliary space.

\begin{lem}\label{evaluations on the basis of auxiliary space}
	Let $j\in\N$ and $\varepsilon>0$ with $\varepsilon\le m_{j}^{-1}$. For every $(n_{j},\varepsilon)$-basic s.c.c. $x=\sum_{k\in F}c_ke_k$   the following hold.
	\begin{itemize}
	\item[(i)] For every $f\in W_{\text{aux}}^{(1)}$ 
	\[
		\big|f(x)\big|\le 
	\begin{cases}
		\varepsilon,\quad\quad & f=\pm e_i^*\text{ for some }i\in\N\\
		  \frac{1}{w(f)},\quad\quad & w(f)\ge m_{j}\\
		 \frac{2}{w(f)m_{j}},\quad\quad & w(f)<m_{j}
	\end{cases}
	\]
	\item[(ii)] If $f\in W_{\text{aux}}^{(1)}$ with a tree analysis $(f_\alpha)_{\alpha\in\mathcal{A}}$ such that $w(f_\alpha)\neq m_{j}$ for all $\alpha\in \mathcal{A}$ and $\varepsilon<m_j^{-2}$ then $|f(x)|<2m_j^{-2}$.
	\end{itemize}
\end{lem}
\begin{proof} We may assume that $\supp (f)\subset F$ and $f(e_i)\ge 0$ for every $i\in\N$. If $f=\pm e_i^*$ for some $i\in F$ then $|f(x)|=c_i<\varepsilon$, since $x$ is an $(n_{j},\varepsilon)$-basic s.c.c. and $\{i\}\in \S_0$.  

Suppose that $m_{j}\le w(f)$. Then $\|f\|_{\infty}\le1/w(f)$ and hence
\[
\big|f(x)\big|\le \|f\|_{\infty}\|x\|_1\le \frac{1}{w(f)}.
\]
In the case where $w(f)=m_i<m_{j}$, let $f=m^{-1}_{i}\sum_{l=1}^df_l$ with $(f_l)_{l=1}^d$ an $\S_{n_i+1}$-admissible sequence in $W_{\text{aux}}^{(1)}$. For $l=1,\ldots,d$, define $D_l=\{k\in F:f_l(e_k)>m^{-1}_{j}\}$ and $D=\cup_{l=1}^dD_l$.  Then, \cite[Lemma 3.16]{AT} implies that $D_l\in \S_{(\log_2(m_{j})-1)(n_{j-1}+1)}$ for each $l=1,\ldots,d$ and hence, since $(f_l)_{l=1}^d$ is $S_{n_{j-1}+1}$-admissible (recall that $i<j$ since $m_i<m_j$) and $D_l\subset\supp (f_l)$, $l=1,\ldots,d$, we conclude that the sequence $(D_l)_{l=1}^d$ is $S_{n_{j-1}+1}$-admissible and 
\[
D=\cup_{l=1}^dD_l\in S_{n_{j-1}+1}*\S_{(\log_2(m_{j})-1)(n_{j-1}+1)}= \S_{\log_2(m_{j})(n_{j-1}+1)}
\]
Since $x$ is an $(n_j,\varepsilon)$-basic s.c.c. and $\log_2(m_{j})(n_{j-1}+1)<n_j$, the above implies that $\sum_{k\in D}c_k<\varepsilon$ and thus
\begin{align*}
	f(x)&=\frac{1}{m_i}\sum_{l=1}^df_l(\sum_{k\in F}c_ke_k)=\frac{1}{m_i}(\sum_{l=1}^df_l|_D(\sum_{k\in F}c_ke_k)+\sum_{l=1}^d f_l|_{\N\setminus D}(\sum_{k\in F}c_ke_k))\\
	&\le	\frac{1}{m_i}(\sum_{k\in D}c_k+\frac{1}{m_{j}})\le \frac{1}{m_i}(\e+\frac{1}{m_j})\le \frac{2}{m_im_{j}}.
\end{align*}

Finally, if there is a tree analysis $(f_\alpha)_{\alpha\in\mathcal{A}}$ of $f$ with $w(f_\alpha)\neq m_j$ for every $\alpha\in\mathcal{A}$, \cite[Lemma 3.16]{AT} implies that $D=\{k\in F:f(e_k)>m_j^{-2}\}\in \S_{(2\log_2(m_j)-1)(n_{j-1}-1)}$ and since $(2\log_2(m_j)-1)(n_{j-1}-1)<n_j$ we have that $\sum_{k\in D}c_i<\e$. Hence, we conclude that
\[
f(x)=\sum_{k\in D}c_kf(x_k)+\sum_{k\in F\setminus D}c_kf(x_k)\le \e + \frac{1}{m^2_j}<\frac{2}{m_j^2}.
\]
\end{proof}

\subsection{The Basic Inequality}

\begin{prop}[Basic Inequality]\label{basic inequality}
Let $(x_k)_{k\in I}$ be a $(C,(j_k)_{k\in I})$-RIS in $\mathfrak{X}_{\text{awi}}^{(1)}$ with $4\le \min\supp (x_{\min I})$, $(a_k)_{k\in I}$ be a sequence of nonzero scalars and $f\in W_{(1)}$ with $I_f\neq \emptyset$. Define $t_k=\max\supp (x_k)$, $k\in I$. Then there exist  
\begin{itemize}
\item[(i)] $g\in W_{\text{aux}}^{(1)}\cup\{0\}$ with $w(g)=w(f)$ if $g\neq 0$ and $\{k:t_k\in\supp (g)\}\subset I_f$,
\item[(ii)] $h\in\{\text{sign}(a_k) e^*_{t_k}:k\in I_f\}\cup\{0\}$ with $k_0\in I_f$ and $k_0 < \min\supp(g)$ if $h=\text{sign}(a_{k_0})e^*_{t_{k_0}}$ and
\item[(iii)] $j_0\ge \min \{j_k:k\in I_f\}$
\end{itemize}
such that
\begin{equation}\label{basic inequality eq}
|f(\sum_{k\in I_f}a_kx_k)|\le C(1+\frac{1}{\sqrt{m_{j_{0}}}})[h+g(\sum_{k\in I_f}a_ke_{t_k})].
\end{equation}
\end{prop}
\begin{proof}
Recall that $W_{(1)}$ is the increasing union of the sequence $(W^n_{(1)})_{n=0}^\infty$ defined in Remark \ref{remarks on the norming set W_{(1)}}. We prove the statement by induction on $n=0,1,\ldots$ for every $f\in W^n_{(1)}$ and every RIS.

For $n=0$ and $f\in W^0_{(1)}$, the fact that $I_f\neq\emptyset$ implies that $I_f=\{k_0\}$, i.e., $f= e^*_{t_{k_0}}$ or $f=-e^*_{t_{k_0}}$ for some $k_0\in I$. In either case, it is immediate to check that $h=\text{sign}(a_{k_0})e^*_{t_{k_0}}$, $g=0$ and $j_0=j_{k_0}$ are as desired.

Fix $n\in\N$ and assume that the conclusion holds for every $f\in W^n_{(1)}$ and every RIS. Pick an $f\in W^{n+1}_{(1)}$ with $f={m_i}^{-1}\sum_{l=1}^df_l$, where $(f_l)_{l=1}^d$ is an $\S_{n_i}$-admissible sequence in $W^n_{(1)}$. We will first treat the two extreme cases, namely, the cases where $i\ge \max \{j_k:k\in I_f\}$ and $i<\min\{j_k:k\in I_f\}$.

For the first case, set $k_0=\max I_f$ and $j_0=j_{k_0}$ and choose $k_1\in I_f$ that maximizes the quantity $|a_{k}|$ for $k\in I_f$. Then, since $(x_k)_{k\in I}$ is a RIS, items (i) and (ii) of Definition~\ref{RIS definition} yield that
\begin{align*}
|f(\sum_{k\in I_f\setminus\{k_0\}}a_kx_k)|&\le \frac{1}{m_i}\max\supp (x_{k_0-1})\|\sum_{k\in I_f\setminus\{k_0\}}a_kx_k\|_\infty\\ &\le \frac{ \max\supp (x_{k_0-1}) }{ m_{j_{k_0}} }C|a_{k_1}|\le\frac{C}{\sqrt{m_{j_{k_0}}}}|a_{k_1}|\nonumber
\end{align*}
and thus
\begin{align}\label{basic 2}
|f(\sum_{k\in I_f}a_kx_k)|&\le \frac{C}{\sqrt{m_{j_{k_0}}}}|a_{k_1}| + |f(a_{k_0}x_{k_0})|\\
&\le \frac{C}{\sqrt{m_{j_{k_0}}}}|a_{k_1}|+C|a_{k_1}|=C(1+\frac{1}{\sqrt{m_{j_{k_0}}}})|a_{k_1}|\nonumber\\
&= C(1+\frac{1}{\sqrt{m_{j_{k_0}}}})\text{sign}(a_{k_1})e^*_{t_{k_1}}(\sum_{k\in I_f}a_ke_{t_k}).\nonumber
\end{align}
That is, $h=\text{sign}(a_{k_1})e^*_{t_{k_1}}$, $g=0$ and $j_{k_0}$ yield the conclusion. 

For the second case, the inductive hypothesis implies that, for every $l=1,\ldots,d$ with $I_{f_l}\neq\emptyset$, there are $g_l$, $h_l$ and $j_{0,l}$ as in (i) - (iii) of the statement, that satisfy the conclusion for the functional $f_l$. Define $J_f=\{k\in I_f: f(x_k)\neq 0\}\setminus \cup_{l=1}^dI_{f_l}$. Then, for every $k\in J_f$, Definition \ref{RIS definition} (iii) yields that
\begin{align*}
|f(a_{k}x_{k})|\le \frac{C}{m_i}|a_{k}|=\frac{C}{m_i}\text{sign}(a_{k})e^*_{t_{k}}(\sum_{k\in I_f}a_{k}e_{t_{k}}).
\end{align*}
and hence we calculate
\begin{align}\label{basic 3}
|f(\sum_{k\in I_f}a_kx_k)|&\le |f(\sum_{k\in \cup_{l=1}^dI_{f_l}}a_kx_k)|+|f(\sum_{k\in J_f}a_kx_k)|\\
&\le\frac{C}{m_i}\sum_{k\in J_f}\text{sign}(a_{k})e^*_{t_{k}}(\sum_{k\in I_f}a_{k}e_{t_{k}})+\frac{C}{m_i}\sum_{l=1}^d[(1+\frac{1}{\sqrt{m_{j_{0,l}}}})(h_l+g_l)](\sum_{k\in I_{f_l}}a_ke_{t_k})\nonumber\\
&\le C(1+\frac{1}{\sqrt{m_{j_{\min I_f}}}})[\frac{1}{m_i}(\sum_{k\in J_f}\text{sign}(a_{k})e^*_{t_{k}}+\sum_{l=1}^dh_l+g_l)](\sum_{k\in I_f}a_ke_{t_{k}}).\nonumber
\end{align}
Define 
\[
g=\frac{1}{m_i}(\sum_{k\in J_f}\text{sign}(a_{k})e^*_{t_{k}}+\sum_{l=1}^dh_l+g_l).
\]
Moreover, for each $l=1,\ldots,d$, define
\[
K_l=\big\{k\in J_f:\min\{ l'=1,\ldots,d:\supp (x_k)\cap \range (f_{l'})\neq\emptyset\}=l\big\}
\]
and
\[
I_l=\{t_k:k\in K_l\}\cup\{\supp (h_l)\}\cup\{ \min\supp (g_l) \}.
\]
Let us make the following remarks. First, observe that $\#K_l\le 2$. In particular, consider the case where $K_l=\{k_1,k_2\}$ for some $l=1,\ldots,d$. Then, $k_1<\min I_{f_l}\le \max I_{f_l}< k_2$ and since $\supp (h_l)\cup \supp (g_l)$ is a subset of  $\{t_k:k\in I_{f_l}\}$ we have $t_{k_1}<\supp (h_l)<\supp (g_l)<t_{k_2}$. Moreover,  if $l<d$ and $\range (f_{l+1})\cap \supp (x_{k_2})\neq\emptyset$, then $k_2\notin K_{l+1}$ and clearly $k_2<I_{l+1}$. In the case where $K_l$ is a singleton for some $l=1,\ldots,d$, then either $\supp (h_l)<\supp (g_l)<k$ or $k<\supp (h_l)<\supp (g_l)$ holds for $K_l=\{k\}$. Hence we conclude that $I_1<\cdots<I_d$. Moreover, let us finally note that  $\min \supp (f_l)\le I_l$ and $\#I_l\le 4$ for every $l=1,\ldots,d$. For each $l=1,\ldots,d$, let $K_l=\{k_1^l,k^l_2\}$, where $k^l_2$ or $k^l_2$ can be ommited if necessary. Then,
\begin{equation}\label{basic 4}
g=\frac{1}{m_i}\sum_{l=1}^d\Big(\text{sign}(a_{k^l_1})e^*_{t_{k^l_1}}+h_l+g_l+\text{sign}(a_{k^l_2})e^*_{t_{k^l_2}}\Big).
\end{equation}
We will show that the sequence $(e^*_{t_k})_{k\in J_f} {}^\frown (h_l)_{l=1}^d  {}^\frown  (g_l)_{l=1}^d$ is $\S_{n_i+1}$-admissible, when the functionals ordered as implied by \eqref{basic 4}, i.e., according to the minimum of their supports. This yields that $g\in W_{\text{aux}}^{(1)}$ and thus $h=0$, $g$ and $j_0=j_{\min I_f}$ satisfy the conclusion, as follows from (\ref{basic 3}). More specifically, we will show that $\cup_{l=1}^dI_l\in\S_{n_i+1}$.   To this end, note that $(I_l)_{l=1}^d$ is $\S_{n_i}$-admissible, since $(f_l)_{l=1}^d$ is $\S_{n_i}$-admissible, $I_1<\cdots <I_d$ and $\min\supp (f_l)\le I_l$ for every $l=1,\ldots,d$. Thus $\cup_{l=1}^dI_l\in \S_{n_i}*\mathcal{A}_4$, since $\#I_l\le 4$ for all $l=1,\ldots,d$. Using item (ii) of Lemma \ref{lemma sn conv am} and the fact that $4\le\min\supp (x_{\min I})$ we conclude that $\cup_{l=1}^dI_l\in \S_{n_i+1}$.

Finally, in the remaining case where $\min\{j_k:k\in I_f\}\le i<\max\{j_k:k\in I_f\}$, define $I^1_f=\{k\in I_f:j_k\le i\}$ and $I^2_f=\{k\in I_f:j_k>i\}$ and observe that $I_f=I_f^1\cup I^2_f$, $\max\{j_k:k\in I^1_f\}\le i$ and $i<\min\{j_k:k\in I^2_f\}$. Applying the result of the first case for $(x_k)_{k\in I^1_f}$ and that of the second for $(x_k)_{k\in I^2_f}$ we have
\begin{align}\label{basic 5}
|f(\sum_{k\in I_f}a_kx_k)|&\le |f(\sum_{k\in I^1_f}a_kx_k)|+|f(\sum_{k\in I^2_f}a_kx_k)|\\
&\le C(1+\frac{1}{\sqrt{m_{j_{\max I^1_f}}}})h(\sum_{k\in I_f}a_ke_{t_{k}})+|f(\sum_{k\in I^2_f}a_kx_k)|\nonumber\\
&\le C(1+\frac{1}{\sqrt{m_{j_{\max I^1_f}}}})h(\sum_{k\in I_f}a_ke_{t_{k}})+ C(1+\frac{1}{\sqrt{m_{j_{\min I^2_f}}}})g(\sum_{k\in I_f}a_ke_{t_k})\nonumber\\
&\le C(1+\frac{1}{\sqrt{m_{j_{\max I^1_f}}}})[h+g(\sum_{k\in I_f}a_ke_{t_{k}})]\nonumber
\end{align}
where $h=\text{sign}(a_{k_1})e^*_{t_{k_1}}$, $k_1\in I^1_f$ maximizes the quantity $|a_k|$ for $k\in I^1_f$  and $g\in W_{\text{aux}}^{(1)}$ with $w(g)=w(f)$.
\end{proof}
\begin{rem}\label{remark basic inequality}
Let $(x_k)_{k\in I}$ and $f$ be as in the statement of Proposition \ref{basic inequality}.
\begin{itemize}
\item[(i)]  Define $E=\range (f)$ and note that the sequence $(x'_k)_{k\in I'_f}$ where $x'_k=x_k|_E$, $k\in I'_f$, is also a $(C,(j_k)_{k\in I'_f})$-RIS. Then,
\[
f(\sum_{k\in I}a_kx_k)=f(\sum_{k\in I'_f}a_kx'_k)
\]
and $\{k\in I'_f:\supp (x'_k)\subset\range (f)\}=I'_f$. Hence, the basic inequality yields $h,g$ and $j_0$ as in items (i) - (iii) such that 
\[
|f(\sum_{k\in I}a_kx_k)|\le C(1+\frac{1}{\sqrt{m_{j_{0}}}})[h+g(\sum_{k\in I'_f}a_ke_{z_k})].
\]
where $z_k=\max\supp (x'_k)$, $k\in I'_f$.
\item[(ii)]	Let $j\in\N$. It follows from the proof of Proposition \ref{basic inequality} that if $f$ has a tree analysis $(f_\alpha)_{\alpha\in\mathcal{A}}$ such that $I_{f_\alpha}=\emptyset$ for every $\alpha\in\mathcal{A}$ with $w(f_\alpha)=m_j$, then the functional $g\in W_{\text{aux}}^{(1)}\cup\{0\}$ that the basic inequality yields for $(x_k)_{k\in I}$ and $f$ has a tree analysis $(g_\beta)_{\beta\in\mathcal{B}}$ with $w(g_\beta)\neq m_j$ for every $\beta\in\mathcal{B}$, whenever $g\neq 0$.
\end{itemize}

\end{rem}

\subsection{Evaluations on Standard Exact Pairs}

We prove the following lemma, which yields Proposition \ref{evaluations SEP} as an immediate corollary. 

\begin{lem}\label{sep = basic + aux}
For every $(C,m_{j_0})$-SEP $(x,f)$ the following hold.
\begin{itemize}
	\item[(i)] For every $f'\in W_{(1)}$ 
	\[
		\big|f'(x)\big|\le 
	\begin{cases}
		\frac{C}{m_{j_0}}(1+\frac{1}{\sqrt{m_{j_0}}}),\quad\quad & f'=\pm e_i^*\text{ for some }i\in\N\\
			C(1+\frac{1}{\sqrt{m_{j_0}}})[\frac{1}{m_{j_0}}+\frac{m_{j_0}}{w(f')}],\quad\quad & w(f')\ge m_{j_0} \\
		 C(1+\frac{1}{\sqrt{m_{j_0}}})[\frac{1}{m_{j_0}}+\frac{2}{w(f')}],\quad\quad & w(f')<m_{j_0}
	\end{cases}
	\]
	\item[(ii)] If $f'\in W_{(1)}$ with a tree analysis $(f'_\alpha)_{\alpha\in\mathcal{A}}$ such that $I_{f'_\alpha}=\emptyset$ for every $\alpha\in\mathcal{A}$ with $w(f'_\alpha)= m_{j_0}$, then \[|f'(x)|\le\frac{3C}{m_{j_0}}(1+\frac{1}{\sqrt{m_{j_0}}}).\]
	\end{itemize}
\end{lem}
\begin{proof}
	Let $(x_k)_{k=1}^n$ be a  $(C,(j_k)_{k=1}^n)$-RIS witnessing that $(x,f)$ is a $(C,m_{j_0})$-SEP. Applying Proposition \ref{basic inequality} , we obtain $h$ and $g$ as in items (i) and (ii) respectively, that satisfy \eqref{basic inequality eq} for $x$ and $f'$, namely,
	\[
	|f'(x)|\le Cm_{j_0}(1+\frac{1}{\sqrt{m_{j_0}}})[h(\tilde{x})+g(\tilde{x})],
	\]
where $\tilde{x}=\sum_{k\in I}a_ke_{z_k}$, $z_k=\max\supp (x_k|_{\range (f')})$ and $I=\{k=1,\ldots,n:\supp (x_k)\cap \range (f')\neq\emptyset\}$. Note that $\tilde{x}$ is a $(n_{j_0},m^{-2}_{j_0})$-b.s.c.c. and hence, since $\supp (h)\in\S_0$, we have $h(\tilde{x})<m_{j_0}^{-2}$.

 To prove (i), first observe that if $g=0$, which is the case e.g. when $f'=\pm e^*_i$ for some $i\in\N$, then we already have established a valid upper bound for $|f'(x)|$. Hence, suppose that $g\neq 0$. Then, using Lemma \ref{evaluations on the basis of auxiliary space} and the fact that $w(g)=w(f')$, we obtain the following upper bounds for $g(\tilde{x})$
 	\[
		g(\tilde{x})\le 
	\begin{cases}
		  \frac{1}{w(f')},\quad\quad & w(f')\ge m_{j_0}\\
		 \frac{2}{w(f')m_{j_0}},\quad\quad & w(f')<m_{j_0}
	\end{cases}
	\]
which yield the desired upper bounds for $|f'(x)|$.
 
 Finally, item (ii) of Remark \ref{remark basic inequality} implies that $g$ admits a tree analysis $(g_\beta)_{\beta\in\mathcal{B}}$ such that $w(g_\beta)\neq m_{j_0}$ for every $\beta\in\mathcal{B}$. We derive the desired upper bound using item (ii) of Lemma \ref{evaluations on the basis of auxiliary space}, which yields that $|g(\tilde{x})|\le 2m^{-2}_{j_0}$.
\end{proof}

\section{Appendix B}

We prove another version of the basic inequality that reduces evaluations on standard exact pairs of $\mathfrak{X}^{(2)}_{\text{awi}}$ to evaluations on the basis of an auxiliary space. The results are almost identical those of Appendix A and we include them for completeness.

\subsection{The Auxiliary Space}

\begin{dfn}
	Let $W_{\text{aux}}^{(2)}$ be the minimal subset of $c_{00}(\N)$ such that
	\begin{itemize}
		\item[(i)] $\pm e^*_i$ is in $W_{\text{aux}}^{(2)}$ for all $i\in\N$ and
		\item[(ii)] whenever $j\in\N$, $(f_i)_{i=1}^d$ is an $\S_{n_j+1}$-admissible sequence in $W_{\text{aux}}^{(2)}$ and $\lambda_1,\ldots,\lambda_d\in\Q$ with $\sum_{i=1}^d\lambda_i^2\le1$, then $f=2m_{j}^{-1}\sum_{i=1}^d\lambda_if_i$ is in $W_{\text{aux}}^{(2)}$.
	\end{itemize}
\end{dfn}

\begin{rem}
For each $f\in W_{\text{aux}}^{(2)}$,   the weight of $f$ is defined as $w(f)=0$ if $f=\pm e^*_i$ for some $i\in\N$ and $w(f)=m_j/2$ in the case where $f=2m_j^{-1}\sum_{i=1}^d\lambda_if_i$.
\end{rem}

The following lemma is a slightly modified version of \cite[Lemma 3.16]{AT}. We use it to prove Lemma \ref{evaluations on the basis of auxiliary space l2}.

\begin{lem}\label{3.16AT l2 version}
Let $f\in W_{\text{aux}}^{(2)}$ with a tree analysis $(f_\alpha)_{\alpha\in\mathcal{A}}$.
\begin{itemize}
\item[(i)] For all $j\in\N$ we have \[\{k\in\supp(f):w_f(e^*_k)<m_{j}\}\in\S_{(\log_2(m_j)-1)(n_{j-1}+1)}.\]
\item[(ii)] If $j\in\N$ is such that $w(f_\alpha)\neq m_j$ for each $\alpha\in\mathcal{A}$ then \[\{k\in\supp(f):w_f(e^*_k)<m_{j}^{2}\}\in\S_{(2\log_2(m_j)-1)(n_{j-1}+1)}.\]
\end{itemize}
\end{lem}
\begin{proof} The proof is similiar to \cite[Lemma 3.16]{AT}.
\end{proof}

Next, we prove a lemma  similar to \ref{evaluations on the basis of auxiliary space}, for the evaluations of functionals in $W^{(2)}_{\text{aux}}$ on the $\ell_2$ version of basic special convex combinations.

\begin{lem}\label{evaluations on the basis of auxiliary space l2}
	Let $j\in\N$ and $\varepsilon>0$ with $\varepsilon\le m_{j}^{-2}$. For every $(2,n_{j},\varepsilon)$-basic s.c.c. $x=\sum_{k\in F}c_ke_k$   the following hold.
	\begin{itemize}
	\item[(i)] For every $f\in W_{\text{aux}}^{(2)}$ 
	\[
		\big|f(x)\big|\le 
	\begin{cases}
		  \frac{1}{w(f)},\quad\quad & w(f)\ge m_{j}/2\\
		 \frac{2}{w(f)m_{j}},\quad\quad & w(f)<m_{j}/2
	\end{cases}
	\]
	\item[(ii)] If $f\in W^{(2)}_{\text{aux}}$ with a tree analysis $(f_\alpha)_{\alpha\in\mathcal{A}}$ such that $w(f_\alpha)\neq m_{j}$ for all $\alpha\in \mathcal{A}$ and $\varepsilon<m_j^{-4}$ then $|f(x)|<2m_j^{-2}$.
	\end{itemize}
\end{lem}
\begin{proof}
Without loss of generality, we may assume that $\supp(f)\subset F$ and $f(e_k)\ge 0$ for every $k\in F$. If $m_{j}/2\le w(f)$, then $\|f\|_{2}\le 1/w(f)$ and hence
\[
\big|f(x)\big|\le \|f\|_{2}\|x\|_2\le \frac{1}{w(f)}.
\]
Suppose now that $m_i<m_{j}$ and let $f=2m^{-1}_{i}\sum_{l=1}^d\lambda_l f_l$, where $(f_l)_{l=1}^d$ is an $\S_{n_i+1}$-admissible sequence in $W^{(2)}_{\text{aux}}$. For $l=1,\ldots,d$, define \[D_l=\{k\in \supp(f_l):w_{f_l}(e^*_k)<m_{j}\},\quad F_l=\supp(f_l)\setminus D_l.\]  Then, Lemma \ref{3.16AT l2 version} (i) implies that $D_l\in \S_{(2\log_2(m_{j})-1)(n_{j-1}+1)}$ for each $l=1,\ldots,d$ and hence, since $(f_l)_{l=1}^d$ is $S_{n_{j-1}+1}$-admissible (recall that $i<j$ since $m_i<m_j$) and $D_l\subset\supp (f_l)$, $l=1,\ldots,d$, we have
\[
D=\cup_{l=1}^dD_l\in S_{n_{j-1}+1}*\S_{(2\log_2(m_{j})-1)(n_{j-1}+1)}= \S_{2\log_2(m_{j})(n_{j-1}+1)}.
\]
Therefore, since $x$ is an $(2,n_j,\varepsilon)$-basic s.c.c. and $2\log_2(m_{j})(n_{j-1}+1)<n_j$, we have $\sum_{k\in D}c^2_k<\varepsilon$. Moreover, observe that for $l=1,\ldots,d $ and $k\in F_l$
\[
f_l(e_k)= \frac{\lambda_{f_l,\alpha_k}}{w_{f_l}(e^*_k)}\le\frac{\lambda_{f_l,\alpha_k}}{m_j}
\]
where $a_k$ is the node in the induced tree analysis of $f_l$ with $f_{l,\alpha_k}=e^*_k$, and
\[
\sum_{l=1}^d\lambda_l^2\sum_{k\in F_l}\lambda_{f_l,\alpha_k}^2\le1.
\]
We then calculate, using the Cauchy-Schwarz inequality
\begin{align*}
	f(x)&=\frac{{2}}{m_i}(\sum_{l=1}^d\lambda_lf_l|_D(\sum_{k\in F}c_ke_k)+\sum_{l=1}^d \lambda_lf_l|_{\N\setminus D}(\sum_{k\in F}c_ke_k))\\
	&=\frac{{2}}{m_i}(\sum_{l=1}^d\lambda_l\sum_{k\in D_l}\frac{c_k\lambda_{f_l,\alpha_k}}{w_{f_l}(e^*_k)}+\sum_{l=1}^d\lambda_l\sum_{k\in F_l}\frac{c_k\lambda_{f_l,\alpha_k}}{w_{f_l}(e^*_k)})\\
	&\le\frac{{2}}{m_i}(\sum_{l=1}^d\lambda_l\sum_{k\in D_l}c_k\lambda_{f_l,\alpha_k}+\frac{1}{m_j}\sum_{l=1}^d\lambda_l\sum_{k\in F_l}c_k\lambda_{f_l,\alpha_k})\\
	&\le \frac{2}{m_i}(\sum_{l=1}^d \lambda_l (\sum_{k\in D_l}c_k^2)^{\frac{1}{2}}(\sum_{k\in D_l}\lambda_{f_l,\alpha_k}^2)^{\frac{1}{2}}+\frac{1}{m_{j}}\sum_{l=1}^d\lambda_l(\sum_{k\in F_l}c_k^2)^{\frac{1}{2}}(\sum_{k\in F_l}\lambda_{f_l,\alpha_k}^2)^{\frac{1}{2}}) \\
	&\le \frac{{2}}{m_i}((\sum_{l=1}^d\lambda^2_l)^{\frac{1}{2}}(\sum_{k\in D}c_k^2)^{\frac{1}{2}}+\frac{1}{m_j}(\sum_{l=1}^d\sum_{k\in F_l}c_k^2)^{\frac{1}{2}}(\sum_{l=1}^d\lambda^2_l\sum_{k\in F_l}\lambda_{f_l,\alpha_k}^2)^{\frac{1}{2}})\\
	&\le \frac{{2}}{m_i}(\sqrt{\e}+\frac{1}{m_j})\le \frac{4}{m_im_j}.	
\end{align*}

Finally, if there is a tree analysis $(f_\alpha)_{\alpha\in\mathcal{A}}$ of $f$ such that $w(f_\alpha)\neq m_j$ for every $\alpha\in\mathcal{A}$, Lemma \ref{3.16AT l2 version} (ii) implies that \[D=\{k\in \supp(f):w_{f}(e^*_k)<m_j^{2}\}\in \S_{(2\log_2(m_j)-1)(n_{j-1}-1)}\] and since $(2\log_2(m_j)-1)(n_{j-1}-1)<n_j$ we have that $\sum_{k\in D}c^2_k<\e$. Hence, using similar arguments as above, we conclude that
\[
f(x)\le \sqrt{\e} + \frac{1}{m^2_j}<\frac{2}{m_j^2}.
\]
\end{proof}

\subsection{The Basic Inequality}

\begin{prop}[Basic Inequality]\label{basic inequality 2}
Let $(x_k)_{k\in I}$ be a $(C,(j_k)_{k\in I})$-RIS in $\mathfrak{X}^{(2)}_{\text{awi}}$ with $4\le \min\supp (x_{\min I})$, $(a_k)_{k\in I}$ be a sequence of nonzero scalars and $f\in W_{(2)}$ with $I_f\neq \emptyset$. Define $t_k=\max\supp (x_k)$, $k\in I$. Then there exist  
\begin{itemize}
\item[(i)] $g\in W_{\text{aux}}^{(2)}\cup\{0\}$ with $w(g)=w(f)/2$ if $g\neq 0$ and $\{k:t_k\in\supp (g)\}\subset I_f$,
\item[(ii)] $h\in\{\text{sign}(a_k) e^*_{t_k}:k\in I_f\}\cup\{0\}$ with $k_0\in I_f$ and $k_0 < \min\supp(g)$ if $h=\text{sign}(a_{k_0})e^*_{t_{k_0}}$ and
\item[(iii)] $j_0\ge \min \{j_k:k\in I_f\}$
\end{itemize}
such that
\begin{equation*}\label{basic inequality 2 eq}
|f(\sum_{k\in I_f}a_kx_k)|\le C(1+\frac{1}{\sqrt{m_{j_{0}}}})[h+g(\sum_{k\in I_f}a_ke_{t_k})].
\end{equation*}
\end{prop}
\begin{proof}
As in Proposition \ref{basic inequality}, we prove the statement by induction on $n=0,1,\ldots$ for every $f\in W^n_{(2)}$ and every RIS. The case of $n=0$ follows easily.

Fix $n\in\N$ and assume that the conclusion holds for every $f\in W^n_{(2)}$ and every RIS. Pick an $f\in W^{n+1}_{(2)}$ with $f={m_i}^{-1}\sum_{l=1}^d\lambda_lf_l$, where $(f_l)_{l=1}^d$ is an $\S_{n_i}$-admissible AWI sequence in $W^n_{(2)}$ and $\lambda_1,\ldots,\lambda_d\in\Q$ with $\sum_{l=1}^d\lambda_l^2\le1$. The proof of the case where $i\ge \max \{j_k:k\in I_f\}$ is identical to that of Proposition \ref{basic inequality}.

Suppose then that $i< \min \{j_k:k\in I_f\}$. The inductive hypothesis implies that, for every $l=1,\ldots,d$ with $I_{f_l}\neq\emptyset$, there are $g_l$, $h_l$ and $j_{0,l}$ as in (i) - (iii) of the statement, that satisfy the conclusion for the functional $f_l$. Define $J_f=\{k\in I_f: f(x_k)\neq 0\}\setminus \cup_{l=1}^dI_{f_l}$. For every $k\in J_f$, since $i<j_k$, Definition \ref{RIS definition} (iii) yields that
\begin{align*}
|f(a_{k}x_{k})|\le (\sum_{l\in L_k}\lambda_l^2)^{\frac{1}{2}}\frac{C}{m_i}|a_{k}|=(\sum_{l\in L_k}\lambda_l^2)^{\frac{1}{2}}\frac{C}{m_i}\text{sign}(a_{k})e^*_{t_{k}}(\sum_{k\in I_f}a_{k}e_{t_{k}}).
\end{align*}
where
\[
L_k=\{l\in\{1,\ldots,d\}:\supp(x_k)\cap\supp(f_l)\neq\emptyset\}.
\]
Hence we calculate
\begin{align*}
&|f(\sum_{k\in I_f}a_kx_k)|\le|f(\sum_{k\in J_f}a_kx_k)|+ |f(\sum_{k\in \cup_{l=1}^dI_{f_l}}a_kx_k)|\\
&\le\frac{C}{m_i}\sum_{k\in J_f}(\sum_{l\in L_k}\lambda_l^2)^{\frac{1}{2}}\text{sign}(a_{k})e^*_{t_{k}}(\sum_{k\in I_f}a_{k}e_{t_{k}})+\frac{C}{m_i}\sum_{l=1}^d[(1+\frac{1}{\sqrt{m_{j_{0,l}}}})\lambda_l(h_l+g_l)](\sum_{k\in I_{f_l}}a_ke_{t_k})\\
&\le C(1+\frac{1}{\sqrt{m_{j_{\min I_f}}}})[\frac{1}{m_i}(\sum_{k\in J_f}(\sum_{l\in L_k}\lambda_l^2)^{\frac{1}{2}}\text{sign}(a_{k})e^*_{t_{k}}+\sum_{l=1}^d\lambda_lh_l+\lambda_lg_l)](\sum_{k\in I_f}a_ke_{t_{k}}).
\end{align*}
Define 
\[
g=\frac{2}{m_i}(\sum_{k\in J_f}\frac{1}{2}(\sum_{l\in L_k}\lambda_l^2)^{\frac{1}{2}}\text{sign}(a_{k})e^*_{t_{k}}+\sum_{l=1}^d\frac{\lambda_l}{2}h_l+\frac{\lambda_l}{2}g_l).
\]
Then, observe that each $l=1,\ldots,d$, belongs to $L_k$ for at most two $k\in J_f$ and thus, using the same arguments as in Proposition \ref{basic inequality}, we have that $g\in W_{\text{aux}}^{(2)}$ and this completes the proof for case where $i<j_k$ for all $k\in I_f$.

Finally, the proof of the remaining case is the same as in Proposition \ref{basic inequality}.
\end{proof}

\subsection{Evaluations on Standard Exact Pairs} Finally, we prove the following lemma which shows that standard exact pairs are in fact strong exact pairs.

\begin{lem}
For every $(2,C,m_{j_0})$-SEP $(x,f)$ the following hold.
\begin{itemize}
	\item[(i)] For every $g\in W$ 
	\[
		\big|g(x)\big|\le 
	\begin{cases}
			2C(1+\frac{1}{\sqrt{m_{j_0}}})[\frac{1}{m_{j_0}}+\frac{m_{j_0}}{w(g)}],\quad\quad & w(g)\ge m_{j_0} \\
		 2C(1+\frac{1}{\sqrt{m_{j_0}}})[\frac{1}{m_{j_0}}+\frac{2}{w(g)}],\quad\quad & w(g)<m_{j_0}
	\end{cases}
	\]
	\item[(ii)] If $g\in W$ with a tree analysis $(g_\alpha)_{\alpha\in\mathcal{A}}$ such that $I_{g_\alpha}=\emptyset$ for every $\alpha\in\mathcal{A}$ with $w(g_\alpha)= m_{j_0}$, then \[|g(x)|\le\frac{3C}{m_{j_0}}(1+\frac{1}{\sqrt{m_{j_0}}}).\]
	\end{itemize}
\end{lem}
\begin{proof}
Apply the basic inequality and the evaluations of functionals in $W_{\text{aux}}^{(2)}$ on $2$-b.s.c.c. from Lemma \ref{evaluations on the basis of auxiliary space l2}. The proof is identical to that of Lemma \ref{sep = basic + aux}.
\end{proof}


\begin{thebibliography}{99999}

\bibitem{AK}  F. Albiac and N.J. Kalton, {\sl Topics in Banach space theory}, Graduate Texts in Mathematics 233, Springer, New York, 2006.
\bibitem{AA} D. Alspach and S.A. Argyros, {\sl Complexity of weakly null
    sequences}, Diss. Math. 321 (1992), 1--44.
 \bibitem{AAT} S. A. Argyros, A. Arvanitakis and A. Tolias, Saturated extensions, the attractors method and Hereditarily James tree spaces, in Methods in Banach Space Theory (J. M. F. Castillo and W. B. Johnson, eds.), London Mathematical Society Lecture Note Series 337, Cambridge University Press, 2006, pp. 1--90.
 \bibitem{AD} S. A. Argyros and I. Deliyanni, Banach Spaces of the Type of Tsirelson arXiv:9207206 (1992).
\bibitem{AD2} S. A. Argyros and I. Deliyanni,
Examples of asymptotic $\ell_{1}$ Banach spaces,
Trans. Am. Math. Soc. 349, No. 3, 973-995 (1997). 
\bibitem{AGLM} S. A. Argyros, A. Georgiou, A.-R. Lagos and P. Motakis, 
{\sl Joint spreading models and uniform approximation of bounded
operators}, Studia Math. 253 (2020), 57--107.
\bibitem{AGMM} S.A. Argyros, A. Georgiou, A. Manoussakis and P. Motakis, Applying the attractor's method for the complete separation of asymptotic properties (in preparation).
\bibitem{AGM} S.A. Argyros, A. Georgiou, P. Motakis, Non-asymptotic $\ell_1$ spaces with unique $\ell_1$ asymptotic model, Bulletin of the Hellenic Mathematical Society 64 (2020), 32--55
\bibitem{AM} S.A. Argyros, P. Motakis,  On the complete
separation of asymtpotic structures in Banach spaces
, Advances    in Mathematics  362 (2020),
https://doi.org/10.1016/j.aim.2019.106962
\bibitem{AMP} S.A. Argyros, A. Manoussakis, A. Pelczar-Barwacz,
  On the hereditary proximity to $\ell_{1}$,  J. Funct. Anal. 261 (2011), no. 5
\bibitem{AMT} S. A. Argyros, S. Mercourakis and A. Tsarpalias, {\sl Convex
unconditionality and summability of weakly null sequences}, Israel J. Math. 107
(1998), 157--193.
\bibitem{AOST} G. Androulakis, E. Odell, Th. Schlumprecht, N. Tomczak-Jaegermann, On the structure of the spreading models of a Banach space, Canad. J. Math.57(2005), no.4, 673--707.
 \bibitem{AT} S.A. Argyros and A. Tolias, {\sl Methods in the theory of
 hereditarily indecomposable Banach spaces},
Mem. Amer. Math. Soc. 170 (2004), no. 806.
\bibitem{BLMS} F. Baudier, G. Lancien, P. Motakis, and Th. Schlumprecht, A new coarsely rigid class of Banach spaces, arXiv:1806.00702v1 (2018).
\bibitem{BS}  A. Brunel, L. Sucheston, {\sl On B-convex Banach spaces}, Math. Systems Theory 7 (1974), no. 4,
294-299
\bibitem{DM} I. Deliyanni and A. Manoussakis, Asymptotic $\ell_{p}$
  Hereditarily  Indecomposable Banach spaces, Illinois Journal of
  Mathematics (2007) 767-803
\bibitem{FJ} T. Figiel, W.B Johnson,A uniformly convex Banach space
  which contains no $\ell_{p}$, Compositio Mathematica, Volume 29 (1974) no. 2, p. 179-190 
  \bibitem{FOSZ}
D. Freeman, E. Odell, B. Sari, B. Zheng,
{\em On spreading sequences and asymptotic structures},
Trans. AMS, 370, no. 10, (2018) 6933-6953.
  \bibitem{GM} W. T. Gowers and B. Maurey, {\sl The unconditional basic sequence problem}, Journal of AMS, 6 (1993), 851--874.
\bibitem{HOdell} L. Halbeisen,E.Odell, On asymptotic models in Banach spaces.. 
Isr. J. Math. 139, 253-291 (2004)
\bibitem{Jrefl} R.C. James, {\em Bases and reflexivity of Banach spaces}, Ann. of Math. (2) 52 (1950), 518-527
\bibitem{James distortion} R. C. James, {\sl Uniformly non-square Banach spaces}, Annals of Mathematics (1964), Vol. 80, No. 3,
542--550.
 \bibitem{LT} J. Lindenstrauss and L. Tzafriri,
 {\sl Classical Banach spaces I}, Springer-Verlag 92, 1977.
 \bibitem{MMT} B. Maurey, V.D. Milman and N. Tomczak-Jaegermann, {\sl Asymptotic
 Infinite-Dimensional theory of Banach spaces}, Geometric aspects of functional analysis (Israel, 1992--1994),  149--175, Oper. Theory Adv. Appl., 77, Birkh\"auser, Basel, 1995.
 \bibitem{MT} V. Milman and N. Tomczak-Jaegermann, {\sl Asymptotic $\ell_p$
 spaces and bounded distortion}, Banach spaces (Merida
 1992), 173--195, Comtemp. Math. 144, Amer. Math. Soc., Providence, RI,
 1993.




 
\bibitem{S} T. Schlumprecht, An arbitrarily distortable Banach space, Israel J. Math., 76 (1991), 81--95
 \bibitem{T}B.S. Tsirelson, Not every Banach space contains an imbedding of $\ell_{p}$ or $c_{0}$, Funct. Anal. Appl.,8:2(1974),pp. 138-141 
\end{thebibliography}
\end{document}